\documentclass{article}
\usepackage{fullpage}
\usepackage[utf8]{inputenc}
\usepackage{hyperref}
\usepackage{natbib}
\setcitestyle{authoryear, open={((},close={)}}
\usepackage{algorithm}
\usepackage{algorithmic}
\usepackage{comment}
\usepackage{fullpage}
\usepackage{amsmath,amssymb,amsthm}
\usepackage{enumerate}
\usepackage{appendix}
\usepackage{graphicx,subfig}
\usepackage[dvipsnames]{xcolor}
\usepackage{color}
\usepackage{multicol}
\usepackage{hyperref}
\hypersetup{
        colorlinks   = true,
        linkcolor    = purple,
        citecolor    = blue
}
\usepackage{authblk}
\newtheorem{theorem}{Theorem}

\newtheorem{lemma}{Lemma}

\theoremstyle{definition}
\newtheorem{definition}{Definition}
\theoremstyle{remark}
\newtheorem{remark}{Remark}
\theoremstyle{example}
\newtheorem{example}{Example}
\usepackage{bm} 
\newcommand\numberthis{\addtocounter{equation}{1}\tag{\theequation}}
 
\DeclareMathOperator*{\argmin}{arg\,min}
\renewcommand{\t}{^{\top}}
\newcommand{\R}{\mathbb{R}}
\newcommand{\ipq}{{I}_{p,q}}
\newcommand{\inv}{^{-1}}

\newcommand{\eps}{\ensuremath{\varepsilon}}
\newcommand{\p}{\mathbb{P}}
\newcommand{\E}{\mathbb{E}}
\newcommand{\F}{\text{F}}
\renewcommand{\hat}{\widehat}
\renewcommand{\tilde}{\widetilde}
\newcommand{\argmax}{\mathop{\rm arg\max}}
\makeatletter
\DeclareRobustCommand\widecheck[1]{{\mathpalette\@widecheck{#1}}}
\def\@widecheck#1#2{%
    \setbox\z@\hbox{\m@th$#1#2$}%
    \setbox\tw@\hbox{\m@th$#1%
       \widehat{%
          \vrule\@width\z@\@height\ht\z@
          \vrule\@height\z@\@width\wd\z@}$}%
    \dp\tw@-\ht\z@
    \@tempdima\ht\z@ \advance\@tempdima2\ht\tw@ \divide\@tempdima\thr@@
    \setbox\tw@\hbox{%
       \raise\@tempdima\hbox{\scalebox{1}[-1]{\lower\@tempdima\box
\tw@}}}%
    {\ooalign{\box\tw@ \cr \box\z@}}}
\makeatother
\renewcommand{\check}{\widecheck}
\newcommand{\wh}[1]{\widehat{#1}}

\newcommand{\dto}{\quad\xrightarrow{\operatorname{d}} \quad }

\renewcommand{\mathbf}{}
\renewcommand{\bm}{}

\usepackage[capitalize]{cleveref}
\usepackage{mathtools}
\usepackage{stmaryrd}

\title{An Overview of Asymptotic Normality in Stochastic Blockmodels:~Cluster Analysis and Inference}

\author{
    Joshua Agterberg\thanks{Department of Statistics, University of Illinois Urbana-Champaign}
    ~ and ~
    Joshua Cape\thanks{Departments of Statistics, University of Wisconsin--Madison}
}

\date{\today}

\begin{document}

\maketitle

\begin{abstract}
    This paper provides a selective review of the statistical network analysis literature focused on clustering and inference problems for stochastic blockmodels and their variants. We survey asymptotic normality results for stochastic blockmodels as a means of thematically linking classical statistical concepts to contemporary research in network data analysis. Multiple different forms of asymptotically Gaussian behavior arise in stochastic blockmodels and are useful for different purposes, pertaining to estimation and testing, the characterization of cluster structure in community detection, and understanding latent space geometry. This paper concludes with a discussion of open problems and ongoing research activities addressing asymptotic normality and its implications for statistical network modeling.
\end{abstract}

\tableofcontents

\section{Introduction}\label{sec:introduction}
Stochastic blockmodels (SBMs) \citep{holland_stochastic_1983,wasserman_social_1994} are generative models for random graphs in which network connectivity is determined by community (block) structure among nodes. In the simplest setting, nodes (vertices) belong to communities, and connectivity, as quantified via the presence or absence of links (edges) between nodes, is probabilistically governed solely by nodes' community memberships. Stochastic blockmodels and their variants serve as workhorse generative models in the statistical analysis of networks, particularly in the vertex-centric, networks-as-graphs paradigm.

Stochastic blockmodels, instances of which are sometimes referred to as \emph{planted partition models} \citep{schaeffer_graph_2007,lancichinetti_community_2009}, have received significant attention within the mathematics, statistics, computer science, and physics research communities, finding applications in social network analysis, neuroimaging, and elsewhere. Popular variants and extensions include mixed membership \citep{airoldi_mixed_2008} and degree-corrected stochastic blockmodels \citep{karrer_stochastic_2011}, among numerous others \citep{arroyo_overlapping_2022,noroozi_estimation_2021,noroozi_sparse_2021,sengupta_block_2018,zhang_detecting_2020}. Historical developments pertaining to the study of networks and SBMs are cataloged in, for example, \citet{fortunato_community_2010,fortunato_community_2016,lee_review_2019,wasserman_social_1994}. Stochastic blockmodels have themselves already been the subject of survey-length treatment, notably in \citet{abbe_community_2018} which summarizes phase transition phenomena for community recovery. Stochastic blockmodels also arise in survey-length treatments of research focused on dynamic networks \citep{kim_review_2018} as well as latent space modeling and analysis \citep{athreya_statistical_2018,matias_modeling_2014}, among others \citep{gao_minimax_2021,goldenberg_survey_2009,zhao_survey_2017}.

The study of networks at times challenges convention and intuition developed in classical statistical inference settings. Typically, one observes a single or multiple networks, each describing \emph{pairwise interactions} between vertices, in contrast to observing independent Euclidean-valued feature vectors for each sample unit. When modeling network data, several natural questions arise. Do multiple networks (graphs) constitute multiple observations, or are individual vertices within a single network the units of interest? Correspondingly, is the inference task focused on vertex-specific behavior by, for example, seeking to infer individual edge probability parameters, or does one wish to understand the behavior of certain summary statistics computed from the entire network? Further, where or how could asymptotic normality appear, and to what end? Possible answers depend on how statistical units of interest are conceptualized. This survey contributes to the existing literature by providing a unified treatment of asymptotic normality results for stochastic blockmodels, as a means of contextualizing and detailing advances in \emph{statistical network analysis}.

A sizable body of literature considers spectral-based matrix factorization approaches to network analysis by leveraging the convenient representation of networks in terms of their adjacency or Laplacian matrix. Such approaches have led to widespread use of \emph{spectral embeddings}, namely eigendecomposition-based low-dimensional Euclidean representations of vertices in graphs. This survey discusses existing asymptotic normality results for different matrix factorization approaches to network analysis as well as applications of these results to parameter estimation and hypothesis testing. In part, to the extent possible, this survey seeks to elucidate the interplay between underlying matrix structure and asymptotic normality behavior. Loosely speaking, a tug-of-war exists between network sparsity, the antagonist, and network size, the protagonist. The emergence of asymptotic normality depends on the extent of community separation (between-community heterogeneity), sample size, and sparsity, with rate of convergence, covariance matrix, and centering determined by these aspects.

Numerous existing works in the literature address statistical consistency under SBMs and in more general network settings, without addressing asymptotic distributional properties. To maintain focus and a reasonable page count, this survey restricts its attention to only the most closely related works. In particular, this review does not delve into the rich literature on information-theoretical thresholds and phase transition phenomena for random graphs \citep{abbe_community_2018} nor does it consider graphons as graph limit objects \citep{lovasz_large_2012}. Many of the results covered herein hold for more general models, and while we aim to mention this when appropriate, we emphasize that the focus herein is deliberately on stochastic blockmodels.

\subsection{Notation}
Given any positive integer $n$, define $[n] \coloneqq \{1,\dots,n\}$. We write $X \sim F$ to denote that the random variable $X$ is distributed according to the distribution $F$, and for a symmetric matrix $A$ with independent random Bernoulli entries on its upper triangle, we write $A \sim P$ to denote that $\mathbb{E}(A) = P$.  We use lowercase letters to denote vectors, writing the $d$-dimensional vector $x$ coordinate-wise as $x=(x_{1},\dots,x_{d})\t$. We use $1$ both as a scalar and to denote the vector of all ones when its dimension is clear from context. We write $\mathcal{S}_{K-1}$ to denote the interior of the unit simplex in $K$ dimensions, i.e.,~$x \in \mathcal{S}_{K-1}$ if and only if simultaneously $x\t 1 = 1$ and entrywise $x_{i} > 0$ for all $1 \le i \le d$. We use upper-case letters to denote matrices, say $M$, where $M_{i *}$ and $M_{* j}$ denote the $i$-th row and $j$-th column of $M$, respectively. We use $\|M\| \coloneqq \max_{i} \sigma_{i}(M)$ for the largest singular value of $M$ and $\|M\|_{\F} \coloneqq \sqrt{\sum_{i} \sigma_{i}^{2}(M)}$ for the Frobenius matrix norm. We use uppercase letters to denote random vectors, with single indices indicating individual observables.

We use the asymptotic notation $f(n) = O(g(n))$ when there exists a constant $C > 0$ such that $f(n) \leq C g(n)$ holds for all $n$ sufficiently large. Similarly, $f(n) = \Theta(g(n))$ indicates that simultaneously $f(n) = O(g(n))$ and $g(n) = O(f(n))$. Herein, $f(n) = \omega(g(n))$ indicates $g(n)/f(n) \to 0$ as $n \rightarrow \infty$, whereas $f(n) = o(g(n))$ indicates $f(n)/g(n) \rightarrow 0$ as $n \rightarrow \infty$. In addition, we write $f(n) \gtrsim g(n)$ if $g(n) = O(f(n))$ and $f(n) \asymp g(n)$ if $f(n) = \Theta(g(n))$, and $f(n) \ll g(n)$ if $f(n)/g(n) = o(1)$. Finally, we write $f(n) = \tilde O(g(n))$ when there exists a constant $c > 0$ not depending on $n$ such that $f(n) = O(\log^c(n) \cdot g(n))$.

\section{Stochastic Blockmodels}\label{sec:sbms}
This section briefly reviews the basics of stochastic blockmodeling wherein graphs (networks) constitute the observable data. Arising at the confluence of \emph{blockmodeling} and \emph{stochastic modeling}, SBMs date back several decades to at least \cite{holland_stochastic_1983,wasserman_social_1994}. Significant attention has been devoted to SBMs as concrete instances of more general \emph{latent space network models} \citep{hoff_latent_2002}.

\begin{definition}[Stochastic blockmodel --- undirected, hollow]
    \label{def:SBM_hollow}
    Let $B \in [0,1]^{K \times K}$ be a symmetric matrix of connectivity probabilities, let $\rho = \rho_{n} \in (0,1]$, and let $\pi = (\pi_{1}, \dots, \pi_{K})\t \in \mathcal{S}_{K-1}$. We say that $A$ is the adjacency matrix for a stochastic blockmodel graph on $n$ vertices with sparsity term $\rho$ if the following hold. First, the block assignment map $\tau: [n] \to [K]$ satisfies $\p[\tau(i) = k] = \pi_{k}$ for $1 \le k \le K$ and for each $1 \le i \le n$. Second, conditional on $\tau$, the upper triangular entries $A_{ij}$ are independently generated in the manner
    \begin{align*}
        A_{ij}
        =
        \begin{cases}
            \operatorname{Bernoulli}(\rho B_{\tau(i),\tau(j)})
            &
            \textnormal{~for~}i < j, \\
            0
            &
            \textnormal{~for~}i = j, \\
            A_{ji}
            &
            \textnormal{~for~}i > j.
        \end{cases}
    \end{align*}
\end{definition}
In words, vertex $i$ belongs to community $l$ with probability $\pi_{l}$, and conditional on the community memberships, the probability of an edge occurring between vertices $i$ and $j$ is determined solely by the community memberships, namely the $(\tau(i),\tau(j))$ entry of $B$. Throughout this work, we implicitly assume that $B$ has distinct rows which prevents multiple communities from exhibiting identical, indistinguishable connectivity behavior, thereby mitigating a potential source of non-identifiability.

\begin{remark}[Presence or absence of self-edges]
    In \cref{def:SBM_hollow}, self-edges (loops) are not permitted (i.e.,~hollowness is enforced), since $A_{ii} = 0$ for all $1 \le i \le n$. If instead self-edges are permitted (i.e.,~hollowness is not enforced), then the main diagonal entries $A_{ii}$ are typically taken to be independent random variables satisfying $A_{ii} \sim \operatorname{Bernoulli}(\rho B_{\tau(i),\tau(i)})$ for all $1 \le i \le n$.
\end{remark}

\begin{remark}[Fixed community memberships and the expected adjacency matrix] \label{rem:fixed_memberships_&_expectation}
   When the community memberships are taken to be fixed (nonrandom), the matrix of edge probabilities for generating an adjacency matrix $A$, denoted by $P \in [0,1]^{n \times n}$, can be expressed as $P \equiv \rho_{n} Z B Z\t$ provided the matrix $Z \in \{0,1\}^{n \times K}$ satisfies $Z_{ij} = 1$ if and only if $\tau(i) = j$. Depending on whether or not the graph corresponding to $A$ is permitted to have self-edges, per \cref{def:SBM_hollow}, either $\E[A] = P - \operatorname{diag}(P)$ or $\E[A] = P$. In certain instances, the distinction between allowing or forbidding the presence of self-edges has a negligible impact on statistical performance guarantees and procedures, while at other times the distinction cannot be ignored. The latter point is addressed further in subsequent sections when needed.

   Irrespective of whether or not self-edges are permitted, $A$ can be viewed as a perturbation of the (approximately) low-rank matrix $P$. Here, low-rankness plays a crucial role in the subsequent analysis of spectral embeddings of $A$.
   
\end{remark}

\subsection{Generalizations of Blockmodels}
The stochastic blockmodel per \cref{def:SBM_hollow} permits heterogeneous nodal connectivity across different blocks (communities) but simultaneously imposes stochastic equivalence on all nodes belonging to the same block. In other words, vertices in a SBM graph can be viewed as each possessing a latent attribute that governs edge formation. A more flexible and arguably realistic modeling approach would allow vertices to possess shared characteristics, to various extents, from among a common set of possible characteristics. Enter the \emph{mixed-membership stochastic blockmodel} (MMSBM) \citep{airoldi_mixed_2008}.

\begin{definition}[Mixed-membership stochastic blockmodel --- undirected, loopy]
    We say that $A$ is the adjacency matrix for a mixed-membership stochastic blockmodel graph on $n$ vertices with sparsity term $\rho = \rho_{n} \in (0,1]$ when the upper triangular entries $A_{ij}$ are independently generated in the manner
    \begin{align*}
        A_{ij} = 
        \begin{cases}
            \operatorname{Bernoulli}(\rho \eta_{i}\t B \eta_{j})
            &
            \textnormal{~for~}i \le j, \\
            A_{ji}
            &
            \textnormal{~for~}i > j,
        \end{cases}
    \end{align*}
    where $B \in [0,1]^{K \times K}$ is a symmetric matrix of connectivity probabilities and $\{\eta_{i}\}_{i\in[n]}$ is a collection of non-negative mixed membership vectors $\eta_{i}\in\R^{K}$ satisfying $\|\eta_{i}\|_{1}=1$. The unobserved vectors $\eta_{i}$ could be specified as deterministic or could be modeled in a generative fashion, for example, i.i.d.~$\eta_{i} \sim \operatorname{Dirichlet}(\alpha)$ for some concentration parameter vector $\alpha \in \R_{+}^{K}$.
\end{definition}

The MMSBM requires that the membership vectors $\eta_{i}$ belong to the unit simplex and
reduces to the SBM when all nodes in the graph are \emph{pure nodes}, namely when their mixed membership vectors are all standard basis vectors (i.e., the vertices of the unit simplex). A different approach to generalizing SBMs is to permit further node-specific connectivity properties, modeled by permitting heterogeneous vertex degree parameters interpretable as a measure of importance in a network, at the expense of introducing additional parameters and model complexity. Thus arises the \emph{degree-corrected stochastic blockmodel} (DCSBM) \citep{dasgupta_spectral_2004,karrer_stochastic_2011}.

\begin{definition}[Degree-corrected stochastic blockmodel --- undirected, hollow]
    We say that $A$ is the adjacency matrix for a degree-corrected stochastic blockmodel graph on $n$ vertices with sparsity term $\rho = \rho_{n} \in (0,1]$ when the upper triangular entries $A_{ij}$ are independently generated in the manner
    \begin{align*}
        A_{ij}
        =
        \begin{cases}
            \operatorname{Bernoulli}(\rho \theta_{i} \theta_{j} B_{\tau(i), \tau(j)})
            &
            \textnormal{~for~}i < j, \\
            0
            &
            \textnormal{~for~}i = j, \\
            A_{ji}
            &
            \textnormal{~for~}i > j,
        \end{cases}
    \end{align*}
    where $B \in [0,1]^{K \times K}$ is a symmetric matrix of connectivity probabilities, $\theta = (\theta_{1},\dots,\theta_{n})\t$ is a vector of node-specific nonnegative degree parameters, and $(\tau(1), \dots, \tau(n))\t$ is the vector of node block labels. The unobserved degree parameters $\theta_{1},\dots,\theta_{n}$ could be specified as deterministic or could alternatively be modeled in a generative fashion, such as i.i.d.~$\theta_{i}\sim\operatorname{Uniform}[a,b]$.
\end{definition}

Further generalizations beyond those discussed above include but are not limited to the \emph{popularity-adjusted stochastic blockmodel} and the \emph{degree-corrected mixed-membership blockmodel} \citep{ jin_fast_2015,karrer_stochastic_2011,sengupta_block_2018,zhang_detecting_2020} as well as \emph{bipartite stochastic blockmodels} and \emph{directed stochastic blockmodels}. For a more in-depth account of these different models, see for example \citet{noroozi_hierarchy_2022}.

The aforementioned models can be viewed as particular instances of the \emph{generalized random dot product graph} \citet{rubin-delanchy_statistical_2022} if one imposes a distribution on the underlying parameters (such as the degree corrections or mixed memberships). In all cases, $P$ is (approximately) low-rank, which facilitates the analysis of spectral algorithms. \cref{def:grdpg} provides further discussion and details.

\section{Community Detection in Stochastic Blockmodels}\label{sec:community_detection}
In network analysis, the problem of \emph{community detection} is to cluster vertices or nodes into \emph{communities} according to their connectivity properties. For stochastic blockmodels, in which nodes belong to latent ground-truth communities or blocks, link connection probabilities are determined solely by their community memberships. Traditionally, edges are observed in stochastic blockmodel graphs, whereas nodes' community memberships and their connectivity probabilities are unobserved, hence the impetus for community detection. As previously remarked, there are at least two competing aspects at play. One aspect is the number of vertices, a measure of sample size, whereby having more vertices leads to observing more pairwise edges which improves inference and community recovery, while the second aspect is the level of global network sparsity, determining how many edges are present in the graph, whereby the presence of few edges amounts to sparsely observed data, making it difficult or even impossible to recover communities, say, better than chance. Various approaches exist for addressing the community detection problem in networks, including but not limited to maximum likelihood estimation, modularity maximization, and spectral clustering. Beyond the present paper, a discussion of these and other methods can be found in the existing literature. See \citet{zhao_survey_2017} for a theoretically focused survey, and see \citet{fortunato_community_2010} for a discussion of practical considerations.

Matrix factorization approaches to community detection are frequently framed as two-step procedures to obtain estimated community memberships. First, starting with an $n \times n$ (similarity) matrix $S$, such as an observed adjacency matrix, obtain an $n \times d$ spectral embedding of $S$ by taking the $d \ll n$ leading eigenvectors, perhaps scaled by their eigenvalues, and second, cluster the $n$ rows of this spectral embedding using a clustering method, such as $K$-means. \cref{alg:spectral_clustering} provides a representative example of spectral clustering algorithms but is by no means exhaustive. For a tutorial on spectral clustering and extensive additional discussion, see \citet{von_luxburg_tutorial_2007}.

The term \emph{spectral embedding} refers to obtaining a low-dimensional representation of the matrix $S$, often through its \emph{spectral decomposition}, namely its eigenvalues and eigenvectors. We will survey different procedures for obtaining spectral embeddings and their statistical properties. Notably, different spectral embeddings lead to different-yet-related forms of asymptotic normality as quantified through the centering, scaling, and associated asymptotic covariance matrices. \cref{sec:populations} discusses the properties of spectral embeddings obtained from the population-level probability matrix with fixed community memberships, and \cref{sec:subspace_perturbation} addresses eigenvector estimation and spectral embeddings obtained from observed adjacency matrices, focusing on consistent estimation of population leading eigenvectors. Discussion of consistent eigenvector estimation and clustering serves as a stepping stone for the presentation of asymptotic normality results in \cref{sec:asymp_norm}. As illustrated throughout this survey, spectral embeddings can be used for other purposes beyond estimating community memberships, such as parameter estimation (\cref{sec:parameter}) or hypothesis testing (\cref{sec:hyp_testing}).

\begin{figure}[h]
    \centering\includegraphics[width=.5\textwidth,keepaspectratio]{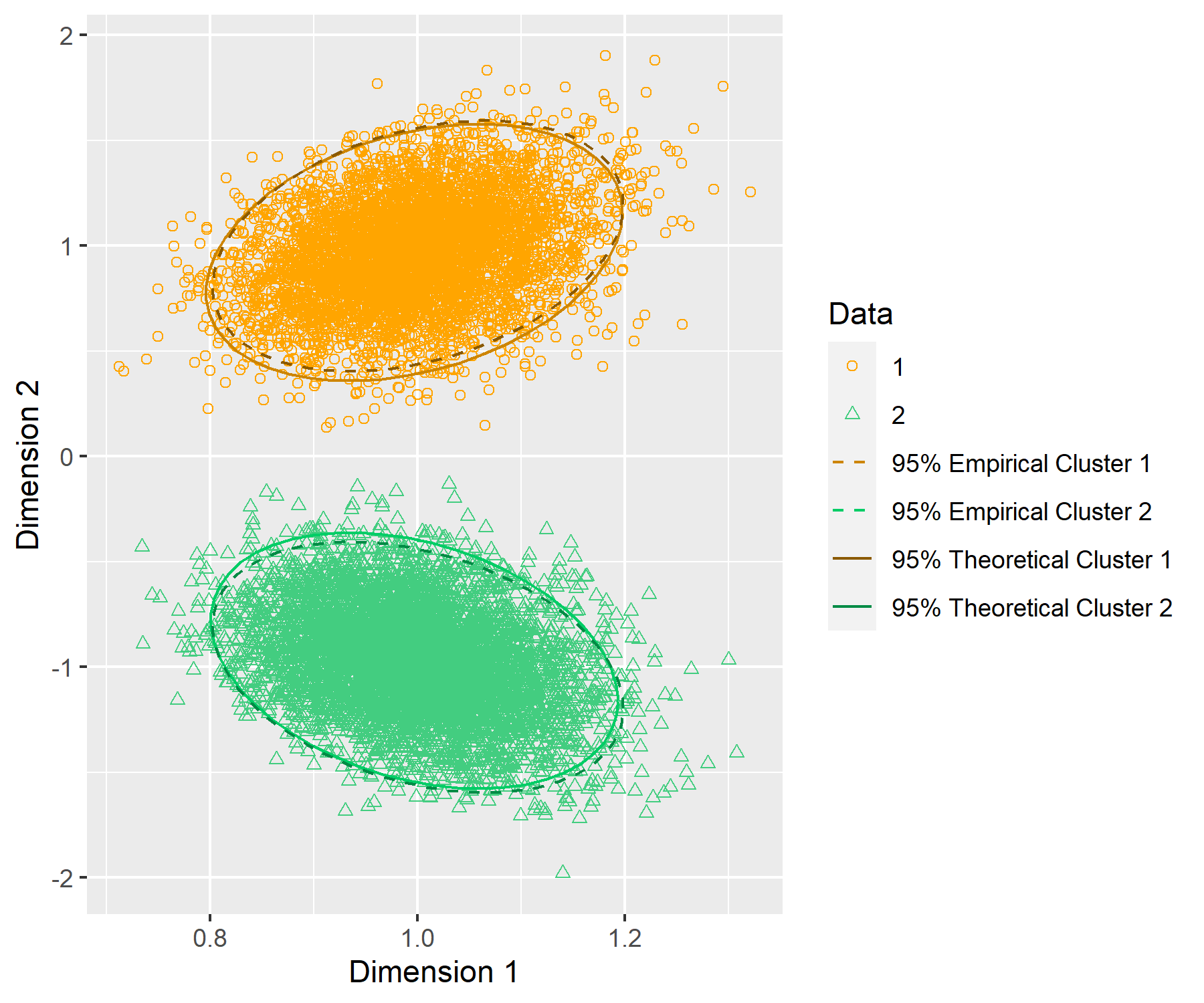}
     \caption{Scatter plot for the adjacency matrix eigenvector components of a SBM with two communities and $n = 10^{4}$ nodes. See \cref{sec:community_detection} for details.} 
     \label{fig:uhat_whet_appetite}
\end{figure}

As a preview, consider the specific matrix
\begin{align*}
    B
    =
    \begin{pmatrix}
        .02 & .01 \\
        .01 & .02
    \end{pmatrix},
\end{align*}
and consider graphs of size $n = 10^{4}$ with $n/2$ vertices in each community, also referred to as a specific \emph{balanced homogeneous planted partition} model. \cref{fig:uhat_whet_appetite} plots the rows of the $n \times 2$ matrix whose columns are eigenvectors for the leading two eigenvalues of a single observed adjacency matrix $A$. The dotted lines depict $95\%$ empirical confidence ellipses for each community assuming knowledge of the communities, and the solid lines depict the $95\%$ theoretical confidence ellipses for each community using the per-block node-specific asymptotic results in \cref{sec:asymp_norm}. The figure illustrates that the rows of the empirical eigenvectors are centered around one of two community-specific centroids. \cref{sec:populations} shows that at the population level there are in fact distinct, separated community-specific centroids, and \cref{sec:subspace_perturbation} shows that the empirical eigenvectors are, in a well-defined sense, concentrating around the true eigenvectors. \cref{fig:uhat_whet_appetite} further suggests to the possibility of approximate asymptotic normality for each vertex embedding which is discussed in greater detail in \cref{sec:asymp_norm}.

\begin{algorithm}[t]
\begin{algorithmic}
    \REQUIRE An $n \times n$ symmetric similarity matrix $S$, an embedding dimension $d$, and the number of communities $K$ (optional).
    \STATE DO: Embed $S$ into $\R^{d}$ by computing the $n \times d$ matrix $\widehat{U}$ whose orthonormal columns are eigenvectors corresponding to the $d$ largest-in-magnitude eigenvalues, arranged in the diagonal matrix $\widehat{\Lambda}$.
    \STATE DO: Cluster the rows of $\widehat{U}$, of $\widehat{U}|\widehat{\Lambda}|^{1/2}$, or of $\widehat{U}\widehat{\Lambda}$, e.g.,~using $K$-means, $K$-medoids, or an expectation maximization-based mixture modeling approach. Alternatively, $K$ could be selected by an information criterion.
    \RETURN A vector of community memberships $\widehat{\tau} \in \{1,\dots,K\}^{n}$ or a membership matrix $\widehat{Z}\in \{0,1\}^{n \times K}$.
\end{algorithmic}
    \caption{A representative spectral clustering algorithm}
\label{alg:spectral_clustering}
\end{algorithm}

We remark that certain properties of spectral embeddings carry over to the model generalizations discussed in \cref{sec:sbms}, at times with minor or major modifications. For example, in the degree-corrected stochastic blockmodel, authors have proposed various row-normalization techniques to mitigate the influence of (nuisance) degree parameters when estimating community membership \citep{du_hypothesis_2023,fan_simple_2022,jin_fast_2015,jin_mixed_2023}. In the mixed-membership stochastic blockmodel, simplices \citep{mao_estimating_2020} and simplicial geometric structures \citep{jin_mixed_2023,mao_overlapping_2018} arise in latent space which can subsequently be used for statistical inference and modeling. The above references elaborate on the relationships between spectral embeddings and population properties in more general models.

\subsection{Population Properties of Spectral Embeddings}\label{sec:populations}
For succinctness and pertinence, this section includes only the adjacency and Laplacian input variants to \cref{alg:spectral_clustering} for which asymptotic normality is discussed in \cref{sec:asymp_norm}. Here and in the remainder of \cref{sec:community_detection}, we consider community memberships to be fixed \emph{a priori}, with $K$ communities and each community $1 \le k \le K$ having size $n_{k}$. In other words, we condition on the community memberships and associated community sizes.

Consider the population adjacency matrix, that is, the matrix $\E[A] = \rho_{n} Z B Z\t = P$ whose independent upper triangular entries satisfy $\p[A_{ij} = 1] = P_{ij}$ for all $i \leq j$. Self-edges are permitted for simplicity, and the surrounding discussion is not materially different for this section if self-edges are disallowed. Recall that $Z$ is the $n \times K$ matrix whose entries satisfy $Z_{ik} = 1$ when vertex $i$ belongs to community $k$, and zero otherwise. Write the (skinny) eigendecomposition $P = U \Lambda U\t$, where $U \in \R^{n \times d}$ is a matrix whose columns are orthonormal eigenvectors for $P$, and $\Lambda$ is the $d \times d$ diagonal matrix whose main diagonal entries are the eigenvalues of $P$. If each community has at least one vertex, then the rank of $P$ is $d = \operatorname{rank}(B)$. Consequently, the eigenvectors of $P$ reveal the population-level community memberships as follows.
\begin{lemma}[Restatement of Lemma~2.1 of \citet{lei_consistency_2015}] \label{lem:Udifference}
    Suppose $B \in [0,1]^{K \times K}$ is rank $K$, and let $U$ be the leading $K$ eigenvectors of $P = \rho_{n} Z B Z\t$. There exists a $K \times K$ matrix $R$ such that $U = Z R$. Furthermore, in Euclidean norm, $\| R_{k*} - R_{l*}\| = \sqrt{n_{k}\inv + n_{l}\inv}$ for all $k \neq l$.
\end{lemma}
Given two vertices $i$ and $j$, it follows that the row differences of $U$ satisfy
  \begin{align}
    \|U_{i*} - U_{j*} \| 
    &=
      \begin{cases}
        \sqrt{n_{\tau(i)}\inv + n_{\tau(j)}\inv}
        &
        \text{~if~}\tau(i) \neq \tau(j), \\
        0
        &
        \text{~if~}\tau(i) = \tau(j).
      \end{cases}
  \label{eq:dif}
  \end{align}
When $B$ is instead rank-degenerate, the row differences can at times be shown to satisfy a similar property.
\begin{lemma}[Restatement of Lemma 2 of \citet{zhang_randomized_2022}]
    Suppose that $P = \rho_n ZBZ\t$, where $B$ satisfies 
   \begin{align*}
      \min_{k,l} \left\{B_{kk} + B_{ll} - 2 B_{kl}\right\}
      &
      \geq
      \eta > 0.
  \end{align*}
  If $\tau(i) \neq \tau(j)$, then
   \begin{align}
      \| U_{i*} - U_{j*} \|
      \geq
      \sqrt{\frac{\rho_n \eta}{\lambda_{\max}(P)}}.
      \label{eq:dif2}
  \end{align}
\end{lemma}
More general conditions for rank-degenerate matrices $B$ are possible such that a lower bound of the form \cref{eq:dif2} holds. See Lemma~1 of \citet{zhang_randomized_2022} for further details. If $B$ is fixed as $n \to \infty$, then $\lambda_{\max}(P) \asymp n\rho_n$ and $\eta \asymp 1$, so in this regime the bound is of order $n^{-1/2}$, which matches the lower bound in \cref{lem:Udifference} for full-rank matrices $B$.

In either the full rank or rank-degenerate case, there are $K$ unique rows in $U$ corresponding to each of the different communities. Consequently, there are only $K$ unique rows of the matrix $U|\Lambda|^{1/2}$, hence $U|\Lambda|^{1/2} = Z \tilde R$ for some $K \times d$ matrix $\tilde R$ with unique rows corresponding to community memberships.

An alternative approach is to consider the normalized Laplacian matrix, with population version given by $\mathcal{L}(P) = (\operatorname{diag}(P1))^{-1/2}$ $P(\operatorname{diag}(P1))^{-1/2}$. 
The normalized Laplacian has a long history in spectral graph theory, capturing relationships between spectral matrix properties and graph connectivity \citep{chung_spectral_1997}. Alternatively, $I - (\operatorname{diag}(P1))^{-1/2} P (\operatorname{diag}(P1))^{-1/2}$, but this formulation simply shifts the eigenvalues without affecting the eigenspaces. Writing $\mathcal{L}(ZBZ\t)$ with eigendecomposition $U^{(L)} \Lambda^{(L)} (U^{(L)})\t$, the following result provides an analogue of \cref{lem:Udifference} for the normalized Laplacian.

\begin{lemma}[Lemma~3.1 of \citet{rohe_spectral_2011}]
    If $B$ is full rank, then there exists $R^{(L)} \in \mathbb{R}^{K \times K}$ such that $U^{(L)} = Z R^{(L)}$.  
\end{lemma}

The Euclidean norm of the difference between the rows of the matrix $R^{(L)}$ is not explicitly given above, but examining the proof, matching notation, and comparing it to \cref{lem:Udifference} reveals that \cref{eq:dif} holds. Namely, $\|R_{k*}^{(L)} - R_{l*}^{(L)}\| = \sqrt{n_k\inv + n_l\inv}$. Hence, for vertices $i$ and $j$,
\begin{align}
    \|U_{i*}^{(L)} - U_{j*}^{(L)} \| 
    &=
    \begin{cases}
        \sqrt{n_{\tau(i)}\inv + n_{\tau(j)}\inv}
        &
        \textnormal{~if~}\tau(i) \neq \tau(j), \\
        0
        &
        \textnormal{~if~}\tau(i) = \tau(j).
    \end{cases} \label{eq:dif3}
\end{align}

In summary, knowledge of the eigenvector matrix $U^{(L)}$ is again sufficient to recover the community memberships, since there are $K$ distinct rows of the eigenvector matrix $U^{(L)}$ and repeated rows corresponds to the same community. As with the population adjacency matrix, directly examining the proof also reveals that $U^{(L)} |\Lambda^{(L)}|^{1/2} = Z \tilde R^{(L)}$, where $\tilde R^{(L)} \in \R^{K \times K}$.

The above observations together show that the spectral embedding of each population (expected) similarity matrix reveals community information, namely that the population spectral embedding is a fixed transformation of the community membership matrix $Z$. Provided the empirical similarity matrices approximate their population counterparts, in a sense to be quantified in the following subsection, then the corresponding empirical spectral embeddings will approximate their population counterparts.

\subsection{Subspace Perturbation Approach to Community Detection} \label{sec:subspace_perturbation}
Several concepts of consistent community detection are studied in the literature on large graphs. \emph{Partial recovery} refers to recovering vertex memberships better than chance with probability tending to one. \emph{Weak recovery} refers to recovering all but a vanishing fraction of vertex memberships with probability tending to one. \emph{Exact recovery} refers to (asymptotically) recovering all vertex memberships exactly with probability tending to one. For additional context, \citet{gao_minimax_2021} surveys minimax rates of estimation and community detection in networks, while \citet{abbe_community_2018} surveys the information-theoretic and computational limits of community detection. The present survey focuses on the exact recovery regime which admits a rich landscape of asymptotic normality results.

The spectral embeddings of population-level similarity matrices distinguish community memberships, so consistency can be achieved provided that the empirical spectral embeddings are sufficiently close to their respective population counterparts. We will first focus on the setting where the spectral embedding is the $n \times d$ matrix of leading $d$ eigenvectors of the similarity matrix, though the underlying methods used to prove some of the results here can also be used to derive results for scaled eigenvectors.

Below, let $\hat U$ denote the leading eigenvectors computed from an observed, generic similarity matrix, and let $U$ denote the leading eigenvectors of the underlying population similarity matrix. The population and empirical eigenvectors of the adjacency matrix are of dimension $n \times d$, where $d$ is the rank of $B$ with $d \ll n$, where $d \leq K$ though $B$ need not be full rank. Consequently, the matrix of eigenvectors grows in size as $n$ increases, therefore the convergence of empirical eigenvectors must be defined with respect to an appropriate metric. Since orthogonal projection matrices uniquely represent orthogonal projections onto lower-dimensional subspaces, metrics on projection matrices offer a natural distance between subspaces. At the same time, cluster information is contained in the rows of the eigenvector matrix $U$. It is therefore natural to consider $\inf_{W \in \mathbb{O}(d)} \| \hat U - UW \|_{\eta}$, where $\| \cdot \|_{\eta}$ denotes a (often unitarily invariant) norm on matrices such as the spectral or Frobenius norm, and $\mathbb{O}(d)$ denotes the set of $d \times d$ orthogonal matrices, namely square matrices $W$ satisfying $W W\t = I_d$. For both the spectral norm and Frobenius norm, classical matrix perturbation theory \citet{bhatia_matrix_1997} provides upper and lower bounds of the form
\begin{align}
    c_{\eta}
    \|\hat U \hat U\t - U U\t \|_{\eta}
    &\leq 
    \inf_{W \in \mathbb{O}(d)} \| \hat U - UW \|_{\eta}
    \leq
    C_{\eta}
    \|\hat U \hat U\t - U U\t \|_{\eta}, 
    \label{spectral} 
\end{align}
where $c_{\eta}$ and $C_{\eta}$ are universal constants depending only on the choice of norm $\eta$. For example, $c_{\operatorname{F}} = 1/\sqrt{2}$ and $C_{\operatorname{F}} = 1$. See Lemma~1 in \citet{cai_rate-optimal_2018} or Lemma~2.1.3 in \citet{chen_spectral_2021}.
Of note, these distances are also closely related to the $\sin\Theta$ distances between subspaces, which are defined via $\|\sin\Theta(\hat U, U)\|_{\eta} = \| (\bm{I} - UU\t) \hat U \|_{\eta}$. Similar to \eqref{spectral}, classical matrix perturbation theory shows the equivalence between $\sin\Theta$ distances and the metrics above. In the setting herein, all of these metrics are equivalent up to constant factors, meaning they generate the same metric topology on the set of $K$-dimensional subspaces.

The Frobenius and spectral norm differences between $\hat U$ and $U$ quantify overall error. Bounding these differences can be shown to yield weak recovery (also called weak consistency) in certain SBMs when using $(1+\epsilon)$ approximate $k$-means clustering, based on the paradigm in \cite{lei_consistency_2015}. By employing bounds of the above form, \cite{lei_consistency_2015} proves that a vanishing fraction of vertices are misclustered with high probability for certain large SBM graphs.

To establish strong recovery or perfect clustering, one requires uniform control over the rows of the difference between $\hat U$ and $U$. Such uniform control can be established using the $\ell_{2,\infty}$ norm and the entrywise $\max$ norm. Traditionally, the $\ell_{2,\infty}$ norm of a matrix $M \in \R^{n\times d}$ is defined as $\|M\|_{2,\infty} \coloneqq \sup_{\|x\|_2 = 1} \| Mx \|_{\infty} \equiv \max_{1\leq i\leq n} \| M_{i*} \|$, i.e.,~the maximum row-wise Euclidean norm. Consequently, numerous works seek to obtain high-probability bounds for the expression
\begin{align}
    \inf_{W \in \mathbb{O}(d)}
    \|&\hat U - UW \|_{2,\infty}
      & =
    \inf_{W \in \mathbb{O}(d)}
    \max_{1 \leq i \leq n}
    \|(\hat U - UW)_{i*} \| \label{twoinfty1}
\end{align}
which is typically achieved indirectly by bounding the proxy quantity
\begin{align*}
    \|\hat U - UW^{\star}\|_{2,\infty}, 
    \quad
    &\text{~where~}
    \quad
    W^{\star}
    =
    \underset{W \in \mathbb{O}(d)}{\arg\inf}
    \|\hat U - U W\|_{\operatorname{F}}.
\end{align*}
Above, $W^{\star}$ denotes the solution to the orthogonal Procrustes problem under the Frobenius norm.

The recent monograph \citet{chen_spectral_2021} presents a broadly applicable entrywise bound which is restated here in the particular setting of stochastic blockmodels.

\begin{theorem}[Specialization of Theorem~4.2.1 in \citet{chen_spectral_2021} for SBMs] \label{thm:twoinfty} Fix $K$ and $B = \bm{\nu\nu}\t$ for some $K \times d$ matrix $\nu$. Suppose that $n\rho_n \gtrsim \log(n)$ and that $n_k \asymp n_l$ for all $1 \leq l,k\leq K$. Suppose $A \sim P$ with $P = \rho_n ZBZ\t$, and let $\hat U$ and $U$ denote the leading eigenvectors of $A$ and $P$, respectively. There exists an orthogonal matrix $W$ such that with probability exceeding $1 - O(n^{-5})$, 
    \begin{align*}
         \max_{1 \leq i \leq n}
         \|(\hat U - U W)_{i*} \|
         =
         O\left(\frac{\sqrt{d \log(n)}}{n \sqrt{\rho_n}}\right).
    \end{align*}
\end{theorem}
Obtaining the above result follows from matching notation, observing that the eigenvalues of $\rho_n ZBZ\t$ grow at order $n\rho_n$, noting that the standard deviation of each entry is bounded by $\rho_n^{1/2}$, and the fact that $U$ is $\mu$-incoherent with $\mu = O(1)$ when $n_k \asymp n_l$ (see \cref{eq:dif}). A slightly stronger result can also be stated which yields exact recovery down to the optimal information-theoretic threshold including explicit constants; see Theorem~4.5.1 in \citet{chen_spectral_2021}.

The sparsity term $\rho_n$ is used to control graph density, balanced against graph size, where often $\rho_n \to 0$ but $n\rho_n \to \infty$. If $n\rho_n \ll \log(n)$, then it is well-known that the adjacency matrix will no longer concentrate about $P$ in spectral norm due to disconnectedness of the underlying graph, whereas if $n\rho_n \gtrsim \log(n)$, then $\| A - P \| \lesssim \sqrt{n\rho_n}$ with high probability (see \citet{lei_consistency_2015} or Remark~3 of \citet{bandeira_sharp_2016}. See also \citet{lu_spectra_2013} for results concerning the Laplacian). Furthermore, the information-theoretic limit to achieve perfect clustering is in the regime $n\rho_n \asymp \log(n)$ with explicit constant determined in some important special cases \citep{abbe_community_2018}.

A host of works have studied eigenvector perturbations in statistical models, often in the context of network analysis \citep{abbe_entrywise_2020,abbe_ell_p_2022,agterberg_entrywise_2022,cai_subspace_2021,cape_signal-plus-noise_2019,cape_two--infinity_2019,chaudhuri_spectral_2012,eldridge_unperturbed_2018,jin_mixed_2023,lei_unified_2019,lei_consistency_2020,rohe_vintage_2023,sarkar_role_2015,yan_inference_2024,zhang_perturbation_2022}. For the setting of stochastic blockmodels with fixed $B$ matrix, these results all translate to similar statistical consequences. If $n_k \asymp n_l$ for all $k,l$, and $n\rho_n \gtrsim \log(n)$, then the row-wise error between $\hat U$ and $U$, after proper orthogonal alignment, is exceedingly small. Since the rows of $U$ reveal the community memberships per \cref{eq:dif}, these results suggest exact community recovery with probability tending rapidly to one for sufficiently dense graphs, namely when $n\rho_n \gtrsim \log(n)$. Moreover, consistent estimation of the $i$-th row of $U$ uniformly over all rows $1 \le i \le n$ paves the way for further downstream inference and cluster analysis.

\section{Asymptotic Normality in Stochastic Blockmodels}\label{sec:asymp_norm}
Asymptotic normality provides a more refined characterization of cluster geometry and spectral embeddings beyond consistency. Further, distributional theory is conventionally a precursor of inference.

\subsection{Asymptotic Normality of the Adjacency Matrix Eigenvectors}\label{sec:eigenvectors}
We begin by highlighting row-wise asymptotic normality for the $n \times K$ matrix $\hat U$ of leading orthonormal eigenvectors of the stochastic blockmodel adjacency matrix $A$.

\begin{theorem}[Restatement of Theorem 4.4 in \citet{xie_entrywise_2024} for positive semidefinite SBMs] \label{thm:xieCLT}
    Fix $K$ and $B = \bm{\nu} \bm{\nu}\t$ for some $K \times d$ matrix $\bm{\nu}$. Suppose that $n\rho_n \gtrsim\log(n)$ and that $n_k \asymp n_l$ for all $1 \le k, l \le K$. Suppose $A \sim P$ with $P = \rho_n ZBZ\t$, and let $\hat U$ and $U$ denote the leading eigenvectors of $A$ and $P$, respectively. Suppose further that there exists a constant $\delta > 0$ such that $\min_{1 \leq i \leq n} (n \rho_n)^{-1}\sum_{j} P_{ij} \geq \delta$. Define
    \begin{align*}
        \mathbf{\Sigma}_U^{(n)}(k)
        &\coloneqq
        \Delta_n^{-3/2} \Gamma^{(k)}_n \Delta_n^{-3/2},
    \end{align*}
    where
    \begin{align*}
        \Delta_n
        &\coloneqq
        \sum_{l=1}^K \frac{n_l}{n} \bm{\nu}_{l*} \bm{\nu}_{l*}\t, \qquad 
        \Gamma^{(k)}_n
        \coloneqq
        \sum_{l=1}^K \frac{n_l}{n} B_{kl} (1 - \rho_n B_{kl}) \bm{\nu}_{l*} \bm{\nu}_{l*}\t.
    \end{align*}
    Then, there exist sequences of orthogonal matrices $\tilde W$ and $W$ depending on $n$ and a constant $C > 0$ depending only on $d$ such that, for each fixed index $i$, for any convex measurable set $\mathcal{A} \in \mathbb{R}^{d}$, it holds that for $n$ sufficiently large,
    \begin{align*}
        \bigg| \p\bigg\{ n \rho_n^{1/2} \big( \mathbf{\Sigma}_U^{(n)}(k) \big)^{-1/2} \tilde W\t &\big( ( \hat U W )_{i*} - U_{i*} \big) \in \mathcal{A} \bigg\} 
        - \p \{ G \in \mathcal{A}\} \bigg|
        \leq
        \frac{C \log(n\rho_n)}{(n\rho_n)^{1/2}}
    \end{align*}
    where $G$ denotes a $d$-dimensional standard Gaussian random vector.
\end{theorem}

\cref{thm:xieCLT} establishes that the rows of the empirical eigenvectors $\hat U$ are approximately Gaussian about the rows of $U$ and with covariance determined by $\mathbf{\Sigma}_U^{(n)}(k)$. The rate of convergence of the cumulative distribution functions clarifies the role of $n\rho_n$ as the effective sample size. For classical finite-dimensional Berry--Esseen theorems, the rate of convergence to the Gaussian cumulative distribution function is of order $n^{-1/2}$, whereas in \cref{thm:xieCLT}, it is of order $(n\rho_n)^{-1/2}$ modulo logarithmic terms. Of note, the full result in \citet{xie_entrywise_2024} is more general, permitting positive semidefinite $P$ matrices outside of blockmodels, and allowing both $d$ and $B$ to change with $n$.

A key property of the covariance matrix $\mathbf{\Sigma}^{(n)}_U(k) \in \R^{K \times K}$ is that it depends only on the membership $\tau(i) = k$, the community probabilities $\{\pi_k\}_{k=1}^{K}$ and the block probability matrix $B$. Per \cref{sec:populations}, there are at most $K$ unique values of $U_{i*}$, and here, each distinct row $U_{i*}$ gives rise to a potentially different block-specific covariance matrix $\mathbf{\Sigma}^{(n)}_U(k)$. Therefore, unconditionally, the random vector $\hat U_{i*}$ asymptotically behaves like a mixture of Gaussians with component distributions associated to the different community centroids. In particular, when the communities are stochastically generated according to $\p[\tau(i) = k] = \pi_k$, then the matrices $\mathbf{\Delta}_n$ and $\mathbf{\Gamma}^{(k)}_n$ in \cref{thm:xieCLT} have well-defined limits with probability one by the law of large numbers.\cref{thm:xieCLT} generalizes an earlier result in \citet{cape_signal-plus-noise_2019} which is applicable to slightly denser SBMs and directly assumes that that $\mathbf{\Delta}_n$ and $\mathbf{\Gamma}^{(k)}_n$ have well-defined limits, hence so too does $\mathbf{\Sigma}_U^{(n)}(k)$.  More recently, further generalization is available in \citet{xie_higher-order_2024}.

When all community sizes are comparable (which holds when relative community membership proportions are held fixed), then per \cref{eq:dif,eq:dif2}, the row norms of the population eigenvectors all necessarily converge to zero and at the same rate, thereby reflecting \emph{incoherence} \citep{chen_spectral_2021}. This shrinking-to-zero property is reflected in the scaling $n\rho_n^{1/2}$ required to obtain a Gaussian limit. Informally, this scaling decomposes as the product $n^{1/2} \times (n\rho_n)^{1/2}$, where the term $n^{1/2}$ accounts for incoherence and the term $(n\rho_n)^{1/2}$ is the square root of the order of each node's expected degree. In the absence of sparsity, namely in the fully dense regime $\rho_{n} \equiv 1$, the classical $n^{1/2}$ parametric rate is obtained after adjusting for the eigenvector incoherence.

The $k$-th entry of the $i$-th row of $\hat U$ is expressible as the inner product of the $k$-th eigenvector $\hat U_{*k}$ with the standard basis vector $e_i$, namely $e_i\t \hat U e_k = e_i\t \hat U_{* k}$. Recently, \citet{fan_asymptotic_2022} develops asymptotic theory for general linear forms $x\t \hat U_{*k}$ where $x \in \R^{n}$ is a user-specified deterministic unit vector. Theorem~2 therein establishes that, as $n \rightarrow \infty$,
\begin{align*}
    \frac{1}{s_{x,k}} \bigg\{ t_k \bigg(x\t \hat U_{*k} - \mathbf{A}_{x,k,t_k} \mathcal{P}_{k,t_k}^{1/2}\bigg) - \E\left[(b_{x,k,t_k} - x\t U_{*k}U_{*k}\t) (A-\E[A]) U_{*k} \right] \bigg\} 
   & \dto
    \mathcal{N}(0,1),
    \numberthis
    \label{fanCLT}
\end{align*}
where
\begin{align*}
    s_{x,k}^2 &= \operatorname{Var}\left[(b_{x,k,t_k}\t - x\t U_{*k} U_{*k}\t )(A-\E[A]) U_{*k} \right]
\end{align*}
with $b_{x,k,t_k}$, $\mathbf{A}_{x,k,t_k}$, and $\mathcal{P}_{k,t_k}$ defined according to population parameters and the additional terms $x$, $k$, and $t_k$. Explicitly, $t_k$ solves the fixed-point equation
\begin{align*}
    1 + \lambda_k\big\{ \mathcal{R}( U_{*k}, U_{*k},z) - \mathcal{R}(U_{*k},U_{-k},z)  \big[ \Lambda_{-k}\inv + \mathcal{R}( U_{-k}, U_{-k},z) \big]\inv \mathcal{R}(U_{-k}, U_{*k},z) \big\} \overset{\textnormal{set}}{=}
   0
   \numberthis
   \label{fixedpoint}
\end{align*}
over an interval $z \in [a_k,b_k]$ that contains $\lambda_k$ but no other eigenvalues of $\E[A]$, thereby precluding the possibility of repeated eigenvalues. Here, $\mathcal{R}(\cdot,\cdot,\cdot)$ is defined as
\begin{align*}
    \mathcal{R}( M_1, M_2,z)
    &\coloneqq
    - \sum_{l=0,l\neq 1}^L z^{-(l+1)}
    M_1\t \mathbb{E}\left[(A - P)^l\right] M_2,
\end{align*}
where $L$ denotes a sufficiently large integer per the proof, and where $U_{-k}$ and $\Lambda_{-k}$ denote the matrices $U$ and $\Lambda$ but with the entries associated to the $k$-th eigenvector and eigenvalue removed, respectively. See \cref{sec:eigenvalues} for further discussion of the quantity $t_k$. Lemma~3 of \citet{fan_asymptotic_2022} shows that $t_k$ exists, is unique, and satisfies $t_k/\lambda_k \to 1$ as $n \to \infty$. The definitions of $b_{x,k,t_k}$, $\mathbf{A}_{x,k,t_k}$, and $\mathcal{P}_{k,t_k}$ appearing in \cref{fanCLT} are involved yet in special cases admit simple asymptotic expressions. Namely, take $x = U_{*k}$, in which case the following simplification describes the cosines of the angles between the empirical and population eigenvectors.

\begin{theorem}[Restatement of Part~2 of Theorem~2 in \citet{fan_asymptotic_2022} for SBMs] \label{thm:fanCLT2}
    Fix $K$ and $B$. Suppose that $n\rho_n \gg n^{\eta}$ for some constant $\eta > 0$ and that $n_k \asymp n_l$ for all $1 \le k, l \le K$. Suppose $A \sim P$ with $P = \rho_n ZBZ\t$, and let $\hat U$ and $U$ denote the leading eigenvectors of $A$ and $P$, respectively. Suppose further that $P$ has $K$ distinct eigenvalues $\lambda_1, \dots, \lambda_K$ satisfying $\min_{i\neq j} |\lambda_i|/|\lambda_j| \geq 1 + c_0$ for some constant $c_0 > 0$ for $\lambda_i \neq - \lambda_j$. Let $t_k$ be the solution to the fixed-point equation in \cref{fixedpoint}.
    Then, as $n \rightarrow \infty$,
    \begin{align*}
        \frac{2 t_k^2 \left( U_{*k}\t \hat U_{*k} + A_{U_{*k},k,t_k} \mathcal{P}_{k,t_k}^{1/2}\right)}{\operatorname{Var}\left[ U_{*k}\t (A - P)^2 U_{*k}\right]^{1/2}}
        \dto
        \mathcal{N}(0,1),
    \end{align*}
    where the centering term has the expansion
    \begin{align*}
        A_{U_{*k},k,t_k} \mathcal{P}_{k,t_k}^{1/2}
        &=
        - 1
        + \frac{1}{2 t_k^2}  U_{*k}\t \mathbb{E}\left[(A - P)^2\right] U_{*k} + O\bigg( \frac{1}{(n\rho_n)^{3/2}} \bigg).
    \end{align*}
\end{theorem}

\cref{thm:fanCLT2} demonstrates that the cosines of the angles between empirical and population eigenvectors are asymptotically Gaussian with mean one minus a bias term stemming from second-order moment properties via $U_{*k}\t \mathbb{E}\left[(A - P)^2\right] U_{*k}$ and variance determined by fourth-order moment properties via $ \operatorname{Var}\left[U_{*k}\t (A - P)^2 U_{*k}\right]$. The presentation in \cref{thm:fanCLT2} differs from the original statement in \citet{fan_asymptotic_2022} but holds by observing that the eigenvalues of $P = \rho_n ZBZ\t$ are of order $n\rho_n$ when $B$ is fixed, the parameter $\alpha_n$ appearing in \citet{fan_asymptotic_2022} is of order $(n\rho_n)^{1/2}$, and that the eigenvectors $U_{*k}$ are incoherent when $n_k \asymp n_l$. Of note, \citet{fan_asymptotic_2022} covers other situations, including demonstrating different asymptotic phenomena based on whether the deterministic vector $x$ is incoherent.

To elucidate \cref{thm:fanCLT2}, we remark that Lemma 3 of \citet{fan_asymptotic_2022} together with the present assumptions implies that $t_k$ is of order $n\rho_n$. Furthermore, by a slight modification of Remark 1 in \citet{fan_asymptotic_2022}, it follows that the asymptotic variance quantity $\operatorname{Var}[U_{*k}\t (A - P)^2 U_{*k}]$ is of order $n\rho_n$. Thus, the variance quantity satisfies
\begin{align*}
    \frac{\operatorname{Var}\big[ U_{*k}\t (A - P)^2 U_{*k} \big]}{4 t_k^4}
    \asymp
    \frac{n\rho_n}{(n\rho_n)^4}
    =
    \frac{1}{(n\rho_n)^3}.
\end{align*}
Informally, $U_{*k}\t \hat U_{*k} - 1$ is approximately Gaussian with bias of order $(n\rho_n)^{-1}$ and variance of order $(n\rho_n)^{-3}$.

By \cref{spectral}, the eigenvector difference $\hat U - UW$ is closely related to the projection norm difference, and hence, its distributional asymptotics are closely related to those of the matrix $\hat U \hat U\t - UU\t$. The theory in \citet{fan_asymptotic_2022} reflects this relationship, namely the limit results therein for the linear form $x\t \hat U_{*k}$ are based on first developing asymptotic expansions for the term $x\t (\hat U_{*k}\hat U_{*k}\t - U_{*k}U_{*k}\t)y$ for sequences of deterministic unit vectors $x$ and $y$, which is the bilinear form determined by the difference of projections onto the $k$-th dominant eigenspace, in contrast to the row-wise focus jointly for the collection of leading eigenvectors in \cref{thm:xieCLT}.

The results in \citet{fan_asymptotic_2022} permit $P$ to have negative eigenvalues but require that the eigenvalues be distinct so as to explicitly quantify the fluctuations of the individual eigenvectors $\hat U_{*k}$. While \cref{thm:xieCLT} also requires that $P$ is positive semidefinite, it permits repeated eigenvalues, hence it can only quantify the fluctuations of the rows (as opposed to entries) of the matrix of leading eigenvectors $\hat U$ up to orthogonal transformation (as opposed to individual eigenvectors $\hat U_{*k}$). The results in \citet{fan_asymptotic_2022} require the sparsity condition $n\rho_n \gg n^{\eta}$ for some {(fixed)} $\eta > 0$, whereas the results in \cref{thm:xieCLT} require $n \rho_{n} \gtrsim \log(n)$.

Rather than investigating the rows or entries of $\hat{U}$ or $\hat{U}_{*k}\hat{U}_{*k}\t$, \citet{li_two-sample_2018} considers the asymptotic behavior of a Frobenius norm difference statistic between $\hat{U}$ and $U$, modulo orthogonal transformation. For balanced homogeneous stochastic blockmodels, simple analytical expressions hold for the centering and scaling terms, and the overall expression is closely related to the $\sin\Theta$ distance.

\begin{theorem}[Restatement of Theorem~3 in \citet{li_two-sample_2018}] \label{thm:liliCLT}
    Fix $K$, and suppose that the $K \times K$ matrix $B$ has the form 
    \begin{align*}
        B_{kl}
        &=
        \begin{cases}
            a
            &\textnormal{~if~}k = l;
            \\ 
            b
            &\textnormal{~if~}k \neq l,
        \end{cases}
    \end{align*}
    where $0 < b < a < 1$ are constants. Assume that  $n\rho_n \gtrsim n^{1/2 + \eps}$ for some $\eps > 0$. Suppose $A \sim P$ with $P = \rho_n ZBZ\t$, and let $\hat U$ and $U$ denote the leading eigenvectors of $A$ and $P$, respectively. Further suppose that each community has precisely $n/K$ vertices. Then, there exists a sequence of orthogonal matrices $W$ such that, as $n \rightarrow \infty$,
    \begin{align*}
        \frac{ \frac{1}{Kn} \left\| \left(\hat U W - U\right)\Lambda \right\|_{\operatorname{F}}^{2} - \mu(a,b) }{\sigma(a,b)}
        \dto
        \mathcal{N}(0,1),
    \end{align*}
    where $\mu(a,b)$ and $\sigma(a,b)$ are given by
    \begin{align*}
        \mu(a,b)
        &\coloneqq
        \rho_n b(1-b \rho_n )+
        \frac{\rho_n a(1-a \rho_n ) - \rho_n b(1-b\rho_n )}{K},
        \\
        \sigma^2(a,b)
        &\coloneqq
        \frac{2}{nK} \bigg[ (b\rho_n (1-b\rho_n ))^2 +
        \frac{(a\rho_n (1-a\rho_n ))^2 - (b\rho_n (1-b\rho_n ))^2}{K} \bigg].
    \end{align*}
\end{theorem}

In fact, \citet{li_two-sample_2018} allows $a$ and $b$ to change with $n$, though this is omitted above for ease of presentation. Whether \cref{thm:liliCLT} can be generalized to the setting of general $B$ matrices and sparsity of order $O(n^{\eta})$ for $0 < \eta \leq 1/2$, or even $n\rho_n = O(\log(n))$, remains an open problem.

The random variable $\|(\hat U W - U)\Lambda\|_{\operatorname{F}}^{2}$ can be viewed as a columnwise re-weighted Frobenius $\sin\Theta$ distance with respect to the magnitude of the eigenvalues. While \cref{thm:fanCLT2} demonstrates that the cosines of the eigenvectors between $\hat U_{*k}$ and $U_{*k}$ have an asymptotically normal distribution after appropriate centering and scaling, the results of \citet{fan_asymptotic_2022} do not imply those in \citet{li_two-sample_2018}. In particular, the former result requires distinct eigenvalues, whereas the latter model setting gives rise to repeated population eigenvalues. Yet, both theorems reflect similar asymptotic normality properties for the (rescaled) sines or cosines between the empirical and true eigenvectors, suggesting that extensions to the regime $n\rho_n \gtrsim \operatorname{polylog}(n)$ or even $n\rho_n \asymp \log(n)$ might be possible.

\subsection{Asymptotic Normality of the Scaled and Laplacian Eigenvectors}\label{sec:scaled_eigenvectors}
Another graph embedding approach is to scale the leading eigenvectors by the square roots of the corresponding eigenvalue magnitudes. This scaling has a convenient interpretation via \emph{latent space modeling}, namely that each vertex has a latent low-dimensional Euclidean representation (i.e., vector or position). The general latent space model posits that the entries of the $P$ matrix are obtained by first generating $n$ independent latent vectors $X_1,\dots,X_n$ from a distribution $F$ on $\R^d$ and then setting $P_{ij} = \rho_n \kappa(X_i, X_j)$, where $\kappa : \R^{d} \times \R^{d} \to [0,1]$ is a symmetric kernel function. The special case $\kappa(X_i,X_j) = \langle X_i, X_j\rangle$, yields the so-called \emph{random dot product graph} (model) \citep{young_directed_2008}, presented below.

\begin{definition}[Random dot product graph -- undirected, loopy]
    We say that a symmetric adjacency matrix $A \in \{0,1\}^{n \times n}$ is a random dot product graph on $n$ vertices with sparsity term $\rho = \rho_n \in (0,1]$ if, conditional on $X \in \mathbb{R}^{n \times d}$ satisfying $XX\t \in [0,1]^{n \times n}$, its upper triangular entries are generated independently in the manner
    \begin{align*}
        A_{ij}
        &=
        \begin{cases}
            \operatorname{Bernoulli}( \rho X_{i*}\t X_{j*} )
            &
            \textnormal{~if~}i\leq j;
            \\
            A_{ji}
            &
            \textnormal{~if~}i > j.
        \end{cases}
    \end{align*}
    The rows of $X$ may be deterministic or modeled in a generative fashion via i.i.d. vectors $X_{i*} \sim F_X$ for a distribution $F_X$ on $\mathbb{R}^d$ such that $X_{i\star}^{\top}X_{j\star} \in [0,1]$ for all $1 \le i, j \le n$.
\end{definition}

Random dot product graphs have the property that $P = \rho_n XX\t$, where each row of $X$ is interpreted as a latent position in $\R^d$. In this setting, $U\Lambda^{1/2} = \rho_n^{1/2} X \tilde W$ for some orthogonal matrix $\tilde W$ (not necessarily identifiable). Consequently, the scaled eigenvectors correspond to an isometric linear transformation and scaling of the original latent vectors. The latent-space interpretation avoids the previous shrinking-to-zero phenomenon from the eigenvectors in $n$ and introduces additional regularity into the definition of $P$.

Here, $\p[A_{ij} = 1] = \rho_n X_{i*}\t,X_{j*}$ leads to having matrices $B$ that are necessarily positive semidefinite. More recently, \citet{rubin-delanchy_statistical_2022} introduces the \emph{generalized random dot product graph} model, which replaces the inner product with a pseudo-inner product to model large-in-magnitude negative eigenvalues of the adjacency matrix $A$ and hence allow for indefinite $B$ matrices.

\begin{definition}[Generalized random dot product graph --- undirected, loopy] \label{def:grdpg}
    We say that a symmetric adjacency matrix $A$ is a generalized random dot product graph on $n$ vertices with sparsity term $\rho = \rho_n \in (0,1]$ if, conditional on $X \in \mathbb{R}^{n \times d}$ satisfying $X \ipq X\t \in [0,1]^{n \times n}$, its upper triangular entries are generated independently in the manner
    \begin{align*}
        A_{ij}
        =
        \begin{cases}
            \operatorname{Bernoulli}(\rho X_{i*}\t \ipq X_{j*}),
            &
            \textnormal{~for~}i\leq j;
            \\
            A_{ji}
            &
            \textnormal{~for~}i > j.
        \end{cases}
    \end{align*}
    Here, $\ipq$ is the matrix $\operatorname{diag}(I_p,-I_q)$ with $p+q = d$. The rows of $X$ may be deterministic or modeled in a generative fashion via i.i.d. vectors $X_{i*} \sim F_X$ for a distribution $F_X$ on $\R^d$ satisfying $X_{i*}\t \ipq X_{j*} \in [0,1]$ for all $1 \le i, j \le n$.
\end{definition}

For GRDPGs, as for RDPGs, the latent positions can be viewed as following a mixture distribution of point masses, namely $X_{i*} = \bm{\nu}_{k*}$ if $\tau(i) = k$, defined by $B = \bm{\nu} \ipq \bm{\nu}\t$, with non-negative integers $p$ and $q$ representing the number of positive and negative eigenvalues of $B$, respectively. Under the generalized random dot product graph model, the matrix $P$ can be decomposed as $\rho_n X \ipq X\t$, but the definition induces an additional nonidentifiability. If $Q$ satisfies $Q \ipq Q\t = \ipq$, namely if $Q$ is a so-called \emph{indefinite orthogonal matrix}, then $X$ and $XQ$ generate the same probability matrix, implying that the latent positions are only identifiable up to such transformations $Q$. Consequently, analyzing the generalized dot product graph requires additional considerations with respect to these transformations. For random dot product graphs, namely when $q = 0$, these matrices $Q$ are simply \emph{orthogonal matrices}.

A survey of statistical inference under the random dot product graph can be found in \citet{athreya_statistical_2018}. Our focus herein is on asymptotic normality for (generalized) random dot product graphs through the lens of stochastic blockmodels.

\subsubsection{Latent Space Interpretation of the Stochastic Blockmodel}
Interpreting the stochastic blockmodel as a latent space model amounts to imposing an additional assumption on the distribution $F$ if the communities are stochastic or on the rows of the matrix $X$ if the communities are deterministic. Namely, either $F$ corresponds to a mixture of $K$ point masses in $\R^d$, where $d = \operatorname{rank}(B)$, or $X$ has at most $K$ unique $d$-dimensional rows, with each row corresponding to community membership in the latter if the communities are fixed. One approach to specify the distribution $F$ is via the eigendecomposition of the matrix $B$. If $B$ is positive semidefinite, write its eigendecomposition as $B = VDV\t$, and define the matrix $\bm{\nu} \coloneqq VD^{1/2}$. Then, the rows of $X$ can be drawn via $\p[X_{i*} = \nu_{k*}] = \pi_k$ hence $\rho_n XX\t = \rho_n ZBZ\t$, so the stochastic blockmodel is a random dot product graph.

If $B$ is not positive semidefinite, then one must consider a generalized random dot product graph, with $p$ and $q$ equal to the number of positive and negative eigenvalues of $B$, respectively. In this case, take $\bm{\nu} \coloneqq V|D|^{1/2}\ipq$ so that $B = V |D|^{1/2} \ipq |D|^{1/2} V\t$. In fact, the latent space interpretation also applies to the stochastic blockmodel extensions mentioned in \cref{sec:subspace_perturbation}. For the mixed-membership stochastic blockmodel, the latent space interpretation is that the latent positions belong to a simplex in $d$ dimensions, and for the degree-corrected stochastic blockmodel, the latent positions are one of $K$ rays emanating from the origin, where each ray corresponds to a different community and location along the ray is governed by the node degree parameters.

Denote the spectral decomposition of $A$ by $A = \hat U \hat \Lambda \hat U\t + \hat U_{\perp} \hat \Lambda_{\perp} \hat U_{\perp}\t$, where $\hat \Lambda$ is the $d \times d$ diagonal matrix of leading eigenvalues in magnitude of $A$, and where $\hat U$ is the $n \times d$ matrix whose columns are the orthonormal eigenvectors corresponding to $\hat \Lambda$. Similarly, write $P = U \Lambda U\t$. Under certain conditions, existing row-wise uniform error bounds establish that $\hat U |\hat \Lambda|^{1/2} \approx U |\Lambda|^{1/2} W$ for a suitable orthogonal transformation $W$. In the special case when $P$ is positive semidefinite with $P = \rho_n XX\t$, then $U\Lambda^{1/2} = \rho_n^{1/2} X \tilde W$ for some orthogonal matrix $\tilde W$, which in particular suggests that $\hat{U} |\hat \Lambda|^{1/2} \approx \rho_n^{1/2} X W'$, where $W'$ is a possibly different orthogonal transformation representing the product of different, underlying orthogonal transformations.

In fact, the matrix of scaled eigenvectors $\hat X \coloneqq \hat U |\hat \Lambda|^{1/2}$ exhibits similar uniform row-wise convergence guarantees to the unscaled eigenvectors \citep{athreya_statistical_2018,rubin-delanchy_statistical_2022}, as follows.

\begin{theorem}[Restatement of Theorem~1 in \citet{rubin-delanchy_statistical_2022}] \label{prdtwoinfty}
    Fix $K$ and $B = \nu \ipq \nu\t$. Assume that each $X_i$ is drawn independently according to $\p[X_i = \nu_{k*}] = \pi_k$, with $\nu_{k*}$ denoting the $k$-th row of $\nu$ and $d = \operatorname{rank}(\nu)$. 
   Suppose $A \sim P$ with $P = \rho_n X \ipq X\t= \rho_n ZBZ\t$. {Let $P$ have eigendecomposition $U \Lambda U\t$, let $\tilde X = U |\Lambda|^{1/2}$, and let $\hat{U}$ and $\hat \Lambda$ denote the leading eigenvectors and eigenvalues (in  magnitude) of $A$.} There exists a universal constant $c > 1$ and an orthogonal matrix $W$ such that provided $n\rho_n \gg \log^{4c}(n)$, then
    \begin{align*}
        \| \hat X W - \tilde X \|_{2,\infty}
        &=
        \tilde O\left(\frac{1}{\sqrt{n}}\right)
    \end{align*}
    holds with probability at least $1 - O(n^{-M})$ for any prespecified choice of constant $M > 0$. The value of $M$ influences the implicit constant in the upper bound.
\end{theorem}

Under the GRDPG model, there exists an indefinite orthogonal matrix $Q$ such that $XQ = \tilde X$.
For sufficiently dense graphs, the convergence of $\hat X$ to $XQ$ resembles the convergence of $\hat U$ to $U$ after accounting for the order of magnitudes. For sufficiently dense graphs, $n\rho_n \gg \operatorname{polylog}(n)$, there exist sequences of orthogonal transformations $W^X$ and $W^U$ depending on $n$ such that with high probability,
\begin{align*}
    \frac{\| \hat{X} W^X - \rho_n^{1/2} X {Q}\|_{2,\infty}}{\|\rho_n^{1/2} XQ\|_{2,\infty}}
    &=
    \tilde O\left( \frac{1}{\sqrt{n\rho_n}}\right), 
        \qquad 
    \frac{\| \hat{U} W^U - U \|_{2,\infty}}{\|U\|_{2,\infty}}
    =
    \tilde O\left( \frac{1}{\sqrt{n\rho_n}}\right),
\end{align*}
establishing that the relative row-wise error for both $\hat{U}$ and $\hat{X}$ decays to zero as the number of vertices grows. \cref{prdtwoinfty} also holds for RDPGs per \citet{xie_entrywise_2024} under the milder assumption $n\rho_n \gtrsim \log(n)$.

\subsubsection{Asymptotic Normality of the Scaled Eigenvectors}
The scaled eigenvectors $\hat{X}$ exhibit a similar central limit theorem as the unscaled eigenvectors, as demonstrated in the following result.

\begin{theorem}[Restatement of Theorem~2 in \citet{rubin-delanchy_statistical_2022} for SBMs] \label{thm:prdCLT}
    Fix $K$ and $B = \nu \ipq \nu\t$. Assume that each $X_i$ is drawn independently according to $\p[X_i = \nu_{k*}] = \pi_k$, with $\nu_{k*}$ denoting the $k$-th row of $\nu$ and $d = \operatorname{rank}(\nu)$. Suppose $A \sim P$ with $P = \rho_n X \ipq X\t= \rho_n ZBZ\t$. {Let $P$ have eigendecomposition $U \Lambda U\t$, let $\tilde X = U |\Lambda|^{1/2}$, and let $\hat{U}$ and $\hat \Lambda$ denote the leading eigenvectors and eigenvalues (in  magnitude) of $A$.}
    There exist sequences of orthogonal matrices $W$ and indefinite orthogonal matrices $Q$ and a universal constant $c > 1$ such that, if $n\rho_n \gg \log^{4c}(n)$, then as $n \rightarrow \infty$,
    \begin{align*}
        \bigg(n^{1/2} Q\t \big((\hat X W )_{i*} - \rho_n^{1/2} \bm{\nu}_{k*} \big) \bigg)
        \mid
        \tau(i) = k
       &\dto
        \begin{cases}
            \mathcal{N}(0, \mathbf{\tilde{\Sigma}}_S(k))
            &
            \textnormal{~if~}\rho_n \equiv 1, \\
            \mathcal{N}(0, \mathbf{\tilde{\Sigma}}_{S,o(1)}(k)) 
            &
            \textnormal{~if~}\rho_n \to 0, \end{cases}
    \end{align*}
    where the asymptotic covariance matrices satisfy
    \begin{align*}
        \mathbf{\tilde{\Sigma}}_S(k)
        =
        \ipq \Delta\inv \Gamma^{(k)} \Delta\inv \ipq, \qquad 
        \mathbf{\tilde \Sigma}_{S,o(1)}(k)
        &=
        \ipq \Delta\inv \Gamma^{(k)}_{o(1)} \Delta\inv \ipq,
        \numberthis
        \label{eq:cltrdpg}
    \end{align*}
    where 
    \begin{align*}
        \Delta
        &\coloneqq
        \sum_{l=1}^K \pi_l \bm{\nu}_{l*} \bm{\nu}_{l*}\t,
        \qquad 
        \Gamma^{(k)}
        \coloneqq
        \sum_{l=1}^K \pi_k B_{kl} \big( 1 - B_{kl} \big) \bm{\nu}_{l*} \bm{\nu}_{l*}\t, 
        \qquad
        \Gamma^{(k)}_{o(1)}
        \coloneqq
        \sum_{l=1}^K \pi_k B_{kl} \bm{\nu}_{l*} \bm{\nu}_{l*}\t.
    \end{align*}
\end{theorem}

A precursor to \cref{thm:prdCLT} appears in \citet{athreya_limit_2016} which proves a central limit theorem for the scaled eigenvectors of random dot product graphs ({$q = 0$}) under the assumption that the latent positions are sampled from a distribution $F$ on $\R^d$, that the second moment matrix $\E_{X_{1*} \sim F}[X_{1*}X_{1*}\t]$ has distinct eigenvalues, and that $\rho_n = 1$. This result is extended in \citet{athreya_statistical_2018} to the setting where the second moment matrix is permitted to have repeated eigenvalues. For positive semidefinite stochastic blockmodels with equal  community membership probabilities, the assumption of distinct eigenvalues of $\E[X_{1*} X_{1*}\t]$ is equivalent to the assumption that $B$ has distinct eigenvalues, but for unequal community membership probabilities the eigenvalues of $\E[X_{1*} X_{1*}\t]$ need not be the same as those of $B$. In fact, $\E[X_{1*} X_{1*}\t] = \sum_{k=1}^{K} \pi_k \bm{\nu}_{k*} \bm{\nu}_{k*}\t$, where the matrix $\bm{\nu}$ is the $K \times d$ matrix of distinct latent vectors associated to the corresponding stochastic blockmodel.
These results were later extended to the case $\rho_n \to 0$ in \citet{tang_limit_2018} prior to the further extension in \citet{rubin-delanchy_statistical_2022}.

In \cref{thm:prdCLT}, the asymptotic covariance matrices for both $\rho_n = 1$ and $\rho_n \to 0$ are extensions of those in the positive semidefinite case, reflecting the roles of the leading eigenvalue signs through the matrix $\ipq$ and introducing an indefinite orthogonal matrix $Q$. A key property of indefinite orthogonal matrices is that they preserve \emph{volume}, namely $|\det(Q_X \ipq Q_X\t)| = |\det(\ipq)| = 1$, so that the asymptotic Gaussianity is retained in the sense that the volume of the ellipse generated by these asymptotic covariance matrices is invariant to these transformations. See \citet{rubin-delanchy_statistical_2022} for more discussion about these matrices and properties pertaining to latent position random graphs.

\cref{thm:prdCLT} demonstrates that the rows of the matrix $\hat X$ are approximately Gaussian and centered around the rows of $\rho_n^{1/2} Q\t \nu_{k*}$. As before, there are only $K$ unique centroids corresponding to the possible community memberships. Results presented earlier only involve orthogonal transformations and so the cluster shape (ellipse) is determined by the asymptotic covariance matrix in the statement of the result, whereas in \cref{thm:prdCLT} the cluster shape is influenced by both the limiting asymptotic covariance as well as the indefinite orthogonal matrix $Q$.

As in \cref{thm:xieCLT} for the unscaled eigenvectors, the asymptotic covariance matrix $\mathbf{\Sigma}_{S}(k)$ in \cref{thm:prdCLT} depends on the community probabilities $\pi_k$, the membership of the $i$-th vertex, and the matrix $B$. Provided $\mathbf{\Sigma}_U(k) = \lim_{n\to\infty} \mathbf{\Sigma}_U^{(n)}(k)$ exists, then \cref{thm:xieCLT} has asymptotic covariance matrix of the form
\begin{align*}
    \mathbf{\Sigma}_U(k)
    &=
    \Delta^{-3/2} \Gamma^{(k)} \Delta^{-3/2}.
\end{align*}
Consequently, when $B$ is positive semidefinite, \cref{thm:prdCLT} demonstrates that the asymptotic covariance matrices for the scaled and unscaled eigenvectors are related in the manner
\begin{align*}
    \mathbf{\Sigma}_S(k)
    &=
    \Delta^{1/2} \mathbf{\Sigma}_U(k) \Delta^{1/2}, 
\end{align*}
which reflects the additional scaling of the eigenvectors by the eigenvalues.

More recently, \citet{xie_entrywise_2024} extends this result to the regime $n\rho_n \asymp \log(n)$, building upon work in \citet{xie_efficient_2021} for sequences of deterministic matrices $P = \rho_n XX\t$, rather than stochastically generated community memberships.

\begin{theorem}[Restatement of Theorem~4.4 in \citet{xie_entrywise_2024} for SBMs] \label{thm:xieCLT2}
    Suppose the conditions of \cref{thm:xieCLT} hold. Write
    \begin{align*}
        \mathbf{\Sigma}_S^{(n)}(k)
        \coloneqq
        \Delta_n^{-1} \Gamma^{(k)}_n \Delta_n^{-1},
    \end{align*}
    where $\Gamma^{(k)}_n$ and $\Delta_n$ are defined in \cref{thm:xieCLT}. Suppose $\tau(i) = k$ with $n_k \asymp n_l$ for all $1 \le k, l \le n$. There exist sequences of orthogonal matrices $\tilde W$ and $W$ depending on $n$ and a constant $C > 0$ depending only on $d$ such that for any convex measurable set $\mathcal{A} \in \R^d$, it holds that
    \begin{align*}
        \bigg| \p\bigg\{ n^{1/2} \big(  \mathbf{\Sigma}_S^{(n)}(k) &\big)^{-1/2} \tilde W\t \big( (\hat X W)_{i*} - \rho_n^{1/2}\bm{\nu}_{k*} \big) \in \mathcal{A} \bigg\} - \p(G \in \mathcal{A}) \bigg\} \bigg|
        \leq
        \frac{C \log(n\rho_n)}{(n\rho_n)^{1/2}},
    \end{align*}
    where $G$ denotes a $d$-dimensional standard Gaussian random vector.
\end{theorem}

Taken together, \cref{thm:prdCLT,thm:xieCLT2} demonstrate a similar phenomenon as for the unscaled eigenvectors, namely that the estimated scaled eigenvectors are approximately Gaussian around block-specific centroids after appropriate (indefinite) orthogonal transformation. For the scaled eigenvectors, the cluster centroids do not depend on $n$ besides through the sparsity factor $\rho_n$, which, for practical purposes, could be estimated consistently in certain settings after imposing an identifiability condition, for example, $\max_{i} B_{ii} = 1$ or $\|B\| = 1$.

\subsubsection{Asymptotic Normality of the Laplacian Eigenvectors}
A central limit theorem for the scaled eigenvectors $\widecheck X$ of the Laplacian under the random dot product graph model is illustrated in the following theorem.

\begin{theorem}[Restatement of Corollary~3.3 in \citet{tang_limit_2018}] \label{thm:minhCLT2}
    Fix $K$ and $B = \bm{\nu} \bm{\nu}\t$. Assume that each $X_i$ is drawn independently according to $\p[X_i = \nu_{k*}] = \pi_k$, with $\nu_{k*}$ denoting the $k$-th row of $\nu$, with $\nu$ a $K \times K$ full rank matrix.   
    Suppose $A \sim P$ with $P = \rho_n X X\t= \rho_n ZBZ\t$. 
  Let $\check X$ be the spectral embedding of $\mathcal{L}(A) =\operatorname{diag}( A 1)^{-1/2} 
    \;A\;\operatorname{diag}( A 1)^{-1/2}$, given by the leading $K$ eigenvectors scaled by the square roots of the corresponding eigenvalues. Define
    \begin{align*}
        \tilde \Delta
        &\coloneqq
        \sum_{l=1}^K \pi_l \frac{\bm{\nu}_{l*} \bm{\nu}_{l*}\t}{ \sum_{l'} \pi_{l'} B_{ll'}}; \\
        \tilde \Gamma^{(k)}
        &\coloneqq
        \sum_{l=1}^K \pi_l \bigg( \frac{B_{kl} (1 - B_{kl})}{\sum_{l'} \pi_{l'} B_{kl'}} \bigg)
        \bigg( \frac{ \bm{\nu}_{l*}}{\sum_{l'} \pi_{l'} B_{ll'}} - \frac{\tilde \Delta\bm{\nu}_{k*}}{2 \sum_{l'} \pi_{l'} B_{kl'}} \bigg)
        \bigg( \frac{ \bm{\nu}_{l*}}{\sum_{l'} \pi_{l'} B_{ll'}} - \frac{\tilde \Delta\bm{\nu}_{k*}}{2 \sum_{l'} \pi_{l'} B_{kl'}} \bigg)\t , \\
        \tilde \Gamma^{(k)}_{o(1)}
        &\coloneqq
        \sum_{l=1}^K \pi_l \bigg( \frac{B_{kl}}{\sum_{l'} \pi_{l'} B_{kl'}} \bigg) 
        \bigg( \frac{ \bm{\nu}_{l*}}{\sum_{l'} \pi_{l'} B_{ll'}} - \frac{\tilde \Delta\bm{\nu}_{k*}}{2 \sum_{l'} \pi_{l'} B_{kl'}} \bigg)
        \bigg( \frac{ \bm{\nu}_{l*}}{\sum_{l'} \pi_{l'} B_{ll'}} - \frac{\tilde \Delta\bm{\nu}_{k*}}{2 \sum_{l'} \pi_{l'} B_{kl'}} \bigg)\t.
    \end{align*}
    Suppose that $n\rho_n \gg \log^4(n)$. Then, there exists a sequence of orthogonal matrices $W$ such that
    \begin{align*}
        \bigg( n \rho_n^{1/2} \bigg\{ \big( \check X W \big)_{i*} - \frac{\bm{\nu}_{k*}}{\sqrt{\sum_{l} n_l B_{lk}}} \bigg\} \bigg)
        \mid
        \tau(i) = k
        &\dto
        \begin{cases}
            \mathcal{N}\big(0, \mathbf{\Sigma}_{\mathcal{L}}(k) \big) 
            &
            \textnormal{~if~}\rho_n \equiv 1, \\
            \mathcal{N}\big(0, \mathbf{\Sigma}_{\mathcal{L},o(1)}(k) \big)
            &
            \textnormal{~if~}\rho_n \to 0, \end{cases}
    \end{align*}
    where the asymptotic covariance matrices are defined via
    \begin{align*}
        \mathbf{\Sigma}_{\mathcal{L}}(k)
        &\coloneqq
        \tilde \Delta\inv \tilde \Gamma^{(k)} \tilde \Delta\inv, \qquad 
        \mathbf{\Sigma}_{\mathcal{L},o(1)}(k)
        \coloneqq
        \tilde \Delta\inv \tilde \Gamma^{(k)}_{o(1)} \tilde \Delta\inv.
    \end{align*}
\end{theorem}

The above result has since been extended to the generalized random dot product graph in \citet{rubin-delanchy_statistical_2022} and to sequences of deterministic matrices $P = \rho_n XX\t$ in \citet{xie_efficient_2021}.

The denominator term $\sum_{l'} \pi_{l'} B_{ll'}$ in \cref{thm:minhCLT2} reflects the expected edge frequency for vertices in community $l$, and $\frac{\bm{\nu}_{k*}}{\sqrt{ \sum_{l} n_l B_{lk}}}$ represents a re-scaling of the latent vector associated with community $k$ that involves the (random) block sizes $\{n_{l}\}_{l=1}^{K}$, via accounting for nodal degree heterogeneity in the normalized Laplacian. The scaling to obtain asymptotic normality in \cref{thm:minhCLT2} is $n\rho_n^{1/2}$, demonstrating a similar shrinking-to-zero phenomenon as the unscaled eigenvectors of the adjacency matrices.

In \cref{thm:minhCLT2}, the covariance matrices $\mathbf{\Sigma}_{\mathcal{L}}(k)$ and $\mathbf{\Sigma}_{\mathcal{L},o(1)}(k)$ have more complicated expressions than those in \cref{thm:xieCLT} or \cref{thm:prdCLT}, though they are still fully determined by community membership probabilities $\pi_k$, the matrix $B$, and the membership of the $i$-th vertex. Given the differences between $\mathbf{\Sigma}_{S}(k)$ and $\mathbf{\Sigma}_{\mathcal{L}}(k)$, \citet{tang_limit_2018} considers model settings for which using the scaled eigenvectors of the Laplacian matrix is preferred to using the adjacency matrix for spectral clustering, on the basis of Chernoff information, an idea that was further explored in \citet{cape_spectral_2019}.

In \citet{modell_spectral_2021}, the authors consider the scaled eigenvectors $\hat X^+$ of the \emph{random walk Laplacian}, defined as
\begin{align*}
    \mathcal{L}_{\operatorname{rw}}(A)
    &\coloneqq
    D\inv A,
\end{align*}
where $D = \operatorname{diag}(A1)$ is the diagonal matrix of node degrees. Of note, $\mathcal{L}_{\operatorname{rw}}(A)$ has real eigenvalues since it is \emph{similar} to the positive semidefinite matrix $\mathcal{L}(A)$.  
\begin{theorem}[Restatement of Theorem~4.2 in \citet{modell_spectral_2021}]
    {Suppose the conditions of \cref{thm:prdCLT} hold}. 
    Define $\tilde \Delta$ as in \cref{thm:minhCLT2}, and let $\hat X^{+}$ denote the leading eigenvectors of $\mathcal{L}_{\operatorname{rw}}(A)$ scaled by the square roots of their corresponding eigenvalues. Denote
    \begin{align*}
        \tilde \Gamma^{(k)}_{\operatorname{rw}}  
        &\coloneqq
        \sum_{l=1}^K \pi_l \bigg( \frac{B_{kl} (1 - B_{kl})}{\sum_{l'} \pi_{l'} B_{kl'}} \bigg) 
        \bigg( \frac{ \bm{\nu}_{l*}}{\sum_{l'} \pi_{l'} B_{ll'}} - \frac{\tilde \Delta\bm{\nu}_{k*}}{ \sum_{l'} \pi_{l'} B_{kl'}} \bigg)
        \bigg( \frac{ \bm{\nu}_{l*}}{\sum_{l'} \pi_{l'} B_{ll'}} - \frac{\tilde \Delta\bm{\nu}_{k*}}{ \sum_{l'} \pi_{l'} B_{kl'}} \bigg)\t ,
        \\
        \tilde \Gamma^{(k)}_{\operatorname{rw},o(1)} 
        &\coloneqq
        \sum_{l=1}^K \pi_l \bigg( \frac{B_{kl}}{\sum_{l'} \pi_{l'} B_{kl'}} \bigg)
        \bigg( \frac{ \bm{\nu}_{l*}}{\sum_{l'} \pi_{l'} B_{ll'}} - \frac{\tilde \Delta\bm{\nu}_{k*}}{ \sum_{l'} \pi_{l'} B_{kl'}} \bigg)
        \bigg( \frac{ \bm{\nu}_{l*}}{\sum_{l'} \pi_{l'} B_{ll'}} - \frac{\tilde \Delta\bm{\nu}_{k*}}{ \sum_{l'} \pi_{l'} B_{kl'}} \bigg)\t.
    \end{align*}
    There exists an sequence of indefinite orthogonal matrices $Q$ and a constant $c \geq 1$ such that provided $n\rho_n \gg \log^{4c}(n)$, then
    \begin{align*}
        \bigg( n^{3/2} \rho_n \bigg\{ (\hat X^+ Q)_{i*} - \frac{\bm{\nu}_{k*}}{\rho_n^{1/2}\sum_{l} n_l B_{kl}} \bigg\} \bigg)
        \mid
        \tau(i) = k 
        &\dto
        \begin{cases}
            \mathcal{N}\big( 0, \mathbf{\Sigma}_{\operatorname{rw}}(k) \big)
            &
            \textnormal{~if~}\rho_n \equiv 1 ,\\
            \mathcal{N}\big( 0, \mathbf{\Sigma}_{\operatorname{rw},o(1)}(k) \big)
            &
            \textnormal{~if~}\rho_n \to 0,
        \end{cases}
    \end{align*}
    where the asymptotic covariance matrices are given by
    \begin{align*}
        \mathbf{\Sigma}_{\operatorname{rw}}(k)
        &\coloneqq
        \frac{\ipq \tilde \Delta\inv \tilde \Gamma_{\operatorname{rw}}(k) \tilde \Delta\inv \ipq}{\sum_{l'} \pi_{l'} B_{kl'}}, \qquad 
        \mathbf{\Sigma}_{\operatorname{rw},o(1)}(k)
        \coloneqq
        \frac{\ipq \tilde \Delta\inv \tilde \Gamma_{\operatorname{rw},o(1)}(k) \tilde \Delta\inv \ipq}{\sum_{l'} \pi_{l'} B_{kl'}}.
    \end{align*}
\end{theorem}

The performance of the random walk-based embedding for clustering under the degree-corrected stochastic blockmodel is studied in \citet{modell_spectral_2021}. The proofs therein use the tools developed in \citet{tang_limit_2018,rubin-delanchy_statistical_2022}, together with the fact that if $\check u$ is an eigenvector of $\mathcal{L}(A)$, then $D^{-1/2} \check u$ is an eigenvector of $\mathcal{L}_{\operatorname{rw}}(A)$. The limiting covariance matrix resembles that in \cref{thm:minhCLT2} because of the close relationship between the eigenvectors of $\mathcal{L}(A)$ and those of $\mathcal{L}_{\operatorname{rw}}(A)$.

To supplement the above discussion, \cref{table:asymptoticnormality} summarizes the literature discussed thus far.

\begin{center}
\begin{table*}[h!]
\begin{tabular} {|p{1.25in}||p{2in}|p{1in}|p{1.25in}|}
 \hline
 \textbf{Paper}
 &
 \textbf{Asymptotic normality}
 &
 \textbf{Sparsity}
 &
 \textbf{Model}\\
 \hline
 \citet{athreya_limit_2016}
 &
 Rows of scaled eigenvectors of $A$
 &
 $\rho_n \equiv 1$
 &
 $P = XX\t$, i.i.d. rows of $X$, distinct eigenvalues \\
 \hline
 \citet{cape_signal-plus-noise_2019}
 &
 Rows of unscaled eigenvectors of $A$
 &
 $n\rho_n \gg \log^4(n)$
 &
 $P = \rho_n XX\t$ \\
 \hline
 \citet{fan_asymptotic_2022}
 &
 Linear forms of distinct unscaled eigenvectors of $A$
 &
 $n\rho_n \gg  n^{\eta}$, $\eta > 0$
 &
 $P$ low-rank and distinct eigenvalues \\
 \hline
 \citet{li_two-sample_2018}
 &
 Squared Frobenius norm difference of scaled eigenvectors $\| (\hat UW - U ) \Lambda \|_{\operatorname{F}}^{2}$
 &
 $n\rho_n \gg n^{1/2}$
 &
 $P = ZBZ\t$, with balanced homogeneous SBM $B$\\
 \hline
 \citet{modell_spectral_2021}
 &
 Rows of $D^{-1/2} \check U$, where eigenvector matrix $\check U$ of the Laplacian and degree matrix $D$ 
 &
 $n\rho_n \gg \log^4(n)$
 &
 $P = \rho_n X \ipq X\t$, i.i.d. rows of $X$\\
 \hline
 \citet{rubin-delanchy_statistical_2022}
 &
 Rows of scaled eigenvectors of $A$ and the Laplacian matrix 
 &
 $n\rho_n \gg \log^4(n)$
 &
 $P = \rho_n X \ipq X\t$, i.i.d. rows of $X$\\
 \hline
 \citet{tang_limit_2018}
 &
 Rows of unscaled eigenvectors of the Laplacian matrix
 &
 $n\rho_n \gg \log^4(n)$
 &
 $P = \rho_n XX\t$, i.i.d. rows of $X$ \\
 \hline
 \citet{xie_efficient_2021}
 &
 Rows of scaled eigenvectors of $A$
 &
 $n\rho_n \gg \log^4(n)$ & $P = \rho_n XX\t$ \\
 \hline 
 (above)
 &
 Rows of one-step updated scaled eigenvectors (see \cref{sec:parameter})
 &
 $n\rho_n^5 \gg \log^2(n)$
 &
 $P = \rho_n XX\t$ \\
 \hline 
 (above)
 &
 Rows of scaled eigenvectors of the Laplacian matrix
 &
 $n\rho_n \gg \log^4(n)$
 & $P = \rho_n XX\t$
\\
 \hline 
 (above)
 &
 Rows of one-step updated scaled eigenvectors (see \cref{sec:parameter}) of the Laplacian matrix
 &
 $n\rho_n^4 \gg  \log^2(n)$
 &
 $P = \rho_n XX\t$ \\
 \hline
 \citet{xie_entrywise_2024}
 &
 Rows of scaled and unscaled eigenvectors of $A$
 &
 $n\rho_n \gtrsim \log(n)$ & $P = \rho_n XX\t$ \\
 \hline
 (above)
 &
 Rows of one-step updated scaled eigenvectors (see \cref{sec:parameter})
 &
 $n\rho_n \gg \log(n)$ & $P = \rho_n XX\t$ \\
 & & & \\
 \hline
\end{tabular}
\caption{Asymptotic normality results for scaled and unscaled eigenvectors}
\label{table:asymptoticnormality}
\end{table*}
\end{center}

\subsection{Relationships Between Limiting Results for Eigenvectors}\label{sec:relationship}
The results in \cref{sec:eigenvectors,sec:scaled_eigenvectors} illustrate a common theme, namely that each row vector of the embedding is (marginally) approximately Gaussian around one of $K$ different centroids in $\R^d$, where the location of the centroid depends on the particular embedding specification. In addition, aggregating across the embedding also yields asymptotic normality, per \cref{thm:fanCLT2,thm:liliCLT}. Underlying many of the proofs in the aforecited works, as well as the entrywise perturbation results discussed in \cref{sec:subspace_perturbation}, are the perturbation expansions
\begin{align*}
    \hat U W - U
    &=
    (A - \E[A])U \Lambda\inv + R_U
    \numberthis \label{asympexp} \\
    \hat U \hat \Lambda^{1/2} W' - U \Lambda^{1/2}
    &=
    (A - \E[A])U \Lambda^{-1/2} + R_X,
\end{align*}
where $R_U$ and $R_X$ are residual matrices whose $\ell_{2,\infty}$ norm is comparatively small with high probability for large graphs. Similar characterizations hold for Laplacian spectral embeddings for networks (e.g., Theorem~3.1 in \citet{tang_limit_2018} or Theorem of 3.1 of \citet{ke_optimal_2022}) and beyond network models (e.g., Lemma~4.18 in \citet{chen_spectral_2021}). The primary technical difficulty is that isolating leading-order terms requires care, though there are several different approaches for obtaining these and related expansions. For example, \citet{fan_asymptotic_2022} uses techniques based on contour integration and resolvents in random matrix theory, \citet{cape_signal-plus-noise_2019} uses deterministic matrix perturbations and series expansions, and \citet{xie_entrywise_2024} uses a \emph{leave-one-out technique}; see also \citet{abbe_entrywise_2020,chen_spectral_2021}.

The appearance of the orthogonal matrix $W$ stems from the fact that the eigenvalues $\Lambda$ need not always be distinct. If the eigenvalues are distinct, then Lemma~9 in \citet{fan_simple_2022} shows that a similar entrywise expansion holds for the individual eigenvectors (see also Theorem~3 in \citet{fan_asymptotic_2022}).

\begin{lemma}[Restatement of Lemma~9 in \citet{fan_simple_2022} for SBMs] \label{thm:simple}
    Fix $K$ and $B$. Suppose that $n\rho_n \gtrsim n^{\eta}$ for some constant $\eta > 0$, and suppose  that $n_k \asymp n_l$ for all $1\le k,l \le K$. Let $A \sim P$ with $P = \rho_n ZBZ\t$, where $P$ has $K$ distinct eigenvalues $\lambda_1, \dots, \lambda_K$ satisfying $\min_{i\neq j} |\lambda_i| / |\lambda_j| \geq 1 + c_0$ for some constant $c_{0} > 0$ for $\lambda_i \neq - \lambda_j$. Then, it holds that
    \begin{align*}
          t_k (e_i\t \hat U_{*k} - e_i\t U_{*k} )
          &=
          e_i\t (A - \E[A]) U_{*k}
          +
          o_{\p}( n^{-1/2}),
    \end{align*}
    where $t_k$ is the solution to the fixed-point equation in \cref{fixedpoint}.  
\end{lemma}

The expansions in \cref{asympexp} can be understood as holding for a diagonal (orthogonal) matrix $W$ with main diagonal elements $\pm 1$ when the eigenvalues are sufficiently separated. Explicitly, this leading-order expansion demonstrates that, under certain conditions, the entries of the $k$-th eigenvector have a variance that is asymptotically proportional to the $k$-th eigenvalue, a phenomenon that has also been shown in non-network settings (e.g.,~\citet{agterberg_entrywise_2022}). \cref{thm:simple} demonstrates that when the eigenvalues are distinct, then the asymptotic variance is proportional to the quantity $t_k^2$, while the difference between $\lambda_k^2/t_k^2$ and $1$ is asymptotically negligible for sufficiently dense networks. Additional details regarding $t_k$ are provided in \cref{sec:eigenvalues}.

The convergence of $\| (\hat{U} W - U)\Lambda\|_{\operatorname{F}}^{2}$ to a Gaussian random variable is established in \citet{li_two-sample_2018}. This expression can be viewed as a modification of the $\sin\Theta$ (Frobenius) distance between the eigenvectors, cf.~\cref{sec:eigenvectors}. Though there do not appear to be existing distributional results for $\sin\Theta$ distances customized to network problems, existing works do give concentration inequalities for the random variable $\|\hat U |\hat \Lambda|^{1/2} W - U |\Lambda|^{1/2}\|_{\operatorname{F}}$. Namely, see \citet{ tang_semiparametric_2017,tang_limit_2018,agterberg_nonparametric_2020} for applicable results to (generalized) random dot product graphs.

Asymptotic multivariate normality results by themselves do not address consistent estimation of the limiting covariance matrix for each row, which are needed to devise asymptotically valid confidence regions or provide data-driven uncertainty quantification. In \cref{sec:hyp_testing}, we review hypothesis testing settings where consistent estimation of the limiting covariance is possible.

\subsection{Cluster Structure: A Simple Numerical Example} \label{sec:example}
The asymptotic normality results presented in previous subsections have limiting covariance matrices that depend on the $B$ matrix and the community membership of the $i$-th vertex. The resulting cluster geometry reflects the particular form of these covariance matrices. Here, consider for example the setting where $B_{11} = B_{22} = a$, $B_{12} = B_{21} = b$, for $0 < b < a < 1$, and with equal size communities. In this setting, the leading eigenvector of $P = ZBZ\t$ is proportional to the constant vector of all ones. In addition, the second (un-normalized) leading eigenvector has entries $\pm 1$, with positive and negative signs distinguishing between the two communities. Consequently, since (unit-norm) eigenvectors are only unique up to sign, by symmetry it follows that $\big(\mathbf{\Sigma}_U(1) \big)_{\ell \ell} = \big(\mathbf{\Sigma}_U(2) \big)_{\ell \ell}$ for $\ell = 1$ or $2$, and that $\big(\mathbf{\Sigma}_U(1) \big)_{12} = -\big(\mathbf{\Sigma}_U(2) \big)_{12}$, where $\mathbf{\Sigma}_U(k)$ denotes the asymptotic covariance obtained from \cref{thm:xieCLT}. Consequently, the covariance matrices are of the form
\begin{align*}
    \mathbf{\Sigma}_U(1)
    =
    \begin{pmatrix}
        \alpha & \beta \\
        \beta & \gamma
    \end{pmatrix},
    \qquad 
    \mathbf{\Sigma}_U(2)
    =
    \begin{pmatrix}
        \phantom{-}\alpha & -\beta \\
        -\beta & \phantom{-}\gamma
    \end{pmatrix},
\end{align*}
and written in terms of underlying, implicit parameters $\alpha$, $\beta$, and $\gamma$. \cref{fig:uhat_whet_appetite} showcases this phenomenon, where the covariance for Community One is positive for dimensions one and two, whereas the covariance between dimensions one and two for Community Two is negative. However, this observation does not necessarily result in nearly spherical covariance matrices.

\begin{figure}
    \centering
    \includegraphics[height=.3\textheight,keepaspectratio]{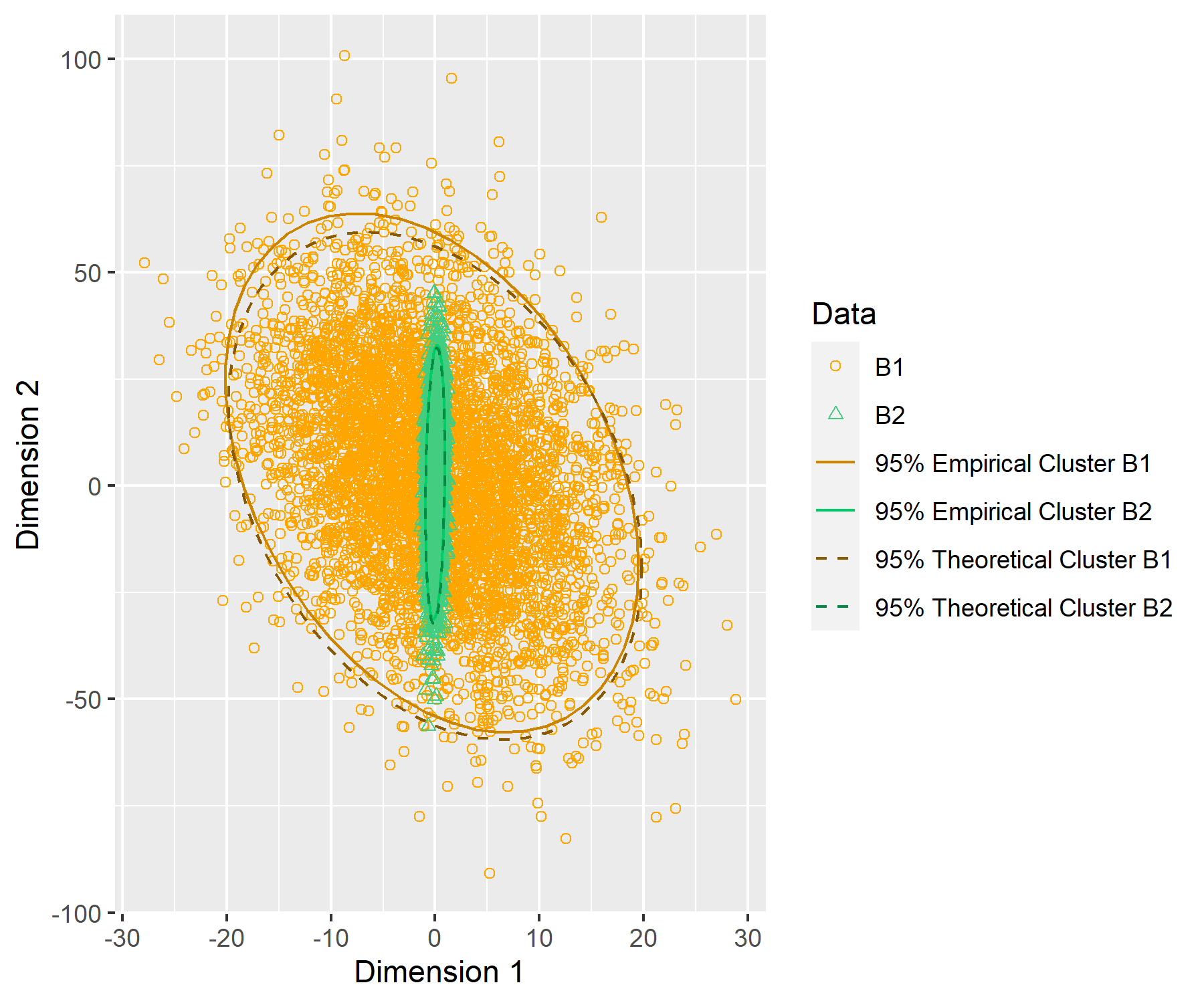}
    \caption{Comparison of the first $n/2$ rows of $n(\hat U W - U)$ for two different $B$ matrices with the same eigenvectors. The condition numbers of the asymptotic covariance matrices are $10.36$ and $1236.59$, respectively.}
    \label{fig:diffcovs}
\end{figure}

Consider two possible connectivity matrices
\begin{equation*}
    B^{(1)}
    =
    \begin{pmatrix}
        .02 & .01 \\
        .01 & .02
    \end{pmatrix},
    \qquad
    B^{(2)}
    =
    \begin{pmatrix}
        .9 & .85 \\
        .85 & .9
    \end{pmatrix},
\end{equation*}
with $n/2 = 5000$ vertices in each community. The leading eigenvectors of $P^{(1)} = Z B^{(1)} Z\t$ and $P^{(2)} = Z B^{(2)} Z\t$ are the same, but the resulting covariance matrices are different as they depend on the entries of $B$. In \cref{fig:diffcovs}, we plot the first $n/2$ rows of $n(\hat U W - U)$ for $B^{(1)}$ and $B^{(2)}$ respectively. Here, $B^{(1)}$ has a condition number of $3$, whereas $B^{(2)}$ has a condition number of $35$, resulting in covariance condition numbers of $10.36$ and $1236.59$, respectively. The cluster associated to $B^{(2)}$ is visibly elongated, even more than the matrix $B^{(1)}$ which consists of smaller connectivity probabilities. Computing the theoretical limiting covariance numerically yields
\begin{align*}
    \mathbf{\Sigma}_U(1)
    &=
    \begin{cases}
    \begin{pmatrix}
        65.556 & -64.667 \\
        -64.667 & 590
    \end{pmatrix}
    &
    \textnormal{~for~}
    B^{(1)}, \\
    \begin{pmatrix}
        0.142 & 0.857 \\
        0.857 & 174
    \end{pmatrix}
    &
    \textnormal{~for~}
    B^{(2)},
    \end{cases}
\end{align*}
which showcases the effect of ill-conditioned matrices $B$ on cluster structure. While the presented asymptotic normality results hold for fixed $B$ matrices as $n \to \infty$, the limiting covariance structure depends highly on the matrix $B$, and even sufficiently dense networks may still exhibit manifestly different cluster structure in latent space.

\subsection{Empirical Eigenvalues} \label{sec:eigenvalues}
The eigenvalues of large-dimensional matrices are well-studied in the field of random matrix theory. This survey focuses on asymptotic normality  concerning the top $d$ in magnitude eigenvalues of large adjacency matrices $A$. Among existing results in the literature, the following classical result establishes the limiting behavior of the leading eigenvalue of Erd\H{o}s--R\'{e}nyi random graphs.

\begin{theorem}[Restatement of Theorem~1 in \citet{furedi_eigenvalues_1981}] \label{thm:furedi}
    Fix $p \in (0,1)$. Let $A$ be a symmetric random matrix with i.i.d. upper triangular entries $A_{ij} \sim \operatorname{Bernoulli}(p)$, $1 \le i \le j \le n$. Let $\hat \lambda_{1}$ denote the largest eigenvalue of $A$. Then, as $n \to \infty$,
    \begin{align*}
        \hat \lambda_1 - np
        \dto
        \mathcal{N}(1-p, 2 p(1-p)).
    \end{align*}
    Observe that $np$ is the largest eigenvalue of $\E[A] = p 11\t$.
\end{theorem}

More recently, \citet{athreya_eigenvalues_2021} extended \cref{thm:furedi} to the general rank $d$ setting for stochastic blockmodels, as follows.

\begin{theorem}[Restatement of Theorem~3 in \citet{athreya_eigenvalues_2021} for SBMs] \label{athreyaeigenvalueCLT}
    Fix $K$ and $B = \bm{\nu} \ipq \bm{\nu}\t$. Suppose that $A \sim P$ with $P = ZBZ\t$ and that $\p[\tau(i) = k] = \pi_k$. Let $\mathbf{v}_i$ denote the unit norm eigenvector satisfying
    \begin{align*}
        \big( \bm{\nu}\t Z\t Z \bm{\nu} \big)^{1/2} \ipq \big( \bm{\nu}\t Z\t Z \bm{\nu} \big)^{1/2} \mathbf{v}_i
        =
        \lambda_i \mathbf{v}_i,
    \end{align*}
    where $\lambda_i$ is the $i$-th nonzero eigenvalue of $ZBZ\t$.  
    Define the vector $\mathbf{\tilde v}_i \coloneqq \big( \bm{\nu}\ Z\t Z \bm{\nu} \big)^{-1/2} \mathbf{v}_i$. Let $\bm{\eta}$ denote the $d$-dimensional vector with components
    \begin{align*}
        \bm{\eta}_i
        &=
        \frac{1}{\lambda_i} \sum_{l=1}^K n_l \mathbf{\tilde v}_i\t \bm{\nu}_l \bm{\nu}_l\t \mathbf{\tilde v}_i \sum_{k=1}^{K} B_{lk}(1- B_{lk}),
    \end{align*}
    and let $\mathbf{C} = (\mathbf{C}_{ij})$ denote the $d \times d$ matrix with elements
    \begin{align*}
        2 \bigg( \sum_{l=1}^K \sum_{k=1}^K n_k n_l \mathbf{\tilde v}_i\t \bm{\nu}_{l*} \bm{\nu}_{l*}\t \mathbf{\tilde v}_j  \mathbf{\tilde v}_i\t \bm{\nu}_{k*} \bm{\nu}_{k*}\t \mathbf{\tilde v}_j B_{kl}(1 - B_{kl}) \bigg).
    \end{align*}
    Then, as $n \rightarrow \infty$,
    \begin{align*}
        \mathbf{C}^{-1/2} \big( \hat \lambda_1 - \lambda_1 - \bm{\eta}_1, \hat \lambda_2 - \lambda_2 - \bm{\eta}_2, \dots, \hat \lambda_{d} - \lambda_d - \bm{\eta}_d \big)  
      &\dto
        \mathcal{N}(0, I_d).
    \end{align*}
\end{theorem}

In fact, the vector $\bm{\eta}$ is also shown to converge when the eigenvalues of $\sum_{k} \pi_k \mathbf{\nu}_{k*}\mathbf{\nu}_{k*}\t \ipq$ are assumed distinct (see Theorem~2 therein). Similar to \cref{thm:furedi}, it can be shown that $\bm{\eta}$ and $\mathbf{C}$ are bounded bias and covariance terms. Explicitly, $\sigma_i^2$ is the limiting value of the variance of the random variable $U_{*i}\t (A - \E[A]) U_{*i}$ after unconditioning on the latent vectors $X_{i*}$. Here $\lambda_i$ is order $n$ as the authors assume that $\rho_n = 1$. The results are derived under the general latent-space framework of the generalized random dot product graph described in \cref{sec:scaled_eigenvectors}.

A closely related result for individual eigenvalues is stated in the following theorem.

\begin{theorem}[Restatement of Theorem~1 in \citet{fan_asymptotic_2022} for SBMs] \label{faneigenvalueCLT}
    Fix $K$ and $B$. Suppose that $n\rho_n \gg n^{\eta}$ for some $\eta > 0$ and that $n_k \asymp n_l$ for all $1 \le k, l \le K$. Suppose $A \sim P$ with $P = \rho_n ZBZ\t$, where $P$ has $K$ distinct eigenvalues $\lambda_1, \dots, \lambda_K$ satisfying $\min_{i\neq j} |\lambda_i| / |\lambda_j| \geq 1 + c_0$ for some constant $c_0 > 0$ for $\lambda_i \neq - \lambda_j$. Then, as $n \rightarrow \infty$,
    \begin{align*}
        \frac{\lambda_k - t_k}{\operatorname{Var}\left[U_{*k}\t ( A - P) U_{*k}\right]^{1/2}}
        \dto
        \mathcal{N}(0,1),
    \end{align*}
    where $t_k$ is the unique solution to the fixed-point equation defined in \cref{fixedpoint}.
\end{theorem}

\cref{athreyaeigenvalueCLT,faneigenvalueCLT}, while related, are difficult to directly compare with one another. For one, unlike \citet{athreya_eigenvalues_2021}, \citet{fan_asymptotic_2022} does not consider a latent-space framework, so the eigenvalues of $n^{-1} P$ need not converge. Since $P = X \ipq X\t$, and the nonzero eigenvalues of a generic matrix product $MN$ are the same as the nonzero eigenvalues of the matrix product $NM$, the nonzero eigenvalues of $P$ are the same as those of $X\t X \ipq$. If each row of $X$ is drawn from a latent position distribution, then $n^{-1} X\t X \to \E[X_{1*} X_{1*}\t]$ almost surely, so the eigenvalues of $n^{-1} P$ are converging to the eigenvalues of $\E[X_{1*} X_{1*}\t] \ipq$. As a consequence, \cref{faneigenvalueCLT} allows for more sparsity by only specifying a sequence of matrices $P = P_n$. Here, $n\rho_n = \omega(n^{\eta})$ for some $\eta > 0$, which is a more stringent condition than the $\log(n)$ information-theoretic threshold required for consistency in community detection but is significantly sparser than the fully dense regime studied in \citet{athreya_eigenvalues_2021}. A further important distinction is that \cref{athreyaeigenvalueCLT} is a \emph{joint} convergence result when the eigenvalues are almost surely asymptotically distinct, whereas \cref{faneigenvalueCLT} holds marginally, for each individual eigenvalue. That being said, the result in \citet{fan_asymptotic_2022} allows for significantly more general forms of the low-rank probability matrix $P$.

Another impediment to explicitly comparing the above results is that, while the bias (centering) term in \citet{athreya_eigenvalues_2021} is computed directly, the bias term in \citet{fan_asymptotic_2022} is expressed implicitly through the solution of a complex analytic fixed-point equation, which is of order $n\rho_n$ from the discussion immediately following \cref{thm:fanCLT2}. On the other hand, the limiting variances are seen to be comparable since both stem from $U_{*k}\t (A - \E [A]) U_{*k}$ (which is of order $\rho_n$).

\cref{athreyaeigenvalueCLT,faneigenvalueCLT} both consider the fluctuations of the eigenvalues at the scaling $O(1)$. Under the latent space framework, the eigenvalues of $\mathcal{K}= n^{-1} X \ipq X\t$ can be viewed as the eigenvalues of an empirical kernel matrix with kernel $X_{i*}\t \ipq X_{j*}$. The results of \citet{koltchinskii_random_2000} imply that the eigenvalues of $\mathcal{K}$ are approximately Gaussian about the eigenvalues of a corresponding integral operator, which in this case are simply the eigenvalues of $\E[X_{1*}X_{1*}\t] \ipq$. Therefore, the fluctuations of the eigenvalues of $\mathcal{K}$ may be approximated by the eigenvalues of the normalized matrix $n^{-1} A$ or $(n\rho_n)^{-1} A$. These eigenvalues are studied in the context of a general low-rank random graph in \citet{lunde_subsampling_2023}.

\begin{theorem}[Restatement of Theorem~3.1 in \citet{lunde_subsampling_2023} for SBMs] \label{thm:lunde}
    Fix $K$ and $B$. Suppose that $n\rho_n \gg n^{1/2}$ and that $\p[\tau(i) = k] = \pi_k$. Let $A \sim P$ with $P = \rho_n ZBZ\t$. Suppose that the eigenvalues $\tilde\lambda_1, \dots, \tilde\lambda_d$ of $\mathbb{E}[X_{1*}X_{1*}\t] \ipq$ are distinct. Then, as $n \rightarrow \infty$,
    \begin{align*}
        n^{1/2} \bigg(\frac{\hat \lambda_1}{n\rho_n} - \tilde \lambda_1 , \dots, \frac{\hat \lambda_d}{n\rho_n} - \tilde \lambda_d \bigg)
        \dto
        \mathcal{N}(0,\mathbf{S}),
    \end{align*}
    where $\mathbf{S}$ denotes the asymptotic covariance matrix.
\end{theorem}

In \citet{lunde_subsampling_2023}, the probabilities are allowed to be drawn from a general latent position distribution but with $n\rho_n \gg n^{1/2}$. Therein, the covariance matrix $\mathbf{S}$ is not discussed in detail, though it is the same as in \citet{koltchinskii_random_2000}, due to the fact that the fluctuations of $(n\rho_n)^{-1} A - \mathcal{K}$ are of smaller order than the fluctuations of $\mathcal{K}$ and its corresponding integral operator. In other words, the leading $d$ eigenvalues of $(n\rho_n)^{-1} A$ and $(n\rho_n)^{-1} P$ have the same limiting distribution. In contrast, the results in \citet{athreya_eigenvalues_2021,fan_asymptotic_2022} concern the leading eigenvalues of $A$.  

The aforementioned results illustrate that for sufficiently dense graphs, the eigenvalues of $A$ are asymptotically biased estimators of the (population) eigenvalues of $P$, but, the eigenvalues of $(n\rho_n)^{-1} A$ are asymptotically unbiased estimators of those of the corresponding kernel operator (not necessarily unbiased for finite $n$). In all cases, these results hold under more general models, namely for both \cref{athreyaeigenvalueCLT,thm:lunde} the edge probability matrix $P$ can be drawn from a generalized random dot product graph, and in \cref{faneigenvalueCLT} the edge probability matrix need only be low-rank. Whether these results can be extended to various sparsity regimes or used explicitly for inference is a topic of ongoing investigation.

\subsection{Parameter Estimation}\label{sec:parameter}
When the community memberships are viewed as stochastic, and when both $\rho$ and $K$ are assumed known, then the model parameters are encapsulated in $(B,\pi)$. Absent additional structural assumptions, $B$ consists of $(K+1)K/2$ unconstrained parameters in $[0,1]$, while $\pi$ consists of $K$ nonnegative parameters in $[0,1]$ that sum to unity. Community memberships are unique only up to a global permutation of community labels (i.e.,~global relabeling), so strictly speaking, estimation of the matrix $B$ holds up to simultaneous row/column permutation, while estimation of the vector $\pi$ holds up to a permutation of its coordinates.

The likelihood function of $(A,\tau)$ and of $A$ when permitting self-edges are given by
\begin{align}
    \mathcal{L}(A, \tau; B, \pi) &= 
    \left( \prod_{1 \le i \le n} \pi_{\tau_{i}} \right)  \nonumber
    \times 
    \left( \prod_{1 \le i \le j \le n} (\rho B_{\tau_{i},\tau_{j}})^{A_{ij}}(1-\rho B_{\tau_{i},\tau_{j}})^{1-A_{ij}}\right),\\
    \mathcal{L}(A; B, \pi)
    &=
    \sum_{\tau \in [K]^{n}}\mathcal{L}(A, \tau; B, \pi). \label{eq:likelihood}
\end{align}
For stochastic blockmodel graphs $A$ governed by parameters $(B,\pi)$, where $\rho$ and $K$ are assumed known, the maximum likelihood estimate of $B$ and $\pi$ is given by
\begin{align}
    (\widehat{B}^{(\operatorname{MLE})}, \widehat{\pi}^{(\operatorname{MLE})})
    =
    \underset{B \in [0,1]^{K \times K}, \pi \in \mathcal{S}_{K-1}}{\operatorname{\textnormal{arg max}}} \mathcal{L}(A; B, \pi)
\end{align}
The asymptotic normality of the maximum likelihood estimate was obtained in Lemma~1 of \citet{bickel_asymptotic_2013} in the asymptotic regime $n\rho_{n} \gg \log n$, with implicit limiting variances (see \citet{tang_asymptotically_2022} for additional discussion).

\begin{theorem}[Restatement of Lemma~1 in \citet{bickel_asymptotic_2013}] \label{thm:BMLE}
    Fix $K$ and $B$. Suppose $A \sim P$ with $P = \rho_n ZBZ\t$ with $\p[\tau(i) = k] = \pi_k$. Let $\hat B^{(\operatorname{MLE})}$ denote the MLE of $B$ obtained from $A$ with $\rho_n$ assumed known. If $n\rho_n\gg \log(n)$, then as $n \rightarrow \infty$, 
    \begin{equation*}
        n \rho_n^{1/2}
        \big(\hat B^{(\operatorname{MLE})} - B \big)
        \hspace{-.5em}
        \dto
        \hspace{-.5em}
        \begin{cases}
            \mathcal{N}(0,\mathbf{\Sigma}_B)
            &
            \textnormal{~if~} \rho_n \equiv 1, \\
            \mathcal{N}(0,\mathbf{\Sigma}_{B,o(1)})
            &
            \textnormal{~if~} \rho_n \to 0,
        \end{cases}
    \end{equation*}
    where 
    \begin{align*}
        (\mathbf{\Sigma}_B)_{kl}
        &=
        \begin{cases}
            \frac{2 B_{kk}(1 - B_{kk})}{\pi_k^2}
            &
            \textnormal{~if~}k = l, \\
            \frac{ B_{kl}(1 - B_{kl})}{\pi_k\pi_l}
            &
            \textnormal{~if~}k \neq l. 
        \end{cases}
        \numberthis
        \label{eq:N_MLE_dense} \\
        (\mathbf{\Sigma}_{B,o(1)})_{kl}
        &= 
       \begin{cases}
            \frac{2 B_{kk}}{\pi_k^2}
            &
            \textnormal{~if~}k = l, \\
            \frac{ B_{kl}}{\pi_k\pi_l}
            &
            \textnormal{~if~}k \neq l.
       \end{cases} 
       \numberthis
       \label{eq:N_MLE_sparse}
    \end{align*}
      In both cases, $\{ n\rho_n^{1/2}~( \hat B^{(\operatorname{MLE})}_{kl} - B_{kl} )\}_{k \leq l}$ are asymptotically independent.
\end{theorem}

The maximum likelihood estimator computed using the likelihood  \cref{eq:likelihood} can be computationally expensive and is in general NP-hard to compute \citep{amini_pseudo-likelihood_2013}. Furthermore, as discussed in \cref{sec:asymp_norm}, when specifying a stochastic blockmodel $B$ matrix having dimension $K \times K$, there is no \emph{a priori} reason to assume that $B$ has full matrix rank equal to $K$. This observation is important, since the presence or absence of a rank constraint can influence the number of free parameters in the SBM parameter estimation problem. For $d \coloneqq \textrm{rank}(B) \le K$, by taking the spectral decomposition of $B$, there are at most $d(2K - d + 1)/2 \le K(K+1)/2$ total parameters, via the eigendecomposition of $B$.

\begin{example}[SBM with two communities and rank one $B$ matrix]
    Consider a $K=2$ SBM with
    $B =
    \left(
    \begin{smallmatrix}
        p^{2} & pq \\
        pq & q^{2}
    \end{smallmatrix}\right)$,
    where $p,q \in (0,1), p \neq q$. The matrix $B$ can be written as the outer product of the vector $(p,q)$ with itself, hence $\textrm{rank}(B)=1$. Here, there are two parameters, rather than $K(K+1)/2 = 3$ parameters, namely $p$ and $q$, with $B_{12}$ depending on $B_{11}$ and $B_{22}$.
\end{example}

The paper \citet{tang_asymptotically_2022} considers entrywise estimation of $B$ by via spectral estimators $\wh{B}_{k\ell}^{(S)}$ formed by averaging over the truncated eigendecomposition $\wh{U}\wh{\Lambda}\wh{U}^{\top}$ of the adjacency matrix $A$ in the following manner. Assume that $d = \textrm{rank}(B)$, $K$, and $\rho_{n}$ are all known. Then, clustering the rows of the eigenvector matrix $\wh{U}$ via $K$-means (or Gaussian mixture modeling) as in \cref{alg:spectral_clustering} yields an estimated block membership vector $\wh{\tau}$ and corresponding assignment vectors $\wh{s}_{k} \in \{0,1\}^{n}$ and block sizes $\wh{n}_{k}$ for all $k \in [K]$, where
\begin{equation*}
    \wh{B}_{k\ell}^{(S)}
    \coloneqq
    \frac{1}{\wh{n}_{k}\wh{n}_{\ell} \rho_{n}}\wh{s}_{k}^{\top}\wh{U}\wh{\Lambda}\wh{U}^{\top}\wh{s}_{\ell}.
    \numberthis
    \label{eq:spectralB}
\end{equation*}
The following result demonstrates the entrywise asymptotic normality of $\hat B^{(S)}$ when $B$ is rank $K$.

\begin{theorem}[Restatement of Corollary~2 in \citet{tang_asymptotically_2022}] \label{thm:Bspectral}
    Fix $K$ and $B$ where $\textrm{rank}(B) = K$. Suppose that $A \sim P$ with $P = \rho_n ZBZ\t$ with $\p[\tau(i) = k] = \pi_k$. Let $\hat B^{(S)}$ be defined as in \cref{eq:spectralB} with $\rho_n$ assumed known. If $n\rho_n \gg n^{1/2}$, then as $n \rightarrow \infty$,
    \begin{align*}
        n \big( \hat B^{(\operatorname{S})}_{kl} - B_{kl} \big)
        \hspace{-.5em}
        \dto
        \hspace{-.5em}
        \begin{cases}
            \mathcal{N}(0,(\mathbf{\Sigma}_{B})_{kl})
            &
            \textnormal{~if~}\rho_n \equiv 1, \\
            \mathcal{N}(0,(\mathbf{\Sigma}_{B,o(1)})_{kl})
            &
            \textnormal{~if~}\rho_n \to 0,
        \end{cases}
        \end{align*}
        for $\mathbf{\Sigma}_B$ and $\mathbf{\Sigma}_{B,o(1)}$ as in \cref{thm:BMLE}.
 \end{theorem}

Consequently, when $n\rho_n \gg n^{1/2}$ and $B$ is full rank, the spectral estimators are (marginally) asymptotically efficient, meaning that the variances in \cref{thm:BMLE} coincide with the variances in \cref{thm:Bspectral}.
More is true, namely Remark~1 in \citet{tang_asymptotically_2022} shows that the half-vectorization of $B^{(\operatorname{S})}$ is asymptotically multivariate (jointly) normal, with given bias, variance, and covariance terms.

When $B$ is not full rank, then the spectral estimator is no longer asymptotically efficient. In particular, \citet{tang_asymptotically_2022} also establishes that $\hat B^{(S)}$ is biased and has an inflated asymptotic variance when $B$ is rank degenerate, providing explicit expressions for the bias and variance terms. To ameliorate the efficiency gap, \citet{tang_asymptotically_2022} further develops a one-step update to obtain efficiency in parameter estimation. Using the (debiased) spectral embedding-based estimate as an initial estimate, the one-step estimator is obtained by performing a single Newton--Raphson step in the negative gradient direction with respect to the likelihood of $B$ (considering $\tau$ fixed), and \citet{tang_asymptotically_2022} demonstrates that this procedure yields an asymptotically efficient estimator. A similar result was also obtained in \citet{xie_euclidean_2023} for an alternative parametrization and different initial estimator to yield asymptotic efficiency in a modified parameter space. The results therein are in terms of an underlying $d(2K-d+1)/2$-dimensional parameter.

In \citet{xie_efficient_2021} the authors consider a notion of \emph{local efficiency} with respect to the random dot product graph model which includes stochastic blockmodels with positive semidefinite $B$ matrices as special cases.

\begin{theorem}[Restatement of Theorems~2 and 5 in \citet{xie_efficient_2021} ] \label{thm:mleclt}
    Fix $K$ and $B$. Suppose that $n_k \asymp n_l$ with $n_k/n \to \pi_k \in(0,1)$ for all $1\le k,l \le n$.  Let $A \sim P$ with $P = \rho_n ZBZ\t$. Suppose that there exists a constant $\delta > 0$ such that $B = \bm{\nu} \bm{\nu}\t$ satisfies $\delta \leq B_{kl} \leq 1 - \delta$ for all $k,l$. Let $P = \rho_n XX\t$, and suppose that $n_k \asymp n_l$ with $n_k/n \to \pi_k \in (0,1)$ for all $1 \le k, l \le n$. Define the set
    \begin{align*}
        \mathcal{X}(\delta)
        \coloneqq
        \left\{
        x:
        x\t u \in [\delta,1-\delta]
        \text{ for all $u \in \mathcal{X}(\delta)$} 
        \right\}.
    \end{align*}
    Define
    \begin{align*}
        \hat X_{i*}^{(\operatorname{MLE})}
        &\coloneqq
        \argmax_{x \in \mathcal{X}(\delta)} \sum_{j\neq i} \{ A_{ij}\log( x\t X_{j*}) +
        (1 - A_{ij}) \log(1 - x\t X_{j*}) \}; \\
        \hat X_{i*}^{(\operatorname{OS})}
        &\coloneqq
        \hat X_{i*} + \bigg\{ \frac{1}{n} \sum_{j=1}^n \frac{ \hat X_{j*} \hat X_{j*}\t}{\hat X_{i*}\t \hat X_{j*}(1 - \hat X_{i*}\t \hat X_{j*})} \bigg\}\inv 
        \bigg\{ \frac{1}{n} \sum_{j=1}^n \frac{(A_{ij} - \hat X_{i*}\t \hat X_{j*}) \hat X_{j*}}{\hat X_{i*}\t \hat X_{j*}(1 - \hat X_{i*}\t \hat X_{j*})} \bigg\},
    \end{align*}
    where $\hat X_{i*}^{(\operatorname{MLE})}$ is computed assuming $\{X_{j*}\}_{j\neq i}$ are known, and $\hat X$ are the scaled eigenvectors of $A$.
    Define
    \begin{align*}
        \mathbf{G}(k)
        \coloneqq
        \sum_{l=1}^{K} \frac{\pi_l}{B_{kl}(1- B_{kl})} \bm{\nu}_{l*} \bm{\nu}_{l*}\t. 
    \end{align*}  
    Then, as $n \to \infty$, if $\rho_n \equiv 1$, 
    \begin{align*}
        &\bigg( n^{1/2} \big( \hat X_{i*}^{(MLE)} - \bm{\nu}_{k*} \big) \bigg)
        \mid
        \tau(i) = k \dto
        \mathcal{N}\big(0,\mathbf{G}(k)\inv\big).
    \end{align*}
    Similarly, if $n\rho_n^5 \gg\log^2(n)$, then there exists a sequence of orthogonal matrices $W$ such that
    \begin{align*}
        \bigg( n^{1/2} \big( ( &\hat X^{(\operatorname{OS})} W )_{i*} - \rho_n^{1/2} \bm{\nu}_{k*} \big) \bigg)
        \mid
        \tau(i) = k
        \dto
       \mathcal{N}(0,\mathbf{G}(k)\inv).
    \end{align*}
    Finally, define $\Sigma_{S}(k) \coloneqq \lim_{n \to \infty} \mathbf{\Sigma}_S^{(n)}(k)$, with $\mathbf{\Sigma}_S^{(n)}(k)$ as in \cref{thm:xieCLT2} (which exists under the present assumptions). Then, 
   $\mathbf{\Sigma}_S(k) - \mathbf{G}(k)\inv$ is positive semidefinite for all $x \in \mathcal{X}(\delta)$.
\end{theorem}

The estimator $\hat X_{i*}^{(\operatorname{MLE})}$ is obtained by maximizing the random dot product graph likelihood assuming that $X_{j*}$ is known for $j \neq i$, which is not possible in practice. On the other hand, the vector $\hat X_{i*}^{(\operatorname{OS})}$ is obtained by starting with $\hat X$ as an initial estimate and then performing a Newton--Raphson step update for $\hat X_{i*}$ with respect to the random dot product graph likelihood. In both cases, these procedures yield estimators  for $\bm{\nu}_{k*}$, modulo orthogonal transformation, that achieve the limiting variance $\mathbf{G}(k)\inv$ which is no larger than the asymptotic variance of $\hat X_{i*}$ with respect to the Loewner order. For the particular rank one stochastic blockmodel with $B$ matrix
\begin{align*}
    B
    \coloneqq
    \begin{pmatrix}
        p^2 & pq \\
        pq & q^2
    \end{pmatrix},
\end{align*}
the authors show that the one-step covariance in the first community is strictly smaller than that of the covariance in \cref{eq:cltrdpg} unless $q = (1-p^2)/p$. However, these results require stronger sparsity assumptions than other results, and it is not clear how the asymptotic covariance matrix is affected with either indefinite stochastic blockmodels or considering the update with respect to the SBM likelihood function, as opposed to the RDPG likelihood function. The above results were later extended in \citet{xie_entrywise_2024} to the setting $n\rho_n \gg \log(n)$.

\subsection{Subgraph Counts}
For network data, the probabilistic concept of \emph{moments} can be understood in the language of \emph{subgraph counts}, namely counts and frequencies of specified subgraphs appearing in a network. The study of the so-called \emph{method of moments} for network data was initiated in \citet{bickel_method_2011} with numerous follow-on works in recent years. Considerable effort in this direction has centered around the latent position random network, defined as follows.

\begin{definition}[Latent position network -- undirected, loopy]\label{def:latent_position}
    We say that an adjacency matrix $A$ is a latent position random network on $n$ vertices with sparsity term $\rho = \rho_n \in (0,1]$ if, conditional on $\{X_{i\star}\}_{i=1}^{n}$, its upper triangular entries are generated independently in the manner
   \begin{align*}
        A_{ij}
        =
        \begin{cases}
            \operatorname{Bernoulli}\left(\rho \cdot f \left( X_{i*}, X_{j*}\right) \right),
            &
            \textnormal{~for~}i \leq j, \\
            A_{ji}
            &
            \textnormal{~for~}i > j,
        \end{cases}
    \end{align*}
    where $X_{i*} \in \mathbb{R}^d$ for all $1 \le i \le n$, and $f: \mathbb{R}^d \times \mathbb{R}^d \to [0,1]$ satisfies $f(x,y) = f(y,x)$. The variables $\{X_{i*}\}_{i=1}^{n}$ may be deterministic or modeled in a generative fashion, e.g., i.i.d. $X_{i*} \sim F_X$ for a distribution $F_X$ on $\R^d$.
\end{definition}

The latent position network with $d = 1$ and $F_X$ being the uniform distribution on $[0,1]$ is often referred to as a so-called \emph{graphon}, though in some settings the function $f$ is itself referred to as a graphon. When $f$ is the (pseudo) inner product, one arrives at the (generalized) random dot product graph. \cref{def:latent_position} therefore encompasses a broad family of models beyond the low-rank network models discussed in previous sections. A latent position network can be associated to an infinite-dimensional version of a generalized random dot product graph via its so-called \emph{graph root distribution} as studied in \citet{lei_network_2021}.  

Next, we highlight the method of moments for networks when, for simplicity, $\rho_n$ is known. Suppose $G$ is a subgraph of interest with $r$ vertices, that is, $G \subset \{(i,j): 1 \leq i < j \leq n\}$ with $|G| = r$. Let $A(G)$ denote the subgraph of $A$ induced by the vertices in $G$. The corresponding network moment are defined as
\begin{align*}
    U_n
    &\coloneqq
    \rho_n^{-r} \p\left[ E( A(G)) = G \right]
    \\
    &=
    \mathbb{E} \left\{ \prod_{( i,j) \in G} f( X_i,X_j) \prod_{(i,j) \in \tilde G} (1 - \rho_n f(X_i,X_j)) \right\},
\end{align*}
where $E(\cdot)$ denotes the edge set of a graph, and $\tilde G$ is the set of edges that do not appear in $G$ but are incident to vertices in $G$.
Define the (normalized) empirical network moment as 
\begin{align*}
    \hat U_n
    &\coloneqq
    \frac{1}{\rho_n^{r} N(G)} {n \choose r}\inv \sum_{H \subset A} \mathbb{I}_{\{H \sim G\}},
\end{align*}
where the notation $H \sim G$ indicates that the graph $H$ is isomorphic to $G$, the summation is over all subgraphs $H$ of $A$ of size $r$, and $N(G)$ denotes the number of graphs isomorphic to $G$. The asymptotic normality of $\hat U_n$ was obtained in Theorem~1 of \citet{bickel_method_2011}.

\begin{theorem}[Theorem~1 in \citet{bickel_method_2011}]
    Suppose that $\int_{0}^{1} \int_{0}^{1} f^2(u,v) \operatorname{d}v \operatorname{d}u < \infty$. Suppose that $G$ is fixed with $\int_{0}^{1} \int_{0}^{1} f^{2r}(u,v) \operatorname{d}u \operatorname{d}v < \infty$. If $G$ is acyclic and $n\rho_n \to C \in (0,\infty]$, then
    \begin{align*}
        \sqrt{n} \big(\hat U_n - U_n \big)
        \dto
        \mathcal{N}(0,\sigma^2(G)), \numberthis \label{uconv}
    \end{align*}
    where $\sigma^2(G)$ depends on $C$ if $ C < \infty$, and otherwise does not. If instead $G$ is cyclic and  $n\rho_n \gg n^{1 - 2/r}$, then \cref{uconv} continues to hold.  
   \end{theorem}

The asymptotic variance, $\sigma^2(G)$, can be difficult to compute and in general depends on the variance of the associated $U$-statistic
\begin{align*}
    \tilde U_n
    &\coloneqq
    \frac{1}{N(G)}
    {n \choose r}\inv
    \sum_{H \sim G} 
    \Bigg\{
    \prod_{(i,j) \in H} f(X_i,X_j) 
   \prod_{(i,j) \in \tilde H} (1 - \rho_n f(X_i,X_j))
    \Bigg\}
    \\
    &=
    \E\big(\hat U_n
    \mid
    \{X_i\}_{i=1}^{n} \big).
\end{align*}
When $n\rho_n$ diverges, $\sigma^2(G)$ is equal to the asymptotic variance of $\sqrt{n}(\tilde U_n - U_n)$, implying that the contribution of the Bernoulli noise of the network $\sqrt{n}(\hat U_n -\tilde U_n)$ is asymptotically negligible relative to the contribution of $\sqrt{n}(\tilde U_n - U_n)$.

In the time since \citet{bickel_method_2011}, a number of papers have studied the theoretical properties of network moments for bootstrapping \citep{bhattacharyya_subsampling_2015,levin_bootstrapping_2019,lunde_subsampling_2023,lin_theoretical_2020,green_bootstrapping_2022,lin_trading_2022,zhang_edgeworth_2022} and for hypothesis testing \citep{gao_testing_2017,maugis_testing_2020,maugis_central_2023,mitra_inference_2024,shao_higher-order_2024}. 
For example, \citet{zhang_edgeworth_2022} develops Edgeworth expansions for network moments, showing that the additional Bernoulli noise contributes a so-called smoothing effect to the higher-order accuracy of the Gaussian distribution. \citet{levin_bootstrapping_2019} considers low-rank networks and shows that by leveraging the low-rank property it is possible to develop a consistent bootstrap procedure for subgraph counts.

\section{Hypothesis Testing and Multiple-Network Normality} \label{sec:hyp_testing}
Network hypothesis testing takes many forms depending on how the statistical unit of interest is specified. Units of interest can include, for example, individual (or several) vertices, representing microscale (local) structure, or individual (or several) entire networks, representing macroscale (global) structure. Each scale or level of resolution has different practical applications and interpretations, and, depending on the inference task at hand, the asymptotic normality results presented in the previous section may be adaptable to designing and analyzing test statistics tailored to inference problems. This section focuses primarily on test statistics motivated by the asymptotic normality results presented in the previous sections.

\subsection{Vertex-Level Tests} \label{sec:vertextests}
One problem of interest is testing whether two vertices share certain properties such as latent community memberships. In this vein, \citet{fan_simple_2022} considers the problem of testing whether two vertices from a general degree-corrected mixed-membership model have the same membership parameters, where the null hypothesis for nodes $i, j$ is
\begin{align*}
    H_0: Z_{i*} = Z_{j*}.
\end{align*}
The following vertex-level testing result is obtained.
\begin{theorem}[Restatement of Theorems~1 and~2 in \citet{fan_simple_2022}] \label{thm:fansimple}
    Fix $K$ and $B$, where $B$ has full rank. Suppose that $n\rho_n \geq n^{\eta}$ for some $\eta > 0$ and that $n_k \asymp n_l$ for all $1\le k,l,\le K$. Let $A \sim P$ with $P = \rho_n ZBZ\t$, where $P$ has $K$ distinct eigenvalues $\lambda_1, \dots, \lambda_K$ satisfying $\min_{i\neq j} |\lambda_i| / |\lambda_j| \geq 1 + c_0$ for some constant $c_0 > 0$ for $\lambda_i \neq - \lambda_j$. Suppose also that the eigenvalues of $(\rho_n)^{-1} \operatorname{Cov}(( e_i - e_j)\t (A - P) U)$ are bounded away from zero and infinity. Let
    \begin{align*}
        \mathbf{\Sigma}_{U}(i,j)
        &\coloneqq
        \operatorname{Cov}\big( (e_i - e_j)\t (A - P) U \Lambda\inv \big).
    \end{align*}
    Define
    \begin{align*}
        T^2_{U}
        &\coloneqq
        \big( \hat U_{i*} - \hat U_{j*} \big)\t \mathbf{\Sigma}_{U}\inv(i,j) \big( \hat U_{i*} - \hat U_{j*} \big).
    \end{align*}
    Then, under the null hypothesis $H_0: Z_{i*} = Z_{j*}$ with $\tau(i) = k$, as $n \rightarrow \infty$, it holds that
    \begin{align*}
        T^2_U
        \dto
        \chi^2_K,
    \end{align*}
    where $\chi^2_K$ denotes the chi-squared distribution with $K$ degrees of freedom. In addition, under the fixed alternative $Z_{i*} \neq Z_{j*}$, for any constant $C > 0$, as $n \rightarrow \infty$, it holds that
    \begin{align*}
        \p\left[ T^2_U > C\right]
        \to
        1.
    \end{align*}
    Finally, these results continue to hold for the plug-in test statistic
    \begin{align*}
        \hat T^2_{U}
        &\coloneqq
        (\hat U_{i*} - \hat U_{j*})\t  \mathbf{\hat{\Sigma}}_{U} \inv(i,j)(\hat U_{i*} - \hat U_{j*}),
    \end{align*}
    provided $\mathbf{\hat{\Sigma}}_U(i,j)$ satisfies
    \begin{align*}
        \rho_n\inv \| \Lambda( \mathbf{\hat{\Sigma}}_U(i,j) - \mathbf{\Sigma}_{U}(i,j) ) \Lambda \|
        \to
        0
    \end{align*}
    in probability as $n \to \infty$.
\end{theorem}

The proofs of the distributional convergence of $T^2_U$ and $\hat T^2_U$ depend on the asymptotics for eigenvectors and eigenvalues from \citet{fan_asymptotic_2022} as discussed in \cref{sec:eigenvectors,sec:eigenvalues}. The test statistic $T^2_U$ and its empirical estimate $\hat T^2_U$ can be understood as the analogue of the classical Hotelling $T^2$ statistic, considering the rows of $U$ as the means, and considering $\mathbf{\Sigma}_U(i,j)$ as the asymptotic covariance of the $i$-th row and $j$-th row of the leading eigenvectors discussed in \cref{sec:asymp_norm}. In fact, under the additional conditions in \cref{thm:xieCLT} the limiting matrix $\mathbf{\Sigma}_U(k) = \lim_{n\to \infty} \mathbf{\Sigma}_U^{(n)}(k)$ satisfies
\begin{align*}
    \lim_{n\to\infty}\mathbf{\Sigma}_U(i,j)
    &=
    \mathbf{\Sigma}_U(\tau(i))) + \mathbf{\Sigma}_U(\tau(j)).
\end{align*}

\cref{thm:fansimple} does not explicitly define an estimator $\mathbf{\hat{\Sigma}}_U(i,j)$. The sample eigenvalues are biased as shown in \cite{fan_asymptotic_2022} (see also \citet{athreya_eigenvalues_2021} and the requisite discussion in \cref{sec:eigenvalues}), hence their estimator is based on a debiasing procedure for the eigenvalues that can be viewed as a single fixed-point iteration for the fixed-point equation that determines proper centering for the empirical eigenvalues from \citet{fan_asymptotic_2022}. See \cref{fixedpoint}.

The aforecited paper also prove a local power result for mixed-membership blockmodels, allowing $Z_{i*} - Z_{j*}$ to converge to zero, and it also extends the result to account for degree corrections using the normalization procedure first proposed in \citet{jin_fast_2015}.

In \citet{du_hypothesis_2023}, the authors propose using the scaled eigenvectors to perform a similar test. Under the latent position framework, they study a hypothesis test for the more general hypothesis
\begin{align*}
    H_0:
    X_i
    =
    X_j
    .
\end{align*}

\begin{theorem}[Restatement of Theorem~3.2 in \citet{du_hypothesis_2023}] \label{thm:duclt}
    Fix $K$ and $B$, and suppose that $A \sim P$ with $P = \rho_n ZBZ\t$. Suppose that $n\rho_n \gg \log(n)$ and that $n_k \asymp n_l$ for all $1 \le k,l,\leq K$.  Suppose that there exists a constant $c > 0$ such that $B_{kl} \geq c$ for all $k$ and $l$, and that $B$ has $p$ positive and $q$ negative eigenvalues with $p + q = d$. Define
    \begin{align*}
        \hat T^2_S
        \coloneqq
        (\hat X_{i*} - \hat X_{j*})\t \mathbf{\hat \Sigma}\inv_{S}(i,j)(\hat X_{i*} - \hat X_{j*}),
    \end{align*}
    where
    \begin{align*}
        \mathbf{\hat{\Sigma}}_{S}(i,j)
        &=
        \mathbf{\hat \Sigma}(\hat X_{i*}) + \mathbf{\hat \Sigma}(\hat X_{j*}), \\
        \mathbf{\hat \Sigma}(\hat X_{i*})
        &\coloneqq 
        \hat \Lambda\inv\bigg[ \sum_{k=1}^{n} \hat X_{k*} \hat X_{k*}\t \hat X_{i*}\t \ipq \hat X_{k*} \big(1 - \hat X_{i*}\t \ipq \hat X_{k*} \big) \bigg] \hat \Lambda\inv.
    \end{align*}
    Then, under the null hypothesis, as $n \rightarrow \infty$,
    \begin{align*}
        \hat T^2_S
        \dto
        \chi^2_d,
    \end{align*}
    and under the fixed alternative, for any choice of constant $C > 0$, as $n \rightarrow \infty$, it holds that
    \begin{align*}
        \p\left[\hat T^2_S >  C \right]
        \to
        1.
    \end{align*}
\end{theorem}

The authors also prove a local power result for $\|X_{i*} - X_{j*}\| \to 0$ and provide analogous results for the Laplacian spectral embedding and for degree-corrected stochastic blockmodels. For the latter, they normalize each row by its Euclidean norm instead of the first component of the leading eigenvector. Furthermore, \citet{du_hypothesis_2023} extends the central limit theorem for $\hat X$ in \citet{rubin-delanchy_statistical_2022} (\cref{thm:prdCLT}) to deterministic sequences $P = \rho_n X \ipq X\t$ as opposed to rows $X_i$ drawn from a latent position distribution.

The matrix $\mathbf{\hat{\Sigma}}_{S}(i,j)$ appearing in \cref{thm:duclt} can be viewed as the plug-in estimate of $\mathbf{\Sigma}_{S}(\tau(i)) + \mathbf{\Sigma}_{S}(\tau(j))$ stemming from the limiting results of \citet{rubin-delanchy_statistical_2022}. Similarly, the limiting covariance $\mathbf{\Sigma}_U(i,j)$ in \cref{thm:fansimple} corresponds to the asymptotic covariance matrix stemming from the distributional theory for the eigenvectors (e.g.,~\cref{thm:simple} or \cref{thm:xieCLT}). Moreover, estimating $\mathbf{\Sigma}_U(i,j)$ requires an additional step to debias the empirical eigenvalues, whereas in \citet{du_hypothesis_2023}, the authors show that estimating $\mathbf{\Sigma}_{S}(i,j)$ requires no additional debiasing procedure for the eigenvalues. Informally, they show that the matrix $\hat X = \hat U |\hat \Lambda|^{1/2}$ already contains information about the spectrum via the eigenvalue scaling, so there is no further need to incorporate additional information to estimate the covariance. Finally, in \citet{du_hypothesis_2023}, the authors allow for repeated eigenvalues of the matrix $P = \rho_n X \ipq X\t$, and sparsity of order $n\rho_n = \omega(\log(n))$, whereas in \citet{fan_simple_2022} the authors require that $P = \operatorname{diag}(\theta) Z B Z\t \operatorname{diag}(\theta)$ has distinct eigenvalues and that the sparsity regime satisfies $n\rho_n = \omega(n^{\eta})$ for some $\eta > 0$.

Theory for two-vertex testing was later extended in \citet{xie_entrywise_2024} to the setting that $n\rho_n \gtrsim \log(n)$, and then further modified to use the one-step updated spectral embeddings as proposed in \citet{xie_efficient_2021}. Explicitly, using the one-step updated term, the author proposes a similar Hotelling-type $T^2$ test statistic using the asymptotic variance $\mathbf{G}(k)\inv$ per \cref{sec:parameter}. Since $\mathbf{G}(k)\inv$ is no larger than the limiting variance $\mathbf{\Sigma}_{S}(k)$, \citet{xie_entrywise_2024} shows that the test statistic using the one-step updated embedding is asymptotically more powerful than the test statistic using $\hat T^2_S$.

\subsection{One-Sample and Two-Sample Network Tests} \label{sec:graphtests}
Network-level tests can be designed to determine whether the network itself or its expected adjacency matrix exhibits certain properties. For low-rank SBM networks, one natural goodness-of-fit consideration is the null hypothesis
\begin{align*}
    H_0^{(\operatorname{gof})}:
    \operatorname{rank}(B)
    =
    d_0
\end{align*}
for a fixed positive integer $d_{0}$, against the alternative $H_A: \operatorname{rank}(B) > d_0$. The asymptotic normality results presented in \cref{sec:asymp_norm} all require that $\operatorname{rank}(B) = d$ is known and fixed (or can be consistently estimated). \citet{han_universal_2023}, proposes a test statistic for the rank of the matrix of interest under generic noise model assumptions that include stochastic blockmodels. Similar to \citet{fan_simple_2022}, the limiting distribution of the test statistic is closely related to the asymptotics considered in \citet{fan_asymptotic_2022} which further emphasizes the need for deriving limit theory for graph embeddings to study test statistics. Their test statistic is motivated as follows. If the underlying matrix $P$ is known, then, the matrix $E \coloneqq A - \E[A]$ is mean zero and has independent entries. If self-edges are permitted, then the diagonals of $E$ constitute $n$ independent random draws from a centered Bernoulli distribution, and therefore a rescaled test statistic using the diagonal entries of $E$ would yield asymptotic normality. In practice, one does not observe the matrix $P$ so it needs to be estimated. \citet{han_universal_2023} proposes taking the truncated eigendecomposition of $A$ and estimating $E$ as the residual. More specifically, define
\begin{align*}
    \hat P^{(d_0)}
    \coloneqq
    \hat U_{[1:d_0]*} \hat \Lambda_{[1:d_0]} \hat U_{[1:d_0]*}\t,
\end{align*}
that is, $\hat U_{[1:d_0]*}$ and $\hat \Lambda_{[1:d_0]}$ consist of the top $d_0$ eigenvectors and eigenvalues, respectively. From the preceding discussion, under the null hypothesis that $d = d_0$, the residual $\hat E^{(d)} \coloneqq A - \hat P^{(d)}$ is anticipated to behave approximately like the mean-zero matrix $E$, and test statistics can be constructed using the entries of $\hat E^{(d)}$. Moreover, under the alternative $d_0 < d$, the matrix $\hat E^{(d_0)}$ should contain some information about the matrix $P$ and will not necessarily be approximately mean zero.

A technical challenge is that the entries of the matrix $\hat E^{(d)}$ are not independent, and therefore the classical central limit theorem does not apply by simply summing up the entries of $\hat E^{(d)}$. The approach in \citet{han_universal_2023} is to subsample the edges of $\hat E^{(d)}$ in an independent fashion. In particular, consider the test statistic
\begin{align}
    T_{n,d_0,m}
    \coloneqq
    \frac{\sqrt{m} \sum_{i\neq j} \hat E^{(d_0)}_{ij} Y_{ij}}{\sqrt{2 \sum_{i\neq j} \left(\hat E^{(d_0)}_{ij}\right)^2}},
    \label{tndm}
\end{align}
where $Y_{ij}$ are Bernoulli random variables satisfying $\p[Y_{ij} = 1] = m^{-1}$ for $i < j$ and $Y_{ij} = Y_{ji}$ otherwise. The following result establishes the asymptotic distribution of $T_{n,d_0,m}$.

\begin{theorem}[Restatement of Theorems~3.1 and~3.2 in \citet{han_universal_2023}] \label{thm:rirs}
    Fix $K$ and $B$. Suppose that $n\rho_n \gtrsim n^{\eta}$ for some $\eta > 0$ and that $n_k \asymp n_l$ for $1\le k,l,\leq K$. Let $A\sim P$ with $P = \rho_n ZBZ\t$ where $P$ has $K$ distinct eigenvalues $\lambda_1, \dots, \lambda_K$ satisfying $\min_{i\neq j} |\lambda_i| / |\lambda_j| \geq 1 + c_0$ for some constant $c_0 > 0$ for $\lambda_i \neq - \lambda_j$. Further suppose that
    \begin{align*}
        n^{1-\eta/2}
        \ll
        m
        \ll
        \frac{n^{1+\eta}}{\log(n)}.
    \end{align*}
    Under the null hypothesis of $\operatorname{rank}(B) = d_0$, as $n \rightarrow \infty$,
    \begin{align*}
        T_{n,d_0,m}
        \dto
        \mathcal{N}(0,1).
    \end{align*}
    If the one-sided alternative hypothesis $d > d_0$ holds and
    \begin{align*}
        \left|\sum_{i\neq j} \sum_{k= d_0+1}^d \lambda_k U_{ik} U_{jk} \right|
        \ll
        n\rho_n,
    \end{align*}
    then as $n \to \infty$,
    \begin{align*}
           T_{n,d_0,m} - \frac{\sqrt{m} \sum_{i\neq j}\sum_{k=d_0+1}^{d} \lambda_k U_{ik} U_{jk} Y_{ij}}{\sqrt{2 \sum_{i\neq j}  \left(\hat E^{(d_0)}_{ij}\right)^2}} 
           &\dto
           \mathcal{N}(0,1).
    \end{align*}
    Finally, if instead
    \begin{align*}
        \left| \sum_{i\neq j} \sum_{k=d_0+1}^d \lambda_k U_{ik} U_{jk} \right|
        \gg
        \sqrt{m} \left( \sqrt{n\rho_n} + \sum_{d_0=1}^{d} |\lambda_k| \right),
    \end{align*}
    then $\p\left[| T_{n,d_0,m}| > C\right] \to 1$ for any pre-specified constant $C > 0$.
\end{theorem}

The authors also determine more general conditions on $m$ such that asymptotic normality holds. The centering in \cref{thm:rirs} reflects the fact that $\hat E^{(d_0)}$ contains information from the low-rank matrix when $d_0 < d$, as exhibited by the appearance of the eigenvalues $\lambda_k$ for $d_0 + 1 \leq k \leq d$. Specialized to the setting of full-rank stochastic blockmodels, their results enable accurately estimating the number of communities $K$, namely the sequential testing procedure $d_0 = 1, 2, \dots, d_{\max}$ asymptotically rejects the null hypothesis at a fixed level $\alpha$ unless $d = d_0$.

The paper \citet{han_universal_2023} tests whether the rank of $P$ is $d_0$ but does not include a goodness-of-fit test for specific submodels such as testing for stochastic blockmodelness. Earlier, \citet{lei_goodness--fit_2016} proposes a test specifically designed for stochastic blockmodels by estimating the matrix $E$ using the estimated communities and then computing the largest eigenvalue of an appropriately rescaled version. The associated limiting null distribution is shown to be Tracy--Widom, though \citet{han_universal_2023} points out that this test empirically performs poorly for sparser graphs. Similarly, \citet{wang_likelihood-based_2017} studies a likelihood-based approach for model selection, though their focus is on estimation (not hypothesis testing).   

The related work \citet{zhang_adjusted_2023} proposes a goodness-of-fit test for degree-corrected stochastic blockmodels under the assumption that entries of $A$ are Poisson random variables. This work leverages high-dimensional asymptotics for testing multinomials to establish that their proposed  test statistic converges in distribution to a Gaussian random variable under the null hypothesis that there are $K$ communities. Again, while the considerations there somewhat outside the scope of this survey, the analysis therein showcases how network asymptotics can differ significantly from classical asymptotics, as demonstrated through the interplay between sample size ($n$), sparsity $(\rho_n)$, and 
other model parameters including $K$.

Via \citet{tang_asymptotically_2022} or \citet{xie_euclidean_2023}, it is also possible to devise single-network tests for specific stochastic blockmodel hypotheses. Consider, for example, the setting where the block memberships are known. To test the null hypothesis $B = B^{(0)}$ for some fixed full rank matrix $B^{(0)}$, the results in \citet{tang_asymptotically_2022} imply that the spectral estimator for $B$ is asymptotically efficient. Moreover, provided the asymptotic covariance can be estimated consistently, a simultaneous confidence interval for the vectorization of $B$ will be asymptotically valid. 

In classical two-sample testing, the user tests whether two sets of observations plausibly arise from the same distribution or have similar distributional properties. By analogy, two-graph testing seeks to determine whether two graphs have similar network-centric distributional properties. One version of the two-graph hypothesis test assumes that there is a known \emph{a priori} correspondence between vertices, which is equivalent to testing whether $P^{(1)} = P^{(2)}$, where $P^{(1)}$ and $P^{(2)}$ denote the probability matrices for two observed adjacency matrices $A^{(1)}$ and $A^{(2)}$, respectively. If the correspondence is not known, then this test is equivalent to testing whether $P^{(1)} = \Pi P^{(2)} \Pi\t$, where $\Pi$ is a permutation matrix. If one assumes the underlying $P$ matrix comes from a stochastic blockmodel with random community assignments, one can test whether the two graphs have the same underlying distributions; this corresponds to testing the hypothesis $(B^{(1)}, \pi^{(1)}) = (B^{(2)}, \pi^{(2})$.

In many application settings, multiple graphs are vertex-matched, namely graphs are observed on the same vertex set with known vertex correspondences. For example, in time series of graphs, one might observe connectivity of the same vertex set evolving over time. Given a specific vertex of interest, its connectivity profile or community membership may change over time. When $B^{(1)} = \gamma B^{(2)}$ for some constant $\gamma$, \citet{li_two-sample_2018} considers testing whether the community memberships between two networks are equal, with null hypothesis
\begin{align*}
    H_0: Z^{(1)} = Z^{(2)}.
\end{align*}
Their test statistic can be viewed as a re-weighted $\sin\Theta$ distance between the estimated eigenvectors.

\begin{theorem}[Restatement of Theorem~5 and Theorem~7 in \citet{li_two-sample_2018}] \label{thm:lilitwosample}  Fix $K$.  
    Suppose that $A^{(1)} \sim P^{(1)}$ and $A^{(2)} \sim P^{(2)}$ with $P^{(1)} = \rho_n Z^{(1)}B(Z^{(1)})\t$ and $P^{(2)} = \rho_n \gamma Z^{(2)}B(Z^{(2)})\t$, where $B$ has the form 
    \begin{align*}
        B_{kl}
        &=
        \begin{cases}
            a
            &
            \textnormal{~if~}k = l, \\
            b
            &
            \textnormal{~if~}k =l,
        \end{cases}
    \end{align*}
    with constant $a$, $b$, and $\gamma$. Suppose that $n\rho_n \geq n^{1/2 + \eps}$ for some $\eps > 0$. Let $\hat U^{(1)}$ be the eigenvectors of $A^{(1)}$ and let $\hat U^{(2)}$ be the eigenvectors of $A^{(2)}$. Let $\hat \Lambda_2$ be the eigenvectors of $A^{(2)}$. Define
    \begin{align*}
        T_{n,K}
        &=
        \frac{1}{Kn} \| (\hat U^{(1)} W_{1,2} - \hat U^{(2)} ) \hat \Lambda_2 \|_{\operatorname{F}}^2,
    \end{align*}
    where $W_{1,2} = \argmin_{W \in \mathbb{O}(d)} \|\hat U^{(1)} W - \hat U^{(2)} \|_{\operatorname{F}}$ is the Frobenius optimal orthogonal matrix aligning $\hat U^{(1)}$ and $\hat U^{(2)}$. Then, under the null hypothesis $H_0: Z^{(1)} = Z^{(2)}$, and provided each community size is order $n/K$, in the large-sample limit,
    \begin{align*}
        \frac{T_{n,K} - \mu_n}{\sigma_n}
        \dto
        \mathcal{N}(0,1),
    \end{align*}
    where $\mu_n$ and $\sigma_n$ have the form
     \begin{align*}
        \mu_n
        &=
        \gamma^2 \rho_n b(1- \rho_n b) + \rho_n b(1-\rho_n b) +
        \frac{1}{K} \bigg( \gamma \rho_n a ( 1 - \rho_n a) + \rho_n a(1-\rho_n a) -
        \big( \gamma^2 \rho_n b(1-\rho_n b) + \rho_n b(1-\rho_n b) \big) \bigg), \\
        \sigma^2_n
        &=
        \frac{2}{nK} \bigg( \big( \gamma^2 \rho_n b(1-\rho_n b) + \rho_n b(1-\rho_n b) \big)^2 \\
        &\qquad\qquad+
        \frac{\big(\gamma^2 \rho_n a(1-\rho_n a) + \rho_n a(1-\rho_n a) \big)^2}{K} -
        \frac{\big( \gamma^2\rho_n b(1-\rho_n b) + \rho_n b(1-\rho_n b) \big)^2}{K} \bigg)+
        o(1).
    \end{align*}
   Let $q_{\alpha/2}$ and $q_{1 - \alpha/2}$ denote quantiles of the standard Gaussian distribution. Under the local alternative $\| Z^{(1)} - Z^{(2)} \|_0 \geq \frac{C K}{n^{1- \eps'} \sqrt{\mu_n}}$, with exactly $n/K$ vertices in each community, it holds that 
    \begin{align*}
        \p\bigg[ q_{\alpha/2} < \frac{T_{n,K} - \mu_n}{\sigma_n} < q_{1- \alpha/2} \bigg]
        \to
        1.
    \end{align*}
\end{theorem}
The result in \citet{li_two-sample_2018} allows $a$ and $b$ to change with $n$ provided $\frac{a-b}{a} \gtrsim \frac{K\log(n)}{n}$, reflecting the setting in which perfect community recovery is information-theoretically possible \citep{abbe_community_2018}. The stated asymptotic normality continues to hold when replacing the parameters with suitable estimators, namely $\mu_n \mapsto \hat \mu_n$ and $\sigma_n \mapsto \hat \sigma_n$. While \citet{li_two-sample_2018} does not explicitly provide such estimators, it is mentioned that they exist using previously established results such as the methodology in \citet{tang_asymptotically_2022}.

The proof of \cref{thm:lilitwosample} applies the analysis leading to \cref{thm:liliCLT}, demonstrating that $\|(\hat U^{(1)} W_1 - U^{(1)})\Lambda^{(1)}\|_{\operatorname{F}}^2$ is asymptotically Gaussian after appropriate normalization and scaling. This argument exploits the particular structure of $B$ together with the approximation
\begin{align*}
    \| (\hat U^{(1)} W- U^{(1)} ) \Lambda^{(1)} \|_{\operatorname{F}}^2
    &\approx
    \| (A^{(1)} - P^{(1)} ) U^{(1)} \|_{\operatorname{F}}^2,
\end{align*}
which follows from the discussion in \cref{sec:relationship}. See also \cref{asympexp}.

The paper \citet{tang_semiparametric_2017} proposes a similar test statistic to \citet{li_two-sample_2018} based on Procrustes alignment of the scaled eigenvectors of each graph and shows that their test statistic concentrates about a constant quantity determined by the distributions of each graph. This constant can also be viewed as the analogue of $\mu_n$ in the limit theorem for \citet{li_two-sample_2018}. \citet{tang_semiparametric_2017} does not study the explicit limiting distribution of their test statistic, and their main results depend primarily on a concentration inequality for the Frobenius norm difference between the scaled eigenvectors. See also \citet{agterberg_nonparametric_2020,athreya_numerical_2023} for similar concentration results.

Several other tests have been proposed for two graphs \citep{agterberg_nonparametric_2020,ghoshdastidar_two-sample_2017,ghoshdastidar_two-sample_2020,tang_semiparametric_2017,tang_nonparametric_2017}, though either they do not have associated asymptotic normality guarantees or are not focused on stochastic blockmodels.

\subsection{Multiple-Graph Asymptotics} \label{sec:twographtests}
The two-sample testing results in \cref{thm:lilitwosample} are closely related to the asymptotic theory derived for a single network in \cref{thm:liliCLT}. Analogously, testing procedures relying on similar theory will typically require first embedding each individual network separately and then aligning them, say via Procrustes rotations. In many practical situations networks may share some structure such as communities but may contain some network-level differences such as edge probabilities, so a \emph{joint graph embedding} may be desirable, namely, an embedding that includes both shared information from each network and individual vertex information within each network.

Assuming each graph has the same   probability matrix, \citet{levin_central_2019} studies a joint embedding constructed from the scaled eigenvectors of the so-called \emph{omnibus} matrix. Given $L$ independent vertex-matched graphs, the \emph{omnibus matrix} $\hat{M}_{\operatorname{OMNI}}$ is defined as
\begin{align*}
  \hat{M}_{\operatorname{OMNI}}
  &\coloneqq
    \begin{pmatrix}
        A^{(1)}
        &
        \frac{1}{2}(A^{(1)} + A^{(2)})
        &
        \cdots
        &
        \frac{1}{2}(A^{(1)} + A^{(L)}) \\
        \frac{1}{2}(A^{(1)} + A^{(2)})
        &
        A^{(2)}
        &
        \cdots
        &
        \frac{1}{2}(A^{(2)} + A^{(L)}) \\
        \vdots
        &
        \vdots
        &
        \ddots
        &
        \vdots \\
        \frac{1}{2}(A^{(1)} + A^{(L)})
        &
        \frac{1}{2}(A^{(2)} + A^{(L)})
        &
        \cdots
        &
        A^{(L)}
    \end{pmatrix}.
\end{align*}
Obtaining an embedding for $\hat{M}_{\operatorname{OMNI}}$ does not require additional Procrustes alignment that would otherwise be needed if each graph were separately embedded. Here, \citet{levin_central_2019} derives asymptotic normality guarantees for the scaled eigenvectors of $\hat{M}_{\operatorname{OMNI}}$ when all the {block} probability matrices are the same.
\begin{theorem}[Restatement of Theorem~1 in \citet{levin_central_2019}]
    \label{thm:omni}
    Fix $K,L$ and $B = \bm{\nu\nu}\t$. Let $A^{(l)} \sim P^{(l)}$, with $P^{(l)} \equiv ZBZ\t$ with $\p[\tau(i) = k] = \pi_k$. Let $\hat U_{\operatorname{OMNI}}$ and $\hat \Lambda_{\operatorname{OMNI}}$ denote the $Ln\times d$ and $d\times d$ matrices of leading eigenvectors and eigenvalues of $\hat{M}_{\operatorname{OMNI}}$ respectively. There exists a sequence of orthogonal matrices $W_*$ such that as $n \to \infty$,
    \begin{align*}
      \bigg( n^{1/2} \big((\hat U_{\operatorname{OMNI}} \hat \Lambda_{\operatorname{OMNI}}^{1/2} W_*)_{(n-1)l+i*} - \bm{\nu}_{k*} \big)\bigg)
      \mid
      \tau(i) = k
           &\dto
      \mathcal{N}(0,\tfrac{L+3}{4L}\mathbf{\Sigma}_{S}(k)),
    \end{align*}
    where $\mathbf{\Sigma}_{S}(k)$ is defined in \cref{thm:prdCLT} with $q = 0$.  
\end{theorem}

The coefficient $\frac{L+3}{4L}$ equals one when $L = 1$, and so the resulting asymptotic covariance matrix here matches the asymptotic covariance matrix for the rows of $\hat X$ for a single network, as in \cref{sec:asymp_norm}. More generally, \citet{pantazis_importance_2022} shows that this coefficient is determined by the choice of weights assigned to the diagonal and off-diagonal block matrices. For example, if, instead of weighting each block-diagonal element by the average of each corresponding matrix (i.e.,~$\frac{1}{2}$ weighting to both $A^{(i)}$ and $A^{(j)}$), one simply uses the average $\frac{1}{L} \sum_{l=1}^{L} A^{(l)}$, then the limiting variance is instead scaled by $\frac{L^2 + L - 1}{L^3} \asymp \frac{1}{L}$, reflecting the fact that variability shrinks as more networks $L$ are observed for this particular choice of weighting. See \citet{pantazis_importance_2022} for further discussion on how the choice of weighting induces correlation between rows of the associated embedding.

\cref{thm:omni} requires $P^{(l)} = ZBZ\t$ for all $l$, which can be restrictive. The asymptotic normality of the omnibus embedding and associated bias was further studied in \citet{draves_bias-variance_2021} when $P_{ij}^{(l)} = X_{i*}\t C^{(l)} X_{j*}$, where $C^{(l)}$ is a diagonal matrix with positive entries, giving the following result.

\begin{theorem}[Restatement of Theorem~1 in \citet{draves_bias-variance_2021} for Stochastic Blockmodels]
    Fix $K$ and the $K \times K$ matrix $\nu$. Let $A^{(l)} \sim P^{(l)}$ with $P^{(l)} = ZB^{(l)} Z\t$ with $B^{(l)} = \bm{\nu} C^{(l)} \bm{\nu}\t$ for some diagonal matrix $C^{(l)}$ with nonnegative diagonal elements, and suppose that $\min_{i} \max_{l} C^{(l)}_{ii} > 0$, with $A^{(l)}$ independent of $A^{(l')}$. 
    Suppose also that $\p[\tau(i) = k] = \pi_k$ and suppose that $\sum_{k} \pi_k\bm{\nu}_{k*} \bm{\nu}_{k*}\t$ is diagonal and full rank. Let $\hat U_{\operatorname{OMNI}}$ and $\hat \Lambda_{\operatorname{OMNI}}$ be as in \cref{thm:omni}. There exists an orthogonal matrix $W_*$ such that with high probability it holds that
     \begin{align*}
        \bigg\{
        \bigg(\big( \hat U_{\operatorname{OMNI}} \hat \Lambda_{\operatorname{OMNI}}^{1/2} W_*\big)_{n(l-1)+i*} - \sqrt{C^{(l)}}\bm{\nu}_{k*} \bigg)
        \mid
        \tau(i) = k
        \bigg\}
       &=
        \big( \tilde C^{(l)} - \sqrt{C^{(l)}} \big) \bm{\nu}_{k*} + N_{n(l-1)+i},
    \end{align*}
    with $\tilde C^{(l)} = \frac{1}{2} \big[ C^{(l)} C^{-1/4} + C^{1/4} \big]$, $C = \frac{1}{L} \sum_{l=1}^{L} (C^{(l)})^2$, and $N_{n(l-1)+i}$ satisfies
    \begin{align*}
        \| N \|_{2,\infty}
        &\lesssim
        \frac{L^{3/2} \log(nL)}{n^{1/2}}.
    \end{align*}
\end{theorem}

The above decomposition resembles the discussion in \cref{sec:relationship}, though there is an additional term of the form $\big( \tilde C^{(l)} - \sqrt{C^{(l)}}\big) \bm{\nu}_{k*}$ appearing on the right hand side, and the term $N$ behaves as a first-order error, similar to the bounds in \cref{prdtwoinfty}. Succinctly, the rows of the scaled eigenvectors of $\hat{M}_{\operatorname{OMNI}}$ are biased away from $\sqrt{C^{(l)}}\bm{\nu}_{k*}$, with bias depending on the additional scaling factor $C^{(l)}$ associated to the $l$-th network. The authors also prove a normality result for the rows of $\hat U_{\operatorname{OMNI}} \hat \Lambda_{\operatorname{OMNI}}^{1/2}$ using this decomposition, via additional expansions for the concentration of the noisy bias about its mean, and the expansion of the term $N_{(l-1) + i}$.

For the purpose of jointly clustering based on multiple graphs, observing more networks improves inference and clustering by boosting signal. Consider the setting where $A^{(1)}, \dots, A^{(L)}$ all satisfy $\E[A^{(l)}] = \rho_n Z B Z\t$ for some $\rho_n$, $Z$, and $B$ not depending on $l$. Let $\bar A \coloneqq \frac{1}{L} \sum_l A^{(l)}$, and let $\hat U_{\operatorname{ave}}$ denote the leading $d$ eigenvectors of $\bar A$. By an application of the $\ell_{2,\infty}$ bound in \citet{chen_spectral_2021} (see \cref{thm:twoinfty}), the concentration of the leading eigenspace improves by a factor of $L^{-1/2}$. Similarly, using the proof methods mentioned in previous sections, asymptotic normality guarantees can be derived for this averaged matrix as both the number of vertices and the number of graphs tends to infinity, thereby providing an initial justification for using more networks in jointly clustering vertices.

The argument in the previous paragraph requires that the networks have the same block probability matrices. 
The \emph{multilayer stochastic blockmodel} posits that between graphs (layers) the memberships remain constant while the probabilities of community membership are allowed to vary between graphs. Various algorithms have been proposed for clustering the networks jointly based on spectral methods \citep{chen_global_2022,lei_bias-adjusted_2022,paul_consistent_2016,paul_spectral_2020} or tensor decompositions \citep{jing_community_2021,lei_consistent_2020,lyu_latent_2023}, though the focus of these is often on weak recovery.

Pertaining to asymptotic normality, \citet{jones_multilayer_2021,arroyo_inference_2021} study algorithms for embeddings in specific multiple graph models that respectively subsume particular instantiations of the multilayer stochastic blockmodel. In \citet{arroyo_inference_2021} the authors study an algorithm for subspace recovery under the assumption that each graph probability matrix shares a common subspace, in the sense that it can be factored as $U S^{(l)} U\t$, where $S^{(l)}$ is a so-called score matrix corresponding to the $l$-th graph, and $U$ is an orthonormal matrix representing the common subspace. For a multilayer stochastic blockmodel, one may take $U = Z(Z\t Z)^{-1/2}$ with $S^{(l)} = (Z\t Z)^{1/2} B^{(l)} (Z\t Z)^{1/2}$.

Let $\hat U^{(l)}$ denote the matrix of leading $d$ eigenvectors of $A^{(l)}$, and define $\hat U_{\operatorname{MASE}}$ to be the leading $d$ eigenvectors of the matrix $\sum_{l=1}^{L} \hat U^{(l)} (\hat U^{(l)})\t$. Note that $\hat U_{\operatorname{MASE}}$ can equivalently be expressed as the left singular vectors of the matrix $[ \hat U^{(1)}, \dots, \hat U^{(L)} ]$. Define
\begin{align*}
    \hat S^{(l)}
    &\coloneqq
    \hat U_{\operatorname{MASE}}\t A^{(l)} \hat U_{\operatorname{MASE}}.
\end{align*}
The authors prove the following asymptotic normality result for their estimates $\hat S^{(l)}$.

\begin{theorem}[Restatement of Theorem~11 in \citet{arroyo_inference_2021} for MLSBMs] \label{thm:MASE1}
    Fix $K$ and $L$. Suppose that $n\rho_n \gg \log(n)$ and that $n_k \asymp n_l$ for all $1\leq k,l\leq K$. Suppose that $A^{(l)} \sim P^{(l)}$ with $P^{(l)} = \rho_n ZB^{(l)} Z\t$, with $B^{(l)}$ full rank and fixed. Define
    \begin{align*}
        \Sigma^{(l)}_{\frac{2k+q(q-1)}{2}, \frac{2k + q(q-1)}{2}} =
        \sum_{s=1}^{n-1} \sum_{t=s+1}^{n} P_{st}^{(l)}(1 - P_{st}^{(l)}) 
        \big[ (U_{sk}U_{tq} + U_{tk}U_{sq})( U_{sk'}U_{tq'} + U_{sk'}U_{tq'}) \big].
    \end{align*}
    Let $\hat S^{(l)}$ be defined as above. 
    There exists a sequence of orthogonal matrices $W$ such that, as $n \to \infty$,
    \begin{align*}
       \bigg(  \Sigma^{(l)}_{\frac{2k+q(q-1)}{2}, \frac{2k' + q'(q'-1)}{2}} \bigg)^{-1/2} \big( W \hat S^{(l)} W\t - S^{(l)} + H^{(l)} )_{kq} 
      &\dto
       \mathcal{N}(0,1),
    \end{align*}

    where $H^{(l)}$ is a stochastic offset term which satisfies $\E \|H^{(l)}\|_{\operatorname{F}} = O(d L^{-1/2})$.
\end{theorem}
The authors study the vectorization of $\hat S^{(l)}$ when the matrix $\Sigma^{(l)}$ is assumed to have its smallest eigenvalue bounded away from zero. \cref{thm:MASE1} can be compared to, for example \cref{faneigenvalueCLT} or \cref{athreyaeigenvalueCLT}, which demonstrate that the empirical eigenvalues of a single random graph are asymptotically biased.
When $n\rho_n \gg n^{1/2}$, \citet{zheng_limit_2022} consider the explicit form of the bias term $H^{(l)}$. 
\begin{theorem}[Restatement of Theorem~3.3 in \citet{zheng_limit_2022}]
    Instate the conditions of \cref{thm:MASE1}, and let $D^{(l)}$ denote the ${n+1 \choose 2} \times {n+1 \choose 2}$ diagonal matrix with diagonal elements of the form $P^{(l)}_{i j}(1 - P^{(l)}_{i j})$ for $i \leq j$. Define
    \begin{align*}
        \tilde D_{i i}^{(l)}
        &\coloneqq
        \sum_{j=1}^{n} P_{i j}^{(l)}(1 - P_{i j}^{(l)}).
    \end{align*}
    Let $\mu^{(l)}$ denote the vector
    \begin{align*}
        \mu^{(l)}
        &=
        \operatorname{vech}\bigg( \frac{1}{L} U\t \tilde D^{(l)} U ( S^{(l)\top})\inv -
        \frac{1}{2L^2} \sum_{l'=1}^{L} S^{(l)} (S^{(l')})\inv U\t \tilde D^{(l')} U(S^{(l')\top})\inv \bigg) \\
        &\quad+
        \operatorname{vech}\bigg( \frac{1}{L}( S^{(l)\top})\inv U\t \tilde D^{(l)} U -
        \frac{1}{2L^2} \sum_{l'=1}^{L} (S^{(l')\top})\inv U\t \tilde D^{(l')} U(S^{(l')})\inv S^{(l)} \bigg),
    \end{align*}
    where $\operatorname{vech}$ is the half vectorization operation that extracts the upper triangle of a matrix, represented as a column vector. There exists a sequence of orthogonal transformations $W$ such that, as $n \to \infty$, 
    \begin{align*}
        \big( \Sigma^{(l)} \big)^{\dagger/2} &\bigg( \operatorname{vech}\big( W \hat S^{(l)} W\t - S^{(l)}\big) - \mu^{(l)} \bigg)
        \dto
        \mathcal{N}(0,\Pi^{(l)}),
    \end{align*}
    where $(\cdot)^{\dagger}$ denotes the Moore--Penrose pseudoinverse and $\Pi^{(l)}$ denotes a symmetric idempotent matrix. If $\Sigma^{(l)}$ is invertible, then $\Pi^{(l)} = I$.
\end{theorem}
In other words, the stochastic bias term $H^{(l)}$ can be replaced with a deterministic bias term $\mu^{(l)}$ when $n\rho_n \gg n^{1/2}$. The additional condition $n\rho_n \gg n^{1/2}$ is used to establish that the random bias term concentrates sufficiently about its expectation $\mu^{(l)}$. A similar requirement appears in \cref{thm:Bspectral}. Informally, this condition is required in order to demonstrate that residual terms of order $(n\rho_n)^{-1}$ are of smaller order than $n^{-1/2}$, the latter of which is the rate of concentration about the expectation and the parametric rate of convergence to asymptotic normality. The authors also uses this result to obtain a Hotelling-type $T^2$ test statistic for testing if $R^{(l)} = R^{(l')}$, namely, if $B^{(l)} = B^{(l')}$.

The authors in \citet{zheng_limit_2022} also demonstrate the following result for the estimates $\hat U_{\operatorname{MASE}}$ obtained via the eigenvectors of $\frac{1}{L} \sum_{l=1}^{L} \hat U^{(l)} \hat U^{(l)\t}$.

\begin{theorem}[Restatement of Theorem~3.2 in \citet{zheng_limit_2022}] \label{thm:MASE2}
    Instate the conditions of \cref{thm:MASE1}, and suppose that $n\rho_n \gg \log^2(n)$ with $L = O(1)$. There exists a sequence of orthogonal matrices $W_*$ such that, as $n \to \infty$,
    \begin{align*}
        \bigg(L^{1/2} n \rho_n^{1/2} \big( (\hat U_{\operatorname{MASE}}W_*)_{i*} - U_{i*} \big) \bigg)
        \mid
        \tau(i) = k 
       & \dto
        \mathcal{N}(0,\mathbf{\Sigma}_{\operatorname{MASE}}(k))
    \end{align*}
    with asymptotic covariance matrix given by
    \begin{align*}
        \mathbf{\Sigma}_{\operatorname{MASE}}(k)
        &\coloneqq
        \frac{n^2 \rho_n}{L} \sum_{l=1}^{L} \big( S^{(l)} \big)\inv U\t \Xi^{(l,i)} U\big(S^{(l)}\big)\inv,
    \end{align*}
    where $\Xi^{(l,i)}$ is the diagonal matrix with diagonal elements $\Xi^{(l,i)}_{kk} = P_{ik}^{(l)}(1 - P_{ik}^{(l)})$.
\end{theorem}
The matrices $S^{(l)}$, $U$, and $\Xi^{(l,i)}$ change with $n$, but, due to the regularity conditions imposed in \citet{zheng_limit_2022}, the matrix $\mathbf{\Sigma}_{\operatorname{MASE}}(k)$ is well-defined. It is also possible to show that
\begin{align*}
    \mathbf{\Sigma}_{\operatorname{MASE}}(k)
    &=
    \frac{1}{L} \sum_{l=1}^{L} \mathbf{\Sigma}_U(k),
\end{align*}
where $\mathbf{\Sigma}_U(k)$ is the {(asymptotic)} covariance matrix appearing in \cref{thm:xieCLT}. To prove \cref{thm:MASE2}, \citet{zheng_limit_2022} demonstrates that $\hat U_{\operatorname{MASE}}$ behaves as an empirical average of the first-order perturbations of each of the individual networks. Namely, with high probability,
\begin{align*}
   \hat U_{\operatorname{MASE}} W - U
   &=
   \frac{1}{L} \sum_{l=1}^{L} ( A - P) U^{(l)} (\Lambda^{(l)})\inv + R_{\operatorname{MASE}}\\
   &=
   \frac{1}{L} \sum_{l=1}^{L} ( A - P) U (S^{(l)})\inv  + R_{\operatorname{MASE}},
\end{align*}
where $R_{\operatorname{MASE}}$ is a residual term similar to the residual terms in \cref{asympexp}.

The authors in \citet{jones_multilayer_2021} consider a multiple random graph model that subsumes the multilayer stochastic blockmodel by allowing for general edge probabilities of the form $P^{(l)} = X \Lambda^{(l)} X\t$, where $\Lambda^{(l)}$ is a general symmetric matrix for the $l$-th graph. For SBMs, one parameterization is $X_i = Z_{i*}$ and $\Lambda^{(l)} = B^{(l)}$, where $B^{(l)}$ is the block probability matrix for the $l$-th graph. Therein, the so-called unfolded adjacency matrix $A_{\operatorname{unfold}}$ is given by
\begin{align*}
    A_{\operatorname{unfold}}
    \coloneqq
    [A^{(1)}, A^{(2)}, \dots, A^{(L)}] \in \{0,1\}^{n \times nL},
\end{align*}
with the goal of estimating the scaled left and right singular vectors of the matrix $\mathbf{P}_{\operatorname{unfold}}$, defined analogously.

\begin{theorem}[Restatement of Theorem~3 in \citet{jones_multilayer_2021} for Multilayer SBMs] \label{thm:jonesCLT}
    Fix $K$ and $L$. Suppose that $n\rho_n \gg \log^c(n) n^{1/2}$ for some universal constant $c>0$, and suppose $\p[\tau(i) = k] = \pi_k$. Let $A^{(l)} \sim P^{(l)}$ with $P^{(l)}= \rho_n Z B^{(l)} Z\t$, where each $B^{(l)} = \nu^{(l)} \ipq (\nu^{(l)})\t$. Let $\hat U_{\operatorname{unfold}} \hat \Lambda_{\operatorname{unfold}} \hat V_{\operatorname{unfold}}\t$ be the top $r$ truncated singular value decomposition of $A_{\operatorname{unfold}}$ above. There exists a sequence of invertible linear transformations $T$ depending on $n$ such that, as $n \to \infty$,
    \begin{align*}
      \bigg(n^{1/2} \big( (\hat U_{\operatorname{unfold}} \hat \Lambda_{\operatorname{unfold}}^{1/2} T)_{i*} - e_k \big) \bigg)
      \mid
      \tau(i) = k
     &\dto
      \mathcal{N}(0,\mathbf{\Sigma}_{\operatorname{unfold}}(k)),
    \end{align*}
    where the asymptotic covariance matrix is given by
\begin{align*}
    &\mathbf{\Sigma}_{\operatorname{unfold}}(k)
    \coloneqq \begin{cases}
        \Delta_{\mathbf{B}}\inv \mathbf{B} ( \Gamma^{(k,1)} \oplus \cdots \oplus \Gamma^{(k,L)} ) \mathbf{B}\t  \Delta_{\mathbf{B}}\inv
        &\textnormal{~if~}\rho_n \equiv 1, \\ 
        \Delta_{\mathbf{B}}\inv \mathbf{B} ( \Gamma^{(k,1)}_{o(1)} \oplus  \cdots \oplus \Gamma^{(k,L)}_{o(1)} ) \mathbf{B}\t \Delta_{\mathbf{B}}\inv
        &\textnormal{~if~}\rho_n \to 0,
    \end{cases}
    \end{align*}
    where
    \begin{align*}
        \mathbf{B}
        &\equiv
        B_{\operatorname{unfold}}
        \coloneqq
        [B^{(1)},\dots, B^{(L)}],
        \\
        \Delta_{\mathbf{B}}
        &\coloneqq
        \mathbf{B} ( \operatorname{diag}(\pi) \otimes \mathbf{I}_L ) \mathbf{B}\t,\\
        \Gamma^{(k,l)}
        &\coloneqq
        \sum_{r=1}^K \pi_r B^{(l)}_{kr} (1 - B^{(l)}_{kr}) \nu^{(l)} (\nu^{(l)})\t, \\
        \Gamma^{(k,l)}_{o(1)}
        &\coloneqq
        \sum_{r=1}^K \pi_r B^{(l)}_{kr} \nu^{(l)} (\nu^{(l)})\t,
    \end{align*}
    where $\oplus$ denotes the direct sum operator, and $\otimes$ denotes the Kronecker product.
\end{theorem}
The paper \citet{jones_multilayer_2021} also establishes similar results for the scaled right singular vectors, and the results therein hold more generally, including for directed, bipartite, and different dimension networks, and for different amounts of individual network sparsity.

\cref{thm:jonesCLT} requires that $n\rho_n \gg \log^{c}(n)n^{1/2}$, which is a stronger requirement than the minimal requirement $n\rho_n \gtrsim \log(n)$ for perfect clustering when specialized to the single graph setting. The appearance of the matrices $\Delta_{\mathbf{B}}\inv$ and $\mathbf{B}$ in the definition of $\mathbf{\Sigma}_{\operatorname{unfold}}(k)$ may informally be viewed as a matrix-weighted average of the individual network covariances. In fact, under the condition that each matrix $B^{(l)}$ is the same, the covariance in \cref{thm:jonesCLT} matches the corresponding variance in the single-network setting (e.g.,~\cref{thm:prdCLT}), though now scaled by a factor of $L^{-1}$, which coincides with the setting where adjacency spectral embedding is applied on the averaged graph.

Many existing asymptotic normality results require that the absolute value of the smallest nonzero eigenvalue of $P$ tends to infinity at a rate faster than $\log(n)$. For the $K$-block balanced homogeneous stochastic blockmodel with $n/K$ vertices in each community, its smallest eigenvalue is $n(a-b)/K$, and hence this requires $n(a-b)/K \gtrsim \log(n)$. In contrast, \citet{huang_spectral_2023} considers the less-studied setting $n(a-b)/K \ll 1$, wherein the $K$-means clustering error need not decrease to zero as $n \to \infty$. 

\begin{theorem}[Restatement of Proposition~3 in \citet{huang_spectral_2023}] \label{thm:huangclt}
    Fix $K$ and $L$. Suppose there are precisely $n/K$ vertices in each community. Let $A^{(l)} \sim P^{(l)}$, with $P^{(l)} = ZB^{(l)}Z\t$ for $1 \leq l \leq L$, where each $B^{(l)}$ is of the form 
    \begin{align*}
        B^{(l)}_{kk'}
        &=
        \begin{cases}
            a^{(l)}
            &
            \textnormal{~if~}k = k', \\
            b^{(l)}
            &
            \textnormal{~if~}k \neq k',
        \end{cases}
    \end{align*}
    where $a^{(l)}$ and $b^{(l)}$ are allowed to depend on $n$. For a given nonnegative vector $w$, define 
    \begin{align}
        \tau_n^{w}
        \coloneqq
        \frac{n \left( \sum_l w_l ( a^{(l)} - b^{(l)})\right)^2}{\sum_l w_l^2 \left( a^{(l)} (1- a^{(l)}) + (K-1) b^{(l)}(1-b^{(l)}) \right)}
        .
        \label{eq:taun}
    \end{align}
     Let $\tau_{\infty}^w \coloneqq \lim_{n\to\infty} \tau_n^{w}$, and suppose that $0 < \tau_{\infty}^w < \infty$. Further suppose that
      \begin{gather*}
            \sum_{l} w_l^2 \left(a^{(l)}(1- a^{(l)}) + (K-1)b^{(l)}(1 - b^{(l)}) \right)
            \gtrsim
            \frac{\log^4(n)}{n}, \\ 
            \left[ \sum_l w_l^2 a^{(l)}(1 - a^{(l)}) \right] \left[ \sum_l w_l^2 b^{(l)} (1 - b^{(l)}) \right]\inv
            \to
            1.
    \end{gather*}
    Let $\eta_1, \dots, \eta_K$ be $(K-1)$-dimensional vectors such that the matrix
    \begin{align*}
        \begin{pmatrix}
            \frac{1}{\sqrt{K}}
            &
            \eta_1\t \\
            \vdots
            &
            \vdots \\
            \frac{1}{\sqrt{K}}
            &
            \eta_K\t
        \end{pmatrix}
        \in \R^{K \times K}
    \end{align*}
    is orthogonal. Define the vectors $\mu_1, \dots, \mu_K$ as
    \begin{align*}
        \mu_k
        &\coloneqq
        \begin{pmatrix}
            1 \\
            \sqrt{\frac{K (\tau_{\infty}^w - K)}{\tau_{\infty}^2}} \eta_k 
        \end{pmatrix}
        \in \R^{K}.
    \end{align*}
    Let $\hat U^{(w)}$ denote the leading $K$ eigenvectors of $\sum_{l} w_l A^{(l)}$. There exists a sequence of orthogonal matrices $W$ such that, as $n \to \infty$,
    \begin{align*}
      \bigg(n^{1/2} \big( \hat U^{(w)} W \big)_{i*} \bigg)
      \mid
      \tau(i) = k 
      \dto
      \mathcal{N}(\mu_k,\mathbf{\Sigma}^w),
    \end{align*}
    where $\mathbf{\Sigma}^w$ is given by 
    \begin{align*}
      \mathbf{\Sigma}^w
      \coloneqq
      \frac{K}{\tau_{\infty}^w}
      \begin{pmatrix}
        0
        &
        0 \\
        0
        &
        I_{K-1}
    \end{pmatrix}.
    \end{align*}
\end{theorem}

The primary difference between the normality result in \cref{thm:huangclt} and the previous asymptotic normality for the unscaled eigenvectors presented in \cref{thm:xieCLT} is that the probability matrix asymptotically degenerates to that of an Erd\H{o}s--R\'{e}nyi random graph (here $a - b \to 0$). The paper \citet{huang_spectral_2023} states a general scenario where the technical conditions are that $a^{(l)} \gtrsim \frac{\log^4(n)}{n}$, and $\frac{n(a^{(l)} - b^{(l)})^2}{a^{(l)}} \asymp 1$, which reflects sufficiently dense networks with average expected degree of order at least $\log^4(n)$ but where the communities are nearly indistinguishable. On the other hand, for the setting $a^{(l)} - b^{(l)} \asymp a^{(l)} \asymp b^{(l)}$, the literature typically emphasizes sparse graphs with distinguishable communities. For the regime studied in \citet{huang_spectral_2023}, the optimal misclustering error need not vanish in the limit, so their methodology is based on optimizing over weight vectors $w$ that minimize the asymptotic Bayes-optimal misclustering error determined by this limiting distribution. Further, any two rows of $\hat U^{(w)}$ are asymptotically independent.

The regime considered in \citet{huang_spectral_2023} is notably different from those most commonly studied in the statistics literature. Several extensions are possible. For example, it remains open to determine whether similar results hold for other types of stochastic blockmodels or multiple stochastic blockmodels with different block probability matrices. The averaging procedure therein requires that $a^{(l)} > b^{(l)}$ for all $l$, so it is not clear how to extend their procedure to settings with $b^{(l)} > a^{(l)}$. More generally, it is an open problem to determine optimal clustering procedures in regimes where the cluster separation vanishes but the networks remain dense. See \citet{kadavankandy_asymptotic_2019} for related considerations in a similar model setting.

The order $n^{1/2}$ scaling in \cref{thm:huangclt} differs from the $n\rho_n^{1/2}$ scaling in \cref{thm:xieCLT}. Further, the limiting covariance matrix is degenerate. When $\tau_{\infty}^w < \infty$, the leading eigenvalues besides the first eigenvalue are {separated from the other eigenvalues but with smaller gaps between distinct eigenvalues,} since $n\rho_n \ll \log(n)$. Therefore, unlike the settings discussed in \cref{sec:asymp_norm}, the matrix $A$ now consists of three parts: the first is a rank-one matrix corresponding to the (large) leading eigenvalue, the second part is a rank $K -1$ matrix corresponding to the bottom (small) $K -1$ nonzero eigenvalues, and the third part is a mean-zero matrix. The rank-one component associated to the leading eigenvalue contributes to the degeneracy in the asymptotic covariance. Since the rank-one component is well-separated from the other eigenvalues, it exhibits qualitatively similar phenomena to the leading eigenvectors when all the eigenvalues are spiked, and hence it is plausible for the leading eigenvector to converge in probability to the all ones vector after scaling by $n^{1/2}$. The second part of the matrix contributes to the centering, namely scaling by $n^{1/2}$ stabilizes these expectations, though they are now slightly biased away from the population eigenvectors. The remaining term in the decomposition of $A$ is a mean-zero (inhomogeneous Wigner) matrix, and the final scaling by $n^{1/2}$ {is determined by the appropriate scaling from the literature on inhomogeneous Wigner matrices in random matrix theory}(e.g.,~\citet{bai_limiting_2012}). All components contribute to the asymptotic normality in \cref{thm:huangclt}.

In preceding discussion, the number of layers $L$ is assumed fixed as $n \rightarrow \infty$, and thus might not apply to the case that $L \asymp n$. 
Several works consider settings where all the network-level block probability matrices $B$ are positive semidefinite, implying that an averaging procedure will suffice \citep{bhattacharyya_spectral_2018,chen_global_2022,paul_consistent_2016,paul_spectral_2020}. In this setting, developing asymptotic normality theory for scaled and unscaled eigenvectors follows with only slight modifications to previous works as long as $Ln\rho_n \gg \log(n)$. For example, focusing on the unscaled eigenvectors, one may develop an analogous first-order expansion to \eqref{asympexp} via
\begin{align*}
    \hat U_{\operatorname{ave}} W - U_{\operatorname{G}}
    &=
    \sum_{l=1}^{L} (A^{(l)} - \mathbb{E} A^{(l)}) U_{\operatorname{ave}} \Lambda_{\operatorname{ave}}\inv + R_U, 
    \numberthis
    \label{asympexp2}
\end{align*}
where $\hat{U}_{\operatorname{ave}}$ are the leading eigenvectors of $\sum_{l} A^{(l)}$, and $U_{\operatorname{ave}}$ and $\Lambda_{\operatorname{ave}}$ are the leading eigenvectors and eigenvalues of $\sum_{l} P^{(l)}$, respectively. Here, $R_U$ is a residual matrix with comparatively small $\ell_{2,\infty}$ norm with high probability provided $Ln\rho_n \gg \log(n)$.

The statistical theory in the setting that each $B^{(l)}$ has a different number of positive and negative eigenvalues has primarily focused on consistency to date \citep{jing_community_2021,lei_bias-adjusted_2022}. For example, \citet{lei_bias-adjusted_2022} focuses on spectral clustering via the matrix $\sum_{l=1}^{L}\mathcal{H}\{ (A^{(l)})^2 \}$, where $\mathcal{H}(\cdot)$ is the \emph{hollowing operator} that sets the diagonal of a matrix to zero. The empirical eigenvectors of this matrix are shown to be consistent with respect to $\|\cdot\|_F$ when $\sqrt{L}n\rho_n \gg \sqrt{\log(n + L)}$. However, there do not appear to exist asymptotic normality results in the literature for this setting, seeing as the hollowing operation introduces non-trivial technicalities in the requisite perturbation analysis. Seemingly the closest such result is for asymmetric matrices from \citet{yan_inference_2024}, who use the algorithm \texttt{HeteroPCA} proposed in \citet{zhang_heteroskedastic_2022} and also analyzed in \citet{agterberg_entrywise_2022}. We omit a detailed discussion for brevity, but in essence the algorithm impute the diagonals starting from the initialization $\sum_{l=1}^{L}\mathcal{H}\big\{ (A^{(l)})^2 \big\}$, where $\mathcal{H}(\cdot)$ is the \emph{hollowing operator} that sets the diagonal of a matrix to zero. 

\begin{theorem}[Restatement of Theorem~5 in \citet{yan_inference_2024} for SBMs] \label{yantheorem}
    Fix $K$.
    Suppose that $n_k \asymp n_l$ for all $1 \le k,l \le K$, that $\sqrt{L}n\rho_n \gg \log(n)$, that $Ln\rho_n \gg \log(n)$, and that $L \leq C n$.
    Let $A^{(l)} \sim P^{(l)}$ consist of independent entries (directed networks), where $P^{(l)} = \rho_n Z B^{(l)} Z\t$, and where each $B^{(l)}$ satisfies $c \leq |\lambda_{\min}(B^{(l)})| \leq \lambda_{\max}(B^{(l)}) \leq C$, for some fixed constants $c, C > 0$.
    Let $\hat U_{\operatorname{HP}}$ denote the output of the \texttt{HeteroPCA} algorithm for $t_0$ iterations, with $t_0 \geq \log( C \sqrt{L} n \rho_n )$, and let $U_{\operatorname{HP}}$ denote the leading $K$ eigenvectors of $\sum_{l} (P^{(l)})^2$. 
    Finally, define $\mathcal{P} \coloneqq [ P^{(1)}, \cdots P^{(L)} ] \in \mathbb{R}^{n \times nL}$ and suppose that $\mathcal{P}$ has singular value decomposition $U_{{\operatorname{HP}}} \Lambda_{{\operatorname{HP}}} V_{{\operatorname{HP}}}\t$. Then, there exists an orthogonal matrix $W$ such that with probability at least $1 - O(\max(n,L)^{-10})$,
    \begin{align*}
        &\hat U_{{\operatorname{HP}}} W - U_{{\operatorname{HP}}}  =
        ( \mathcal{A} - \mathcal{P} ) V_{{\operatorname{HP}}} \Lambda_{{\operatorname{HP}}}\inv  +
        \mathcal{H}\{( \mathcal{A} - \mathcal{P})( \mathcal{A} - \mathcal{P})\t \} U_{{\operatorname{HP}}} \Lambda_{{\operatorname{HP}}}^{-2}
        +
        R
        ,
    \end{align*}
    where the matrix $R$ satisfies
    \begin{align*}
        \| R \|_{2,\infty}
        &\lesssim
        \frac{\log(n)}{\sqrt{L} n^2 \rho_n}
        +
        \frac{\sqrt{\log(n)}}{n\sqrt{Ln\rho_n}}.
    \end{align*}
\end{theorem}

The above presentation differs significantly from the original statement in \citet{yan_inference_2024} but can be established using the fact that the squared singular values of $\mathcal{P}$ are equivalent to the eigenvalues of $\sum_{l} (P^{(l)})^2$, and hence are of order $L n^2 \rho_n^2$, by taking $n_2$ therein to be $nL$, and by using the fact that the eigenvectors of $\sum_{l} (P^{(l))^2}$ are incoherent. Furthermore, the results in \citet{yan_inference_2024} are for a very general model with independent heteroskedastic noise, a special case of the setting considered herein.
For simplicity, we have focused on the setting that $L \leq C n$, but a similar result can be stated even when $L \gg n$ provided one replaces $\log(n)$ with $\log(L)$ and one modifies the bound on the residual term appropriately. In particular, \cref{yantheorem} holds with slight modifications even when $L \to \infty$. However, unlike the asymptotic expansion in \eqref{asympexp}, \cref{yantheorem} contains two terms, the first of which is linear in the noise $\mathcal{A} - \mathcal{P}$, and the second of which is quadratic in this difference. Hence, any subsequent normality result will depend on both the first-order and second-order fluctuations.

\section{Discussion}
This survey illustrates how properties of the underlying block probability matrix $B$ influence the emergence and behavior of asymptotic normality phenomena. For not-too-sparse networks, this survey highlights that the row vectors of data-derived embeddings are asymptotically Gaussian about their population counterparts after proper orientation, with asymptotic covariance affected by the limiting behavior of the sparsity term $\rho_n$. Notably, cluster geometry and point cloud shape in low-dimensional Euclidean space are driven by characteristics of underlying network connectivity. Networks with communities that are difficult to distinguish on the basis of connectivity patterns yield highly elliptical cluster geometry compared to settings where communities are well-separated with sufficiently distinct connectivity properties.

Outside of cluster geometry, aggregating information from spectral embeddings also yields asymptotic normality. Furthermore, spectral embeddings and related quantities are useful for hypothesis testing and for parameter estimation. By understanding the geometry of spectral embeddings, one can devise test statistics with known limiting distributions or asymptotically efficient estimators of population parameters. Information-geometric properties also interface with latent space geometry, in the sense that spectral estimators of population parameters, without refinement, can be statistically inefficiency in rank-degenerate stochastic blockmodels. For multiple networks, the particular form of the asymptotic distribution can also depend on the relationships between different block probability matrices. Latent space geometry can substantially change based on the choice of spectral embedding and data aggregation, manifest in the diverse expressions for resulting asymptotic covariance matrices. For example, contrast \cref{thm:jonesCLT,thm:MASE2}.

Numerous open problems remain in the study of stochastic blockmodels. Below, we mention several pertinent topics and areas of ongoing investigation.
\begin{itemize}
    \item \textbf{Vanishing or weak signals}. Aside from \citet{huang_spectral_2023,kadavankandy_asymptotic_2019}, comparatively few works in the statistics literature address the asymptotic distributional properties of spectral embeddings when the signal strength in the model is either weak or vanishing. Such settings require the further treatment of regularization and normalization when constructing matrix representations of data and specifying the embedding method. More remains to be done.
    
    \item \textbf{Asymptotic normality of other embeddings}. This survey focuses on embeddings based on eigenvectors and eigenvalues of adjacency or Laplacian matrices. Other graph embeddings are also of interest, both theoretically and in applications, such as the eigenvectors of the \emph{regularized Laplacian} matrix \citep{qin_regularized_2013}, or the \texttt{Node2Vec} embedding procedure \citep{grover_node2vec_2016}. Several recent works in this direction include  \citet{barot_community_2021,cape_varimax_2024,kloepfer_delving_2021,lin_ergodic_2023,zhang_theoretical_2024}.
    
    \item \textbf{Choosing between embedding methods}. Different choices of embedding yield different forms of asymptotic normality, depending on whether or not the eigenvectors are scaled by the eigenvalues, or whether one uses the adjacency matrix or Laplacian matrix. Several works have studied the pros and cons of different embedding choices in the presence of network structure (e.g.,~core-periphery; affinity) per the underlying $B$ matrix, including but not limited to \citet{cape_spectral_2019,tang_limit_2018,priebe_two-truths_2019}, though they are limited to asymptotic regimes wherein perfect node clustering is possible. Deeper theoretical guidance to inform practice is needed.
    
    \item \textbf{Multiple testing problems for networks}. Throughout the natural and physical sciences, analysts frequently encounter or construct multiple networks corresponding to different populations, for example, brain graphs for healthy subjects, high-risk subjects, and individuals diagnosed with a particular disease. Multi-sample testing, discrimination, and comparison are consequently of interest beyond the two-vertex or two-graph problems discussed in this survey. To date, ongoing efforts in the study of multilayer or multiplex networks are paving the way for a deeper statistical understanding and applicable statistical toolkit for multiple graph inference \citep{fan_alma_2022,macdonald_latent_2022,pensky_clustering_2024,pensky_dynamic_2019,pensky_spectral_2019,zheng_limit_2022,jones_multilayer_2021,jing_community_2021,lei_bias-adjusted_2022,lei_consistent_2020,chen_hypothesis_2023,lyu_latent_2023}.
    
    \item \textbf{Trade-offs between sparsity and the number of networks}. Several works study community detection with a large or growing number of sparse networks \citep{bhattacharyya_spectral_2018,chen_global_2022,jing_community_2021,lei_consistent_2020,lei_bias-adjusted_2022,lyu_latent_2023,paul_spectral_2020}. Beyond consistency, it is of interest to study joint network embeddings and understand how resulting cluster geometry reflects network sparsity and sample size.
    
    \item \textbf{Covariates, dependencies, and model generalizations}. In various applied domains, networks often have vertex covariates or features, and edges need not be independent \citet{abbe_ell_p_2022,ma_universal_2020,ma_community_2023,mele_spectral_2022,mu_spectral_2022,razaee_matched_2019,xu_covariate-assisted_2022,zhang_node_2019}. It will be interesting to further study models and methods that incorporate and leverage such additional information when available.
\end{itemize}

\section*{Acknowledgements}
JA was partially supported by a fellowship from the Johns Hopkins Mathematical Institute of Data Science (MINDS) via NSF TRIPODS award CCF-1934979, by the Charles and Catherine Counselman Fellowship, and by the Acheson J. Duncan Fund for the Advancement of Research in Statistics. JC was partially supported by the National Science Foundation under grants DMS-1902755, SES-1951005, and DMS-2413552. JC gratefully acknowledges support from the University of Wisconsin–Madison, Office of the Vice Chancellor for Research and Graduate Education, with funding from the Wisconsin Alumni Research Foundation.

\bibliography{sbm_normality.bib}
\end{document}